\documentclass[11pt]{amsart}
\usepackage{amsmath,amsthm,color,amssymb, amscd, amsfonts,url,multicol,comment, xcolor, fancybox}
\usepackage[notcite,notref,final]{showkeys}
\usepackage{pdflscape}

\usepackage[utf8]{inputenc}
\usepackage{imakeidx}

\usepackage{tikz}
 \usetikzlibrary{positioning}
\usetikzlibrary{decorations.pathreplacing}

\usetikzlibrary{calc} 

\usepackage{amsxtra}
\usepackage{amssymb}
\usepackage{mathtools}
\usepackage{mathrsfs}  
\theoremstyle{plain}

\newcommand {\sectionnew}[1]{\section{#1}}

\newtheorem{theorem}{Theorem}[section]
\newtheorem{lemma}[theorem]{Lemma}
\newtheorem{sublemma}[theorem]{Sublemma}
\newtheorem{definition-theorem}[theorem]{Definition-Theorem}
\newtheorem{definition-lemma}[theorem]{Definition-Lemma}
\newtheorem{proposition}[theorem]{Proposition}
\newtheorem{corollary}[theorem]{Corollary}
\newtheorem{conjecture}[theorem]{Conjecture}
 \theoremstyle{definition}
\newtheorem{definition}[theorem]{Definition}
\newtheorem{example}[theorem]{Example}
\newtheorem{remark}[theorem]{Remark}
\newtheorem{notation}[theorem]{Notation}
\newtheorem{routine}[theorem]{Routine}
\newtheorem{problem}[theorem]{Problem}
\newtheorem{hypothesis}[theorem]{Hypothesis}
\newtheorem{question}[theorem]{Question}
\newcommand \bth[1] { \begin{theorem}\label{t#1} }
\newcommand \ble[1] { \begin{lemma}\label{l#1} }
\newcommand \bsubl[1] { \begin{sublemma}\label{sl#1} }

\newcommand \bdl[1] { \begin{definition-lemma}\label{dl#1} }
\newcommand \bpr[1] { \begin{proposition}\label{p#1} }
\newcommand \bco[1] { \begin{corollary}\label{c#1} }
\newcommand \bde[1] { \begin{definition}\label{d#1}\rm }
\newcommand \bex[1] { \begin{example}\label{e#1}\rm }
\newcommand \bre[1] { \begin{remark}\label{r#1}\rm }
\newcommand \bcj[1] { \begin{conjecture}\label{j#1}\rm }
\newcommand \bqu[1]  { \begin{question}\label{n#1}\rm }
\newcommand \bnota[1] { \begin{notation}\label{n#1}\rm }
\newcommand \bro[1] { \begin{routine}\label{n#1}\rm }
\newcommand \bpb[1] { \begin{problem}\label{n#1}\rm }
\newcommand \bhy[1] { \begin{hypothesis}\label{n#1}\rm }
\renewcommand {\eth} { \end{theorem} }
\newcommand {\ele} { \end{lemma} }
\newcommand {\esubl} { \end{sublemma} }

\newcommand {\edl}{ \end{definition-lemma} }
\newcommand {\epr} { \end{proposition} }
\newcommand {\eco} { \end{corollary} }
\newcommand {\ede} { \end{definition} }
\newcommand {\eex} { \end{example} }
\newcommand {\ere} { \end{remark} }
\newcommand {\ecj} { \end{conjecture} }
\newcommand {\equ} {\end{question}}
\newcommand {\enota} { \end{notation} }
\newcommand {\ero} { \end{routine} }
\newcommand {\epb} { \end{problem} }
\newcommand {\ehy} { \end{hypothesis} }
\newcommand \thref[1]{Theorem \ref{t#1}}
\newcommand \dlref[1]{Definition-Lemma \ref{dl#1}}
\newcommand \leref[1]{Lemma \ref{l#1}}
\newcommand \prref[1]{Proposition \ref{p#1}}
\newcommand \coref[1]{Corollary \ref{c#1}}

\newcommand \deref[1]{Definition \ref{d#1}}

\newcommand \lb[1]{\label{#1}}

\makeatletter              
\let\c@equation\c@theorem  
\makeatother
\numberwithin{equation}{section}

\def \Cset {{\mathbb C}}

\def \Zset {{\mathbb Z}}
\def \Nset {{\mathbb N}}

\def \Pset {{\mathbb P}}

\def \PP {{\mathcal{P}}}

\def \hb {{\hbar}}

\def \al {\alpha}
\def \be {\beta}

\def \om {\omega}

\def \ga {\gamma}

\def \sig {\sigma}

\def \sig{\sigma}

\def \mt  {\mapsto}

\def \hra {\hookrightarrow}

\def \rcor {\rangle}
\def \lcor {\langle}
\def \del {\partial}

\def \ol {\overline}
\def \wt {\widetilde}
\def \wh {\widehat}

\def \id { {\mathrm{id}} }

\def \nn  {\mathfrak{n}}
\def \mm  {\mathfrak{m}}
\def \b  {\mathfrak{b}}

\newcommand{\kk}{\Bbbk}

\DeclareMathOperator \Aut { {\mathrm{Aut}} }
\DeclareMathOperator \Der { {\mathrm{Der}} }

\DeclareMathOperator \maxSpec { {\mathrm{maxSpec}}}

\DeclareMathOperator \ad { {\mathrm{ad}} }

\DeclareMathOperator \tr { {\mathrm{tr}} }

\DeclareMathOperator \Irr { {\mathrm{Irr}} }

\newcommand \red {\textcolor{red}}
\newcommand \blue {\textcolor{blue}}

\newcommand{\pcoor}[1]{%
  \begingroup\lccode`~=`: \lowercase{\endgroup
  \edef~}{\mathbin{\mathchar\the\mathcode`:}\nobreak}%
  [
  \begingroup
  \mathcode`:=\string"8000
  #1%
  \endgroup
  ]
}

\input xy
\xyoption{all}

\makeindex[columns=3,
title=Index of Notation]
\begin{document}

\title[Poisson geometry and representations of   Sklyanin algebras]
{Poisson geometry and representations\\  of
PI 4-dimensional  Sklyanin algebras}
\author[Chelsea Walton]{Chelsea Walton}
\address{
Walton: Department of Mathematics \\
Rice University \\
Houston, TX 77005, 
USA
}
\email{notlaw@rice.edu}
\author[Xingting Wang]{Xingting Wang}
\address{
Wang: Department of Mathematics \\
Howard University \\
Washington, DC 20059, 
USA
}
\email{xingting.wang@howard.edu }
\author[Milen Yakimov]{Milen Yakimov}
\address{
Yakimov: Department of Mathematics \\
Northeastern University \\
Boston, MA 02115,
USA
}
\email{m.yakimov@northeastern.edu }
\thanks{Walton and Yakimov are partially supported by NSF grants \#1663775, 1601862, 1901830, 1903192, 2100756, and 2131243. 
Walton is also supported by the Sloan Foundation. Wang is partially supported by Simon collaboration grant \#688403}

\begin{abstract} 
Take $S$ to be a 4-dimensional Sklyanin (elliptic) algebra that is module-finite over its center $Z$; thus, $S$ is PI. Our first result is the construction of a Poisson $Z$-order structure on
$S$  such that the induced Poisson bracket on $Z$ is non-vanishing. We also provide the explicit Jacobian structure of this bracket, 
leading to a description of the symplectic core decomposition of the maximal spectrum $Y$ of $Z$.
We then classify the irreducible representations of $S$ by combining (1) the geometry of the Poisson order structures,
with (2) algebro-geometric methods for the elliptic curve attached to $S$, along with (3) representation-theoretic methods using line and fat point modules of $S$. Along the way, we improve
results of Smith and Tate obtaining a description the singular locus of $Y$ for such $S$.
The classification results for irreducible representations are in turn used to 
determine the zero sets of the discriminants ideals of these algebras~$S$.
\end{abstract}

\maketitle

\bibliographystyle{abbrv}  

\setcounter{tocdepth}{2} \tableofcontents


\vspace{-.4in}

\sectionnew{Introduction}
\lb{Intro}
\subsection{Overview of Sklyanin algebras and results in the paper}
\lb{1.1}

Sklyanin algebras are quadratic algebras that play a major role in the Artin-Schelter-Tate-van den Bergh's
classification of noncommutative projective spaces \cite{AST,AS,ATV1,ATV2} and the Feigin-Odeskii's investigation of elliptic algebras \cite{FO, FO2}. These directions were motivated by Sklyanin's work \cite{Skly1983} on 
quantum integrable systems, in which he introduced the algebras now known as {\it 4-dimensional Sklyanin algebras} and studied their representations with applications to the quantum inverse problem method in quantum and statistical mechanics. Since then, $n$-dimensional versions of Sklyanin's algebras were introduced for $n\geq 3$ and have arisen in numerous areas including: deformation-quantization, non-commutative geometry, and quantum groups and R-matrices; see, for instance, the reviews \cite{Rogalski, Odesskii, StaffordvdB} and the references within.

\smallbreak In terms of their representation theory, Sklyanin algebras fall into one of two classes-- those that are module-finite over their center so that they satisfy a polynomial identity (or, are PI), or otherwise. 

\smallbreak (A) In the generic (non-PI) case, a Sklyanin algebra is not finite-module over its center. Its representation theory can be viewed in parallel with the representations of quantum groups 
at non-roots of unity and Lie algebras over fields of characteristic 0. Sklyanin algebras are more complicated than the other two classes because they do not have PBW bases 
in the usual sense and thus one cannot use the classical methods of standard modules (for instance Verma modules in the case of the categories $\mathcal{O}$). 
However, Artin-Schelter-Tate-van den Bergh have developed powerful 
projective algebro-geometric methods that can be applied to analyze the representations of these algebras via algebraic geometry of elliptic curves. 
Such methods were further developed in works of Levasseur-Smith \cite{LevSmith}, Smith-Staniszkis \cite{SmithStan}, and Staniszkis \cite{Stan},
to achieve results on the noncommutative projective algebraic geometry of irreducible representations of generic 4-dimensional Sklyanin algebras.

 \smallbreak (B) Representations of PI Sklyanin algebras can be viewed in parallel with representations of quantum groups at roots of unity and modular representations 
of Lie algebras. This case is substantially harder than setting (A). {\em{In this paper we unify the algebro-geometric methods for Sklyanin algebras with 
Poisson geometric methods for quantum groups at roots of unity, to}}
\begin{enumerate}
\item {\em{construct nontrivial structures of Poisson orders on all PI 4-dimensional Sklyanin algebras, and}}
\item {\em{explicitly classify their irreducible representations and describe their dimensions.}}
\end{enumerate}

\medbreak Smith \cite{Smith1993} has obtained a number of important results in the direction (2) by exclusively using algebro-geometric methods. Pertaining to direction (1), De Concini-Kac-Procesi 
\cite{DKP1, DKP2} pioneered the applications of Poisson geometry in representation theory of PI algebras, for the cases of quantized universal 
enveloping algebras and quantum function algebras at roots of unity. This approach was axiomatized by Brown-Gordon \cite{BrownGordon} in the theory 
of {\it Poisson orders}, and was applied to other families of algebras with PBW bases, such as the symplectic reflection algebras. A Poisson order is
an algebra $A$ which is module-finite over its center $Z$, and such that $Z$ admits a Poisson structure for which all Hamiltonian 
derivations of $Z$ can be extended to derivations of $A$. The punchline of the construction is that Brown and Gordon used the latter extension  to 
construct an isomorphism $A/( \mm A) \cong A /(\nn A)$ for $\mm, \nn \in \rm{maxSpec} (Z)$ in the same symplectic core, \cite{BrownGordon,DKP1, DKP2}. 
This provides general Poisson geometric tools to organize the irreducible representations of $A$ into families that behave in similar ways (e.g., that have the same dimension).

 \smallbreak We use previous algebro-geometric results on the PI 4-dimensional Sklyanin algebras $S$, \cite{AST,AS,ATV1,ATV2, SmithStafford, SmithTate, Stafford},
to prove that every such algebra has a nontrivial structure of Poisson order. We classify the symplectic cores of $Z(S)$ and obtain from it a concrete description 
of the Azumaya locus of $S$. We then link the developed Poisson geometry back to the algebro-geometric approach to the representation theory of $S$, \cite{LeBruyn,LevSmith,Smith1993}, to classify the irreducible representations of $S$ and to describe their dimensions. See Sections~\ref{sec:1.2} and~\ref{sec:1.3} for more details.

\smallbreak Poisson structures on commutative algebras which arise as semiclassical limits of elliptic algebras were studied by Odesski\u\i-Fe{\u\i}gin \cite{FO}, Pym-Schedler \cite{PymSchedler},
and others. The novelty in our work is the noncommutative extension of Poisson algebras to Poisson orders for all PI 4-dimensional Sklyanin algebras $S$, 
which is the key feature that is needed to approach the representation theory of $S$. 

\smallbreak In our previous work \cite{WWY}, we have obtained similar results for the PI 3-dimensional Sklyanin algebras. The 4-dimensional case turns out to be substantially 
more challenging, and requires new approaches at various steps. We indicate those throughout the paper. On the other hand, we develop methods that can be applied
in far greater generality than 4-dimensional Sklyanin algebras, for instance, 
connected $\mathbb{N}$-graded algebras $T$ containing a regular sequence of central elements $\Omega_1, \dots, \Omega_m$ such that $T/(T \Omega_1  + \cdots + T \Omega_m)$ is isomorphic to a twisted homogeneous coordinate ring of an elliptic curve; such results are highlighted throughout our work.

\smallbreak 
We use the solution of the classification problem (2) for irreducible representations of the algebras $S$ to fully determine the zero sets of all of their discriminant ideals.
This will likely have further applications in relation to recent work of Bell-Zhang \cite{BellZhang}, Ceken-Palmieri-Wang-Zhang \cite{CPWZ}, and Brown-Yakimov \cite{BrownYakimov}.

\subsection{Poisson orders, classification of symplectic cores and Azumaya loci} \label{sec:1.2}
We work over an algebraically closed base field $\kk$ of characteristic~$0$. The projective algebro-geometric data attached to a 4-dimensional Sklyanin algebra $S$ are as follows: 
an elliptic curve $E$ in $\mathbb{P}^3$, an invertible sheaf $\mathcal{L}$ on $E$ and an automorphism of $\sigma$ of $E$ given by a translation of a point $\tau\in E$. See Section~\ref{4DSkly} for details. 

\smallbreak One useful fact we have is that $S$ is module-finite over its center, if and only if the automorphism $\sigma$ attached to $S$ has finite order.  In this case, $S$ is PI and the  PI degree of $S$ (which is a sharp upper bound on the dimension of irreducible representations of $S$) is equal to $|\sigma|$; see Proposition~\ref{pPI}. So, consider the following standing hypothesis and notation used throughout this work.

\bhy{hyp}[$S$, $Z$, $n$, $s$]  Let  $S$ \index{S0@$S$} be a 4-dimensional Sklyanin algebra, and (with the exception of Sections~\ref{sec:back} and~\ref{sec:formal}) assume that $S$ is module-finite over its center $Z$ \index{Z0@$Z$} so that $S$ is PI. For the geometric data $(E, \mathcal{L}, \sigma)$ attached to $S$ mentioned above, let $n$ \index{n5@$n$} denote $|\sigma|$, which is equal to the PI degree of $S$ when $|\sigma|$ is finite. Also, let $s:=n/(n,2)$. 

In this work, we do not consider the cases when $n$ divides 4, as done typically in works on 4-dimensional Sklyanin algebras \cite{LevSmith, SmithStan, Stan}. (See, e.g., \cite[proof of Lemma~2.1]{Stan}.)
\ehy

Smith and Tate \cite{SmithTate} proved that the center $Z(S)$ of every PI 4-dimensional Sklyanin algebra $S$ is generated by four elements 
$z_0, \ldots, z_3$ of degree $n$ and two elements $g_1, g_2$ of degree $2$, subject to two relations $F_1$ and $F_2$ of degree $2n$. 
Set
\[
Y \mbox{ \index{Y@$Y$} } := {\rm{maxSpec}}(Z(S)) \quad \mbox{and} 
\quad Y_{\gamma_1, \gamma_2}  \index{Y2gama@$Y_{\gamma_1, \gamma_2}$}  := Y \cap \mathbb{V}(g_1 - \gamma_1, g_2 - \gamma_2), \quad \text{ for } \gamma_1, \gamma_2 \in \kk. 
\]
Denote by $Y^{sing}$\index{Y1sing@$Y^{sing}$}  the singular locus of $Y$. Further, $(Y_{\gamma_1,\gamma_2})^{sing}$ \index{Y2gamasing@$(Y_{\gamma_1, \gamma_2})^{sing}$}
denotes the singular locus of the subvariety $Y_{\gamma_1,\gamma_2}$ and
$$
(Y^{sing})_{\gamma_1,\gamma_2}:=Y^{sing}\cap Y_{\gamma_1,\gamma_2}.
\index{Y3singgama@$(Y^{sing})_{\gamma_1,\gamma_2}$}
$$
Here, $Y^{sing} = \bigcup_{\gamma_1,\gamma_2 \in \kk} (Y^{sing})_{\gamma_1,\gamma_2}$, which is contained in but not necessarily equal to $\bigcup_{\gamma_1, \gamma_2 \in \kk} (Y_{\gamma_1, \gamma_2})^{sing}$.
Moreover, consider the following notation for $\gamma_1,\gamma_2\in \kk$: \index{Y1Pois@$Y^{symp}_0$} 
$$(Y_{\gamma_1, \gamma_2})^{smooth}:=Y_{\gamma_1, \gamma_2}\setminus (Y_{\gamma_1,\gamma_2})^{sing}, \quad \quad \textstyle Y^{symp}_0:=\bigcup_{\gamma_1,\gamma_2\in \kk} (Y_{\gamma_1,\gamma_2})^{sing}.$$

\smallbreak Our first main result constructs a nontrivial Poisson order on $S$,  classifies the symplectic cores of its center $Z(S)=\kk[z_0,z_1,z_2,z_3,g_1,g_2]/(F_1,F_2)$, and fully determines the Azumaya locus of $S$. 

\bth{thm1}\textnormal{[Theorems \ref{cps-Z(S)} and~\ref{tAzumaya}, Corollary~\ref{cPoissonZero}, Proposition~\ref{pscores}]} For all 4-dimensional Sklyanin algebras $S$ satisfying Hypothesis \ref{nhyp} the following hold:
\begin{enumerate}
\item $S$ admits a nontrivial structure of Poisson order for which 
\begin{enumerate}
\item the induced Poisson structure on the center $Z(S)$ is of Jacobian form in the sense of \eqref{Jacobian-Pbrack} 
in terms of two potentials taken to be $F_1$ and $F_2$, while 
\item $g_1,g_2$ lie in the Poisson center of $Z(S)$. 
\end{enumerate}

\smallskip

\item The corresponding symplectic core stratification of the Poisson variety $Y$ are 
\begin{enumerate}
\item 2-dimensional cores: $(Y_{\gamma_1, \gamma_2})^{smooth}$ for $\gamma_1,\gamma_2\in \kk$;
\item 0-dimensional cores: points in $Y^{symp}_0$; these are the points on which the Poisson bracket on $Y$ in part \textnormal{(1)} vanishes. 
\end{enumerate}

\smallskip

\item $Y^{sing}\subseteq Y^{symp}_0$, with strict containment when $n$ is even, and they both have codimension $\ge 2$ in $Y$;

\smallskip

\item The Azumaya locus of $S$ coincides with the smooth locus $Y\setminus Y^{sing}$ in $Y$. 
\end{enumerate}
\eth

The above result can be thought of as a generalization of the previous work done for PI 3-dimensional Sklyanin algebras. There are a number of complications that arise in the 4-dimensional case compared to the 3-dimensional one. The first one is that in the 4-dimensional case we have $Y^{sing}\subsetneq Y^{symp}_0$ when the order $n$ of $\sigma$ is odd; in other words the set of symplectic points $Y_0^{symp}$ of the Poisson variety $Y$ is larger than its singular locus.
Indeed  in that case $Y^{symp}_0$ is the union of $Y^{sing}$ plus four cuspidal curves meeting at the origin; 
this is proved in Theorem \ref{tSOdd}. 

\smallbreak See Figure~1 in Section~\ref{sec:sing} for an illustration of the $n$ odd case, where $Y^{sing} \subsetneq Y_0^{symp}$ as discussed above. On the other hand, Figure~2 in Section~\ref{sec:sing} illustrates the $n$ even case, when $Y^{sing}$ and $Y_0^{symp}$ coincide.

\smallbreak 
As mentioned in Section~\ref{1.1}, the proof of the existence of nontrivial Poisson order structure on $S$ relies on the previous algebro-geometric results for $S$
and a method of higher order specialization developed in \cite{WWY}.
The precise (Jacobian) form of the Poisson structure on $Z(S)$ is derived from the property that $g_1, g_2$ are in the Poisson center of $Z(S)$ and from
Stafford's result \cite{Stafford} that $S$ is a maximal order.

\smallbreak The proof of \thref{thm1}(4) faces a second major difficulty compared to the 3-dimensional case. In the latter case, the symplectic cores of $Z(S)$ of maximal dimension 
have codimension 1, which was used to show that all of them are in the Azumaya locus of $S$ by incorporating the automorphism group of $S$. This strategy fails in the 
4-dimensional case since the leaves have codimension 2, while the automorphism group of $S$ is still 1-dimensional.
Our strategy for the proof of part (4) is summarized as follows. Firstly, we consider the factors
$$
S_{[\kappa_1:\kappa_2]} := S/ (\kappa_1 g_1 - \kappa_2 g_2) S, \quad \quad \text{for } [\kappa_1:\kappa_2] \in \mathbb{P}^1_\kk
\index{S4kappa@$S_{[\kappa_1:\kappa_2]}$}
$$
which all turn out to be Noetherian PI domains of the same PI degree as $S$. We show that the maximal spectra of their centers are canonically isomorphic to Poisson 
subvarieties of $Y$ with the property that the symplectic leaves in them of maximal dimension have codimension 1 in $Y$. Now, by
applying Brown and Gordon's Poisson order result \cite{BrownGordon} 
(as discussed in Section~\ref{sec:BGor-thm}) to the symplectic core stratification of the centers of the algebras $S_{[\kappa_1:\kappa_2]}$ and using the density of their Azumaya loci, 
we conclude that all the points in the 
2-dimensional cores $(Y_{\gamma_1, \gamma_2})^{smooth}$ lie in the Azumaya locus of $S$. Then the coincidence of the Azumaya and smooth locus of $S$ follows from a general result of Brown and Goodearl~\cite{BrownGoodearl} since the non-Azumaya locus of $S$ contained in $Y^{symp}_0$ has codimension $\ge 2$ in $Y$.
Thus, Poisson geometry was employed primarily to establish the codimension $\geq 2$ fact, for which there are no other known methods.

\smallbreak The novelty of the 4-dimensional Sklyanin algebras of odd PI degree is that their Azumaya loci intersects nontrivially the varieties of symplectic points $Y_0^{symp}$ of $Y$.
\subsection{On irreducible representations of PI 4-dimensional Sklyanin algebras} \label{sec:1.3}
Recall that $n$ denotes the order of the automorphism $\sigma$ of the elliptic curve attached to $S$;
it equals the PI degree of $S$. For $y \in Y$, denote by $\mm_y$ \index{m1y@$\mm_y$} the corresponding maximal ideal of $Z(S)$. It follows from 
\thref{thm1}(4) that 
\begin{itemize}
\item For $y \in Y \backslash Y^{sing}$, $S/\mm_y S \cong M_n(\kk)$, and thus, up to isomorphism, 
there is precisely one irreducible representation of $S$ with central annihilator $\mm_y$ 
and this representation has dimension $n$.
\end{itemize} 
It is a simple fact that
\begin{itemize}
\item For $y = \underline{0}$, $S/\mm_y S$ is a local algebra, and thus, up to isomorphism, 
$S$ has precisely one irreducible representation with central annihilator $\mm_{\underline{0}}$ 
and this representation is the trivial representation of dimension $1$.
\end{itemize}

To achieve results on irreducible representations of intermediate dimension, we first use the explicit structure of the center $Z(S)$ and its link with the geometry of the elliptic curve $E$ to describe the singular locus $Y^{sing}$ of $S$.

\bth{thm2}\textnormal{[Theorems \ref{tSOdd}, \ref{tSEven}]} Let $S$ be a 4-dimensional Sklyanin algebra satisfying Hypothesis \ref{nhyp}. Denote by $E_2$ the subgroup of 2-torsion points of $E$.
\begin{enumerate}
\item If $n$ is odd, then the singular locus of $Y$ is the union of $2(n-1)$ cuspidal curves $C(\omega+k\tau)$,
$$Y^{sing}~=\bigcup_{\substack{\omega\in E_2\\0\le k\le n-2} } C(\omega+k\tau),$$
defined in Lemma \ref{lWcusp}, which only meet at the origin, as depicted in Figure~1.

\smallskip

\item If $n=2s$ is even, then $Y^{sing}$ is the union of two explicitly defined subvarieties $Y_1^{sing}$ and $Y_2^{sing}$ in \thref{SEven}, which only meet at the origin, as depicted in Figure~2. When $\kk=\mathbb C$, they contain the cuspidal curves $$\mathcal C(\omega + k \tau) \quad \quad \text{for } \omega \in E_2, ~~0 \leq k \leq s-1$$ described in Notation \ref{nscriptC}. 
\end{enumerate}
\eth

The proof of Theorem~\ref{tthm2} relies on the geometry of line modules and fat point modules of $S$, extending results of \cite{LevSmith, SmithTate}.

\smallbreak Next, we fully classify the irreducible representations of $S$ of intermediate dimension and describe their dimensions. Along the way,
using the action of Heisenberg group $H_4$ of order 64 on $S$, we determine the exact form of the Smith-Tate 
defining relations $F_1$ and $F_2$ of the center $Z(S)$ when $n$ is even, which is of independent interest (see Proposition~\ref{pZEven} below).

\bth{thm3}\textnormal{[Theorems \ref{tOddS}, \ref{teven}]} Let $S$ be a complex 4-dimensional Sklyanin algebra satisfying Hypothesis \ref{nhyp}. Denote by $\Irr_d S$ the isomorphism classes of all $d$-dimensional irreducible representations of $S$. Recall the notation above.

\begin{enumerate}
\item If PIdeg($S$) = $n$ is odd, then we have the following maps via central annihilators.
{\small \[
\hspace{-1.23in} \xymatrix{
\Irr_n S\ar@{->>}[r]^-{1:1} & Y^{smooth}= (Y\setminus Y^{sing})
}
\]
\[
\xymatrix{
\Irr_{k+1} S\sqcup \Irr_{n-1-k} S \ar@{->>}[r]^-{2:1} & (C(\omega+k\tau)=C(\omega+(n-2-k)\tau)) \setminus\{\underline{0}\}, \quad 0\le k\le n-2
}
\]
\[
\hspace{-2.4in} \xymatrix{
\Irr_1 S\ar@{->>}[r]^-{1:1} & \{\underline{0}\}
}
\]}
\smallskip
\item If PIdeg($S$) = $n=2s$ is even, then we have the following maps via central annihilators. 
\smallskip
{\small \[
\hspace{-1.5in} \xymatrix{
\Irr_n S\ar@{->>}[r]^-{1:1} & Y^{smooth}=(Y\setminus Y^{sing})
}
\]
\[
\hspace{-0.8in} \xymatrix{
\Irr_s S\ar@{->>}[r]^-{2:1} & Y^{sing}\setminus \bigcup_{\omega\in E_2,0\le k\le s-2}\, \mathcal C(\omega+k\tau)
}
\]
\[
\xymatrix{
\Irr_{k+1} S{ \sqcup }\Irr_{s-1-k} S\ar@{->>}[r]^-{4:1} & (\mathcal{C}(\omega+k\tau)=\mathcal{C}(\omega+s\tau+(s-2-k)\tau)) \setminus\{\underline{0}\},  \quad 0\le k\le s-2;
}
\]
\[
\hspace{-2.65in} \xymatrix{
\Irr_1 S\ar@{->>}[r]^-{1:1} & \{\underline{0}\}
}
\]}
\end{enumerate}
\eth

Over complex numbers $\kk = \mathbb{C}$, Sklyanin constructed in \cite{Skly1983} for each $\omega\in E_2$ and $k\in \mathbb N\cup \{0\}$, a representation $V(\omega+k\tau)$ \index{V0omega@$V(\omega+k\tau)$}
over $S$ in a certain $(k+1)$-dimensional subspace of theta functions of order $2k$; see \cite[Section 3]{SmithStan} for details. These representations were proved later by Smith and Staniszkis to be irreducible whenever $k<s$ \cite[Theorem~3.6]{SmithStan}. Our next result shows that every irreducible representation of $S$ of intermediate dimension $<s$ is a scalar twist of $V(\omega+k\tau)^\lambda$ for some $\lambda\in \kk$ such that $V(\omega+k\tau)^\lambda$ equals $V(\omega+k\tau)$ as vector spaces and $s_i\cdot_\lambda v=\lambda^is_iv$ for any homogeneous element $s_i$ in  $S$ of degree $i$ and $v\in V(\omega+k\tau)^\lambda$. 

\bth{thm4}\textnormal{[Theorems \ref{tOddS}, \ref{teven}]}
Let $S$ be a complex 4-dimensional Sklyanin algebra satisfying Hypothesis \ref{nhyp}. Then $S$ has irreducible representations of each dimension $1,2,3, \dots, n ~(=s)$ if $n$ is odd, and of each dimension $1,2,3,\dots,\frac{n}{2}~ (=s), n$ if $n$ is even. Moreover, the nontrivial irreducible representations of intermediate dimension $<s$ are given by scalar twists of $V(\omega+k\tau)$ for all $\omega\in E_2$ and $0\le k\le s-2$. 
\eth

To prove Theorems \ref{tthm3} and \ref{tthm4}, we use the classification of fat points of $S$ by Smith in \cite{Smith1993}, a number of other representation-theoretic results of Levasseur, Le Bruyn, Smith and Staniszkis \cite{LeBruyn,LevSmith,Smith1993,SmithStan}, and 
part (4) of \thref{thm1} that the non-Azumaya part of $Z(S)$ is $Y^{sing}$. 
In particular, we apply the deep connections between fat point modules, $\mathbb{C}^\times\times {\rm PGL}_d$-stabilizers of irreducible representations, and $\mathbb{C}^*$-families of irreducible representations of a graded algebra within the framework of the noncommutative projective algebraic geometry of $S$. 

\smallbreak Finally, \thref{disc} in Section~\ref{sec:discr} contains the aforementioned results on the description of the  discriminant ideals of the algebras $S$.

\bre{uqsl2}
Note that Theorem~\ref{tthm3} makes precise the representation-theoretic connection between the PI 4-dimensional Sklyanin algebras and the quantized enveloping algebra $U_q(SU(2))$ at $q$ a root of unity, as introduced initially in the physics literature. Namely, the main result of Roche-Arnaudon \cite{RocheArnaudon} is that when $q$ is a root of unity of order $m$, then $U_q(SU(2))$ has $d$-dimensional irreducible representations for all $1 \leq d \leq m$ when $m$ is odd, and has $d$-dimensional irreducible representations for all $1 \leq d \leq \frac{m}{2}$ and $d = m$ when $m$ is even. They also show that $U_q(SU(2))$ arises as a `trigonometric limit' of a 4-dimensional Sklyanin algebra $S$. Connections between $S$ and $U_q(SU(2))$ (or rather, $U_q(\mathfrak{sl}_2)$) are also discussed in work of Smith and Stafford \cite[Section~1]{SmithStafford}.
\ere

\medskip
\noindent {\bf Acknowledgements}.  The authors thank S. Paul Smith for generously making available his 1993 work \cite{Smith1993} via private communication and for also making it available for public use recently on the ArXiv e-print service. The authors also thank the anonymous referee for their numerous comments, which enabled us to greatly improve the quality of this manuscript.


\sectionnew{Preliminary results on 4-dimensional Sklyanin algebras $S$} \label{sec:back}

We provide in this section background material on (the noncommutative projective algebraic geometry of) 4-dimensional Sklyanin algebras. This includes a discussion of various ring-theoretic and homological properties of these algebras, which were established by Smith-Stafford \cite{SmithStafford}. We also provide a detailed analysis of the center and symmetries of 4-dimensional Sklyanin algebras, which will play a key role in the rest of the paper. This analysis extends results of Smith-Tate \cite{SmithTate}, Smith-Staniszkis \cite{SmithStan}, and Chirvasitu-Smith \cite{ChirSmith}.
\smallskip

\begin{center}
{\it In this section, we do not assume that the 4-dimensional Sklyanin algebras $S$  are PI.}
\end{center}


\subsection{Noncommutative geometry of $S$}
\label{4DSkly}

Here, we recall the definition and properties of the 4-dimensional Sklyanin algebras and the corresponding twisted homogeneous coordinate rings.

\bde{Skly4} [$S$, $S(\alpha, \beta, \gamma)$, $x_i$] \cite{SmithStafford}
\index{S0@$S$}
\index{S1albega@$S(\alpha,\beta,\gamma)$}
\index{alpha@$\alpha$} \index{beta@$\beta$} \index{g3amma@$\gamma$} 
Take $\alpha, \beta, \gamma \in \kk$ so that
\begin{align}
\label{albega-cond1}
&\alpha + \beta + \gamma + \alpha \beta \gamma =0 \quad \mbox{[$ \Leftrightarrow (1+\alpha)(1+\beta)(1+\gamma) = (1-\alpha)(1-\beta)(1-\gamma)$]}, \\
\label{albega-cond2}
&(\alpha, \beta, \gamma) \not \in \{(-1,1,\gamma),(\alpha,-1,1),(1,\beta,-1)\}.
\end{align}
Then the (regular) {\it 4-dimensional Sklyanin algebra} $S:=S(\alpha, \beta, \gamma)$ are $\kk$-algebras, generated by noncommuting variables 
$x_0$, $x_1$, $x_2$, $x_3$\index{x3s@$x_0, \ldots, x_3$}
of degree one, subject to the following relations  
\[
\begin{array}{llll}
x_0x_1-x_1x_0 ~=~  \alpha(x_2x_3+x_3x_2), &&& x_0x_1+x_1x_0 ~=~ x_2x_3-x_3x_2,\\
x_0x_2-x_2x_0 ~=~  \beta(x_3x_1+x_1x_3), &&& x_0x_2+x_2x_0 ~=~  x_3x_1-x_1x_3,\\
x_0x_3-x_3x_0 ~=~  \gamma(x_1x_2+x_2x_1), &&& x_0x_3+x_3x_0 ~=~  x_1x_2-x_2x_1.\\
\end{array}
\]
\ede

\smallbreak These algebras come equipped with geometric data that is used to establish many of their ring-theoretic, homological, and representation-theoretic properties. To start, take a $\mathbb{N}$-graded algebra $A = A_0 \oplus A_1 \oplus A_2 \oplus \cdots$ with $A_0 = \kk$, and recall that an {\it $A$-point module} is a cyclic, graded left $A$-module with Hilbert series $(1-t)^{-1}$. These modules play the role of points in noncommutative projective algebraic geometry and the parameterization of $A$-point modules is referred to as the {\it point scheme} of $A$. See \cite[Chapter~I, Section~3]{Rogalski} for more details.

\bdl{geomS}\textnormal{[$\widehat{E}$, $E$, $\phi_1$, $\phi_2$,  $v_i$, $e_i$, $\sigma$, $\tau$]} \cite[Section~2]{SmithStafford}  The point scheme of the 4-dimensional Sklyanin algebra $S=S(\alpha, \beta, \gamma)$ is given by the union $\widehat{E}$\index{E0@$E, \widehat{E}$} 
of an elliptic curve $E:= \mathbb{V}\left(\phi_1, \phi_2 \right) \subseteq \mathbb{P}^3_{\pcoor{v_0:v_1:v_2:v_3}}$,\index{E0@$E, \widehat{E}$}\index{vs@$v_0, \ldots, v_3$} 
where
$$
\phi_1 = v_0^2+v_1^2+v_2^2+v_3^2  \quad \text{ and } \quad \phi_2 =  \textstyle \frac{1-\gamma}{1+\alpha}~v_1^2 + \frac{1+\gamma}{1-\beta}~v_2^2 + v_3^2,
\index{p30hi@$\phi_1, \phi_2$}
$$
and the four points 
$$
\{e_0:=[1:0:0:0], ~~e_1:=[0:1:0:0], ~~e_2:=[0:0:1:0],~~ e_3:=[0:0:0:1]\}.
\index{e4@$e_0, e_1, e_2, e_3$}
$$

\smallbreak The automorphism $\sigma= \sigma_{\alpha \beta \gamma}$\index{s6igma@$\sigma$}
of $\widehat{E}$ attached to $S$ fixes each of the four points $e_i$, and on $E$ it is defined on a dense open subset by 
\begin{equation} \label{sigma}
\sigma: \left[ \begin{array}{l} v_0 \smallskip \\ v_1 \smallskip \\ v_2  \smallskip \\ v_3 \end{array}\right] \mapsto
\left[ \begin{array}{c} 
-2\alpha \beta \gamma v_1 v_2 v_3 - v_0(-v_0^2+\beta\gamma v_1^2 + \alpha \gamma v_2^2 + \alpha \beta v_3^2) \smallskip \\ 
2\alpha  v_0 v_2 v_3 + v_1(v_0^2 - \beta\gamma v_1^2 + \alpha \gamma v_2^2 + \alpha \beta v_3^2) \smallskip \\ 
2\beta  v_0 v_1 v_3 + v_2(v_0^2+\beta\gamma v_1^2 - \alpha \gamma v_2^2 + \alpha \beta v_3^2) \smallskip \\ 
2 \gamma v_0 v_1 v_2 + v_3(v_0^2+\beta\gamma v_1^2 + \alpha \gamma v_2^2 - \alpha \beta v_3^2) \end{array}\right]
\cdot
\end{equation}
The automorphism $\sigma$ of $E$ is given by translation by a point of $E$; call this point $\tau$.\index{t1au@$\tau$} 
The triple $\left(E,  ~\mathcal{O}_{\mathbb{P}^3}(1)|_E, ~\sigma\right)$ is referred to as  the geometric data of $S$.
\qed
\edl

Using this  data, we consider a noncommutative coordinate ring of $E$; its generators are sections of the invertible sheaf $\mathcal{O}_{\mathbb{P}^3}(1)|_E$ and its multiplication depends on the automorphism $\sigma$. The general construction is given as follows.

\bde{thcr}[$\mathcal{L}_i$] \cite{ArtinVdB} Given a projective scheme $X$, an invertible sheaf $\mathcal{L}$ on $X$, and an automorphism $\sigma$ of $X$,
the {\it twisted homogeneous coordinate ring} attached to this geometric data is a graded $\kk$-algebra 
$$B(X,  \mathcal{L}, \sigma) = \textstyle \bigoplus_{i \geq 0} B_i, \quad \text{where } B_i := H^0(X, \mathcal{L}_i)$$
 with $\mathcal{L}_0 = \mathcal{O}_X$, $\mathcal{L}_1 = \mathcal{L}$, and $\mathcal{L}_i \index{L1i@$\mathcal{L}_i$} = \mathcal{L} \otimes \mathcal{L}^{\sigma} \otimes \cdots \otimes \mathcal{L}^{\sigma^{i-1}}$ for $i \geq 2$. The multiplication map $B_i \otimes B_j \to B_{i+j}$ is defined by $b_i \otimes b_j \mapsto b_i b_j^{\sigma^i}$ using $\mathcal{L}_i \otimes \mathcal{L}_j^{\sigma^i} = \mathcal{L}_{i+j}$.
 \ede
 
\bnota{BELsig} [$B$, $\mathcal{L}$] Let $B$\index{B0@$B$} 
denote the twisted homogeneous coordinate ring attached to the geometric data $\left(E, \mathcal L, \sigma\right)$ from \dlref{geomS}, where 
$\mathcal L:=\mathcal{O}_{\mathbb{P}^3}(1)|_{E}$.\index{L1@$\mathcal{L}$}
 \enota

It is often useful to employ the following embedding of $B$ into a skew-Laurent extension of the function field of $E$.

\ble{embedB} \cite{AST}
Given $\left(E, \mathcal L,\sigma\right)$   from \dlref{geomS}, extend $\sigma$ to an automorphism of the field  $\kk(E)$ of rational functions on $E$ by $\nu^{\sigma}(p) = \nu(\sigma^{-1} p)$ for $\nu \in \kk(E)$ and $p \in E$. 
For any nonzero section $w$ of $\mathcal L$, that is, any degree 1 element of $B$, take $D$ to be the divisor of zeros of $w$, and let $V$ denote $H^0(E, \mathcal{O}_{E}(D)) \subset \kk(E)$. 

Then, the vector space isomorphism $\nu w \mapsto \nu t$ for $\nu \in V$ extends to an embedding of $B$ in $\kk(E)[t^{\pm 1}; \sigma]$. Here, $t\nu = \nu^{\sigma} t$ for $\nu \in \kk(E)$. \qed
\ele

Now the first step in obtaining  useful  properties of 4-dimensional Sklyanin algebras is to use the result below.
 
\ble{factor-g} \textnormal{[$g_1, g_2$]} \cite[Lemma~3.3, Corollary~3.9, Theorem~5.4]{SmithStafford} 
\index{g12@$g_1,g_2$}
The degree 1 spaces of $S$ and of $B$ are equal and there is a surjective map from $S$ to $B$, whose kernel is generated by the two central degree 2  elements below
\begin{equation}
\label{g-elem}
g_1:= -x_0^2+x_1^2+x_2^2+x_3^2 
\quad \quad \text{ and } \quad \quad
  g_2:=x_1^2 + \textstyle \frac{1+\alpha}{1-\beta}~x_2^2 +  \frac{1-\alpha}{1+\gamma}~x_3^2 \end{equation}
Moreover, $\{g_1, g_2\}$ is a central regular sequence in $S$.
\qed
\ele

Many good ring-theoretic and homological properties of $S$ are obtained by lifting such properties from the factor $B$,  some of which are listed in the following result.

\bpr{hilb}   \cite[Theorem~5.5]{SmithStafford} \cite[Corollary~6.7]{Levasseur}
The 4-dimensional Sklyanin algebras are Noetherian domains of global dimension 4, they satisfy the Artin-Schelter Gorenstein condition, along with the Auslander-regular and Cohen-Macaulay conditions, and they have Hilbert series $(1-t)^{-4}$. \qed
\epr

\smallbreak The representation theory of both $S$ and $B$ depend on the geometric data $\left(E,  \mathcal L, \sigma \right)$, as illustrated by the following result.

\bpr{PI}   \textnormal{[$n$]} Both of the algebras  $S$ and $B=B(E, \mathcal L, \sigma) \cong S/(Sg_1+Sg_2)$ are module-finite over their center, if and only if the automorphism $\sigma$ has finite order. In this case, both $S$ and  $B$ satisfy a polynomial identity (i.e., are PI) and are of  PI degree $n := |\sigma| <\infty$.  
\index{n5@$n$}
\qed
\epr

\begin{proof}
Suppose that the automorphism $\sigma$ has finite order. Since $\mathcal L$ is ample and $\sigma$-ample, the algebra $B$ is module-finite over its center by \cite[Corollary 2.3]{SmithTate}. Moreover, so is $S$ by \cite[Theorem 3.7(c)]{SmithTate}.

On the other hand, suppose that the automorphism $\sigma$ has infinite order. Now the center of $S$ is the polynomial ring $\kk[g_1, g_2]$ by \cite[Proposition 6.12]{LevSmith}. By comparing the Hilbert series of $S$ to that of $\kk[g_1,g_2]$, we get $S$ cannot be module-finite over its center $\kk[g_1,g_2]$. As a consequence, $B$ is not module-finite over its center as $B$ is a quotient ring of $S$ by a regular sequence of homogenous central elements of degree two by \cite[Lemma 3.6(b)]{SmithTate}.

To verify the last statement, note that $B$ has PI degree $n$ since it has a localization isomorphic to $\kk(E)[t^{\pm 1}; \sigma]$ (see Lemma 2.9), 
which in turn also has PI degree $n$. We will also see later in Corollary~\ref{crepB} that all nontrivial irreducible representations of $B$ have dimension $n$, which provides another proof that $B$ has PI degree $n$.

Finally, we see that the PI degree of $S$ is also $n$ as follows. First, PIdeg($S$) $\geq$ PIdeg($B$)  since $B$ is a homomorphic image of $S$, and  recall that PIdeg($B$) $=n$ from above. So, it suffices to show that ${\rm PIdeg}(S) \leq n$. Take the central element defined in \cite[(4-2)]{Smith1993}, which is denoted by $g$  later in Notation~\ref{ncentral}. (It is denoted by ``$c$" in \cite{Smith1993}.) Let $M$ be an irreducible representation of $S$ that is $g$-torsionfree. By \cite[Theorem III.1.7]{BrownGoodearl:book}, we know the Azumaya locus of $S$ is dense in ${\rm maxSpec}(Z)$. So, we only need to show $\dim M\leq n$. Recall  $s = n/(n,2)$ and note that $M$ is the quotient of a $g$-torsionfree {\it fat point module} $F$ of multiplicity $s>1$ by \cite[Lemma~4.1]{SmithStan} and \cite[Theorem~7.7(c)]{Smith1993}. Fat point modules will be discussed in detail in Section~\ref{sec:fat}; they serve as representatives of the simple objects of the quotient category $S$-qgr which is the category of graded $S$-modules modulo those that are bounded above ($M_n = 0$ for all $n \gg 0$). Now if $n$ is odd, then $F$ and its shift $F[1]$ are equivalent in $S$-qgr by \cite[Corollary~8.7]{Smith1993}. So, $\dim(M) \leq {\rm mult}(F) = s = n$ by \cite[discussion after Proposition~3.19]{Smith1993}. On the other hand, if $n$ is even, then $F$ and $F[2]$ are equivalent in $S$-qgr by \cite[Proposition~5.4]{Smith1993}, and hence, ${\rm dim}(M) \leq 2\, {\rm mult} (F) = 2s = n$ by \cite[discussion after Proposition~3.19]{Smith1993}. Thus, $\dim(M) \leq n$, as desired.
\end{proof}

\bre{PIS}
The claim that the PIdeg($S) =n$ first appeared in \cite[Theorem 8.8(2)]{Smith1993}, but the stronger statement that the algebra $S[c^{-1}]$  in \cite[Theorem~8.8(1)]{Smith1993} is Azumaya (which is used to prove \cite[Theorem~8.8(2)]{Smith1993}) is incorrect; there, ``$c$" is our ``$g$" in Notation~\ref{ncentral}. We establish later in Theorem~\ref{tAzumaya} and 
 Theorem~\ref{tSEven} that the Azumaya locus of $S$ is equal to the smooth locus of 
 $Y:=\maxSpec(Z(S))$, and in the case when $n$ is even, the singular locus of $Y$ does not lie in the
  hypersurface $\mathbb{V}(g)$. That is, there exist $g$-torsionfree irreducible representations of $S$ of 
  dimension $<n$ when $n$ is even.
\ere


\subsection{Center of $S$}
As mentioned in the proof of Proposition~\ref{pPI}, the center of a 4-dimensional Sklyanin algebra $S$ is equal  to $\kk[g_1,g_2]$ when $n:=|\sigma|$=PIdeg($S$) is infinite. On the other hand, one expects that both $S$ and $B$ have  a large center when the PI degree of $S$ is finite. We record here several results of Smith-Tate \cite{SmithTate} pertaining to the center $Z$ of $S$ when PIdeg($S$) $= n<\infty$.

\ble{centerB}  \textnormal{[$E''$]} \cite[Corollary~2.8]{SmithTate} 
 Given the geometric data $\left(E, \mathcal L, \sigma \right)$ from 
\dlref{geomS}, suppose that $|\sigma|=:n<\infty$. 
Take $E'':=E/\langle \sigma \rangle$\index{E0primeprime@$E''$} 
so that $E\to E''$ is a cyclic \'{e}tale cover of degree~$n$. Recall \leref{embedB} and let $D''$ be the image of $D$ on $E''$ and let $V''$ denote $H^0(E'', \mathcal{O}_{E''}(D''))$.  

Then, the center of $B$ is the intersection of $B$ with $\kk(E'')[t^{\pm n}]$, which is equal to $\kk[V'' t^n]$, and this is also a twisted homogeneous coordinate ring of $E''$ for an embedding of $E'' \subseteq \mathbb{P}^3$. \qed
\ele

Central elements of $B$ lift to central elements of $S$ as described below. 
We will identify $S_1 \cong B_1$ via the canonical projection mentioned in Lemma~\ref{lfactor-g}.

\bde{good} [$s$] \cite{SmithTate} Take $|\sigma|=:n<\infty$ and let $s$\index{s5@$s$} be the value $n/(n,2)$.  A section 
of $B_1 := H^0(E, \mathcal{L})$ is called {\em{good}} if its divisor of zeros 
is invariant under $\sigma^s$ and
consists of distinct points, whose orbits under the group $\lcor \sigma \rcor$ 
do not intersect.  A {\em{good basis}} of $B_1$ is a basis that consists of good elements so that the $s$-th powers of these elements generate $B_n$ if $n$ is odd or generate $(B^{\langle \sigma^2 \rangle})_{n/2}$ if $n$ is even.
\ede

\bnota{rho} [$\rho$]  Take $|\sigma|=:n <\infty$. As mentioned at \cite[top of page~31]{SmithTate},  $\sigma^s$ fixes the class $[\mathcal L]$ in $\mathrm{Pic}\, E$. So the 
automorphism $\sigma^s$ of $E$ induces an automorphism of $B_1$ via the identification $\mathcal L^{\sigma^s}\cong \mathcal L$. 
This automorphism of $B$ will be denoted by  $\rho$\index{r2ho@$\rho$}. 
\enota

By \cite[Lemma~3.4 and its proof]{SmithTate}, there is a unique lifting of $\rho$ to a graded automorphism of $S$ and $\rho^2$ is the identity. Further,
the elements $x_i^2$ and the central elements $g_1$ and $g_2$ are all $\rho$-invariant.

\ble{actionrho}  \cite[proof of Lemma~3.4, page~46]{SmithTate} When $|\sigma|=:n<\infty$, we can take as  a good basis of $B_1$ the set of generators $x_0, x_1, x_2, x_3$ of $S$ where each $x_i$ is a $\rho$-eigenvector. If $n$ is odd, then $\rho$ is the identity. If $n$ is even, then $\rho = \sigma^{n/2}$ of order 2 and each $\rho$-eigenspace of $B_1$ is two-dimensional eigenspace with one having eigenvalue $1$ and the other having eigenvalue $-1$. \qed
\ele

\bpr{centerS} \textnormal{[$Z$, $z_i$, $F_i$, $u_i$, $\Phi_i$, $\ell_i$,  $h_i$, $f_i$, $E'$]}  \cite[Theorems 3.7,~4.6,~4.9,~4.10]{SmithTate}. 
The center $Z$\index{Z0@$Z$} of $S$, of PI degree $n < \infty$, is given as follows. 

\begin{enumerate}
\item The center $Z$ is generated by four algebraically independent elements $z_0,z_1, z_2, z_3$\index{z3s@$z_0, \ldots, z_3$}
of degree $n$ along with~$g_1,g_2$ in \eqref{g-elem}, 
subject to two relations $F_1$, $F_2$\index{F012@$F_1, F_2$} of degree $2n$. In fact, there is a choice of generators $z_i$ of the form 
$$z_i = x_i^n + \textstyle \sum_{1 \leq j < n/2} c_{ij} x_i^{n-2j}$$
where $\{x_0,x_1, x_2, x_3\}$ is any
good basis of $B_1$ and $c_{ij} \in \kk[g_1,g_2]_{2j}$. 
\smallskip

\item If $n$ is even, then there exist elements $u_0,u_1, u_2, u_3$\index{us@$u_0, \ldots, u_3$}  of degree $n/2$ whose image in the Veronese subalgebra, $B^{(n/2)}$, of $B$ generate the center. Here, $z_i = u_i^2$.
\smallskip

\item If $n$ is odd, then for $i =1,2$,
$$F_i = \Phi_i(z_0,z_1,z_2,z_3) + h_i(g_1,g_2),$$
where $\Phi_1, \Phi_2$\index{P300hi@$\Phi_1, \Phi_2$} 
are the quadratic homogeneous defining polynomials of the elliptic curve 
$$E'' = E/\langle \sigma \rangle \subset \mathbb{P}^3 = \mathbb{P}(H^0(E'', \mathcal{L}'')^*)
$$
with $\mathcal{L}''$ the descent of $\mathcal{L}_n$ to $E''$,  and where $h_1, h_2$\index{h12@$h_1,h_2$} 
are homogeneous degree $s$ forms in variables $g_1, g_2$ having no common factor. 
\smallskip

\item If $n$ is even, then for $i =1,2$,
\[
\begin{array}{c}
\medskip
F_i = \Phi_i(z_0,z_1,z_2,z_3) + \ell_i(z_0,z_1,z_2,z_3)h_i(g_1,g_2)+h_i(g_1,g_2)^2, \quad \text{ and}\\

f_i(u_0, u_1,u_2,u_3) + h_i(g_1,g_2) = 0 \quad \text{ in $Z(S^{(2)})$}.
\end{array}
\]
Here, $E'' = E/\langle \sigma \rangle = \mathbb{V}(\Phi_1, \Phi_2)$, $\ell_1, \ell_2$\index{l212@$\ell_1, \ell_2$} 
are linear forms and $f_1, f_2$\index{f12@$f_1, f_2$} are linearly independent  quadratic forms in the variables $u_i$ defining the elliptic curve 
 $$
 E':=E/\langle \sigma^2 \rangle \subset \mathbb{P}^3 = \mathbb{P}(H^0(E', \mathcal{L}')^*),
 \index{E0prime@$E'$}
 $$
 where $\mathcal{L}'$ is the descent of $\mathcal{L}_s$ to $E'$, and $h_1, h_2$ are homogeneous degree $s$ forms in variables $g_1, g_2$ having no common factor.
 \qed
\end{enumerate}
\epr

We will use the defining relations of the center $Z$ of PI Sklyanin algebras $S$ mentioned above more extensively in Section \ref{sec:sing} to figure out the singular locus of maxSpec($Z$).


\subsection{Symmetries of $S$ and of their centers}
Now we turn our attention to symmetries of 4-dimensional Sklyanin algebras $S$ and of their respective centers $Z$. Recall that we do not assume that $S$ satisfies a polynomial identity in this section. 

\bpr{H4} \textnormal{[$H_4$, $\epsilon_1$, $\epsilon_2$, $\epsilon$, $a$, $b$, $c$, $\xi$]} 
\index{a@$a$} 
\index{bb@$b$} 
\index{c1@$c$}
\index{xi@$\xi$}
\index{epsilon@$\epsilon, \epsilon_1, \epsilon_2$}
\cite[Section~2]{SmithStan} \cite[Section~2.7]{ChirSmith} The group of graded automorphisms of $S$ contains the Heisenberg group $H_4$ of order 64 which is presented as
$$
H_4 = \langle \epsilon_1, \epsilon_2, \epsilon ~:~ \epsilon_1^4 = \epsilon_2^4= \epsilon^4 = 1,~\; \epsilon \epsilon_1 = \epsilon_1\epsilon, ~\; \epsilon \epsilon_2 = \epsilon_2 \epsilon,~\; \epsilon_1 \epsilon_2 = \epsilon \epsilon_2 \epsilon_1 \rangle.
\index{H04@$H_4$}
$$
Here, $\epsilon$ scales the generators of $S$ by $-i$ and 
\[
\begin{array}{lrrrr}
\epsilon_1:(x_0, x_1, x_2, x_3) \mapsto 
(b^{\frac{1}{2}}c^{\frac{1}{2}}\xi^{-1}\;x_1, &b^{\frac{-1}{2}}c^{\frac{-1}{2}}\xi\;x_0, &
b^{\frac{1}{2}}c^{\frac{-1}{2}}\xi\;x_3, &-b^{\frac{-1}{2}}c^{\frac{1}{2}}\xi^{-1}\;x_2),\\  
\epsilon_2:(x_0, x_1, x_2, x_3) \mapsto 
(a^{\frac{1}{2}}c^{\frac{1}{2}}\xi^{-1}\;x_2, &-a^{\frac{1}{2}}c^{\frac{-1}{2}}\xi^{-1}\;x_3, &
a^{\frac{-1}{2}}c^{\frac{-1}{2}}\xi\;x_0, & a^{\frac{-1}{2}}c^{\frac{1}{2}}\xi\;x_1)
\end{array}
\]
for $a^2 = \alpha,~ b^2 = \beta, ~c^2 = \gamma$, $i := \sqrt{-1} = e^{\frac{\pi i}{2}}$, and $\xi = e^{\frac{3\pi i}{4}}$.  In particular, $\epsilon^2$ scales the generators of $S$ by $-1$ and 
$$\epsilon_1^2:(x_0,x_1,x_2,x_3) \mapsto (x_0, x_1,-x_2,-x_3)$$
$$\epsilon_2^2:(x_0,x_1,x_2,x_3) \mapsto (x_0, -x_1,x_2,-x_3).$$

\vspace{-.2in}
 \qed
\epr

The following consequences hold mostly by direct computation. First, consider the following notation.

\bnota{E2} \textnormal{[$E_2$, $\omega$]} Let $E_2$\index{E2@$E_2$} 
be the set of the (four) 2-torsion points of $E$, and let a point of $E_2$ be denoted by $\omega$\index{omega@$\omega$}.
\enota

\bco{fixZ}  \textnormal{[$N_4$]} Let $N_4$\index{N4@$N_4$} 
be the subgroup of $H_4$ generated by $\epsilon^2$, $\epsilon_1^2$, $\epsilon_2^2$, which is normal and is isomorphic to $\mathbb{Z}_2 \times \mathbb{Z}_2 \times \mathbb{Z}_2$. We have the statements below.

\begin{enumerate}
\item The element $\epsilon$ has four 1-dimensional eigenspaces each with eigenvalue $-i$.
\smallskip

\item The element $\epsilon_1^2$ has two 2-dimensional eigenspaces $\langle x_0, x_1 \rangle$ and $\langle x_2, x_3 \rangle$ with eigenvalue $1$ and $-1$, respectively;  the element $\epsilon_2^2$ has two 2-dimensional eigenspaces  $\langle x_0, x_2 \rangle$ and $\langle x_1, x_3 \rangle$ with eigenvalue $1$ and $-1$, respectively; and the element $\epsilon_1^2\epsilon_2^2$ has two 2-dimensional eigenspaces  $\langle x_0, x_3 \rangle$ and $\langle x_1, x_2 \rangle$ with eigenvalue $1$ and $-1$, respectively.
\smallskip

\item The element $\rho$ from Notation~\ref{nrho} is  an element of the quotient group $N_4/\langle \epsilon^2 \rangle$.
\smallskip

\item The subspace $\kk g_1+\kk g_2\subset Z(S)_2$ is a 2-dimensional irreducible representation of $H_4$. In the basis $\{g_1, g_2\}$, the $H_4$-action is given by
$\epsilon \mapsto \text{diag}(-1,-1)$, and
\[
\epsilon_1 \mapsto \frac{i}{bc}\begin{pmatrix} 1 & 1 \\ -1-\beta\gamma & -1\end{pmatrix},\quad \epsilon_2 \mapsto \frac{i}{ac}\begin{pmatrix} 1 & \frac{1+\alpha}{1-\beta} \\  -1-\gamma& -1\end{pmatrix};
\]
this implies that the subgroup $N_4$ fixes $g_1,g_2$. 
\smallskip

\item The group of graded automorphisms of $B$ contains $H_4$.
\smallskip

\item The subspace $\kk z_0+\kk z_1+\kk z_2+\kk z_3 \subset Z(S)_n$ is a 4-dimensional representation of $H_4$. In the basis $\{ z_0, z_1, z_2, z_3 \}$, the $H_4$-action is given by  
\begin{gather*}
\epsilon_1 \mapsto \begin{pmatrix} 0 &  b^{-\frac{n}{2}}c^{-\frac{n}{2}} \xi^{n}   & 0 &0  \\  b^{\frac{n}{2}}c^{\frac{n}{2}} \xi^{-n}  & 0 & 0 & 0\\  0 & 0& 0 & (-1)^nb^{-\frac{n}{2}}c^{\frac{n}{2}} \xi^{-n}  \\ 0 & 0& b^{\frac{n}{2}}c^{-\frac{n}{2}} \xi^{n} &0 \end{pmatrix},\\ 
\bigskip
\epsilon_2 \mapsto \begin{pmatrix} 0 & 0 & a^{-\frac{n}{2}}c^{-\frac{n}{2}} \xi^{n}  & 0 \\ 0 & 0 & 0 &a^{-\frac{n}{2}}c^{\frac{n}{2}} \xi^{n}\\  a^{\frac{n}{2}}c^{\frac{n}{2}} \xi^{-n}& 0& 0 & 0 \\ 0 &  (-1)^na^{\frac{n}{2}}c^{-\frac{n}{2}} \xi^{-n}& 0 & 0
\end{pmatrix},
\end{gather*}
and $\epsilon\mapsto \mathrm{diag}((-i)^n,(-i)^n,(-i)^n,(-i)^n)$; 
this implies that $N_4$ fixes $z_0, z_1, z_2, z_3$ when $n$ is even. 
\end{enumerate}
\eco
\begin{proof} The first statement  about $N_4$ is clear.

\smallbreak
(1) and (2) are clear.

\smallbreak
(3) By Lemma~\ref{lactionrho}, we only treat the case when $n:=|\sigma|$ is even. Here, $\rho = \sigma^{n/2}$ is given by translation by a point $\omega$ of $E$ of order 2. We aim to pick off the graded automorphism of $S_1 = B_1$ to which $\rho$ corresponds. To do this, first identify the generators $x_i$ of $S$ with the coordinates $v_i$ of $E$. Then, take $\kk = \mathbb{C}$ and write $E$ as $\mathbb{C}/(\mathbb{Z}+\mathbb{Z}\Lambda)$ for some $\Lambda \in \mathbb{C}$\index{L3ambda@$\Lambda$} 
with Im($\Lambda$)$>0$ as in \cite[Sections~2.9-2.13]{SmithStafford}. Now $\rho$ corresponds to a translation by
a nontrivial 2-torsion point $\omega$ of $\mathbb{C}/(\mathbb{Z}+\mathbb{Z}\Lambda)$, that is 
$$\omega \in \left\{\textstyle \frac{1}{2}, ~\frac{\Lambda}{2}, ~\frac{1+\Lambda}{2}\right\}.$$
Using the notation of \cite{SmithStafford}, note that $$j: \mathbb{C}/(\mathbb{Z}+\mathbb{Z}\Lambda) \to E, \quad z \mapsto [g_{11}(z) : g_{00}(z) : g_{01}(z) : g_{10}(z)]$$ is an isomorphism where $$g_{pq}(z) = \theta_{pq}(2z) \theta_{pq}(\omega) \gamma_{pq}$$ are holomorphic functions on $\mathbb{C}$ so that the theta functions satisfy the conditions
$$\theta_{pq}(z+1) = (-1)^p \theta_{pq}(z) \quad \text{and} \quad 
\theta_{pq}(z+ \Lambda) =  -\text{exp}(\Lambda + 2 z + q) \theta_{pq}(z)$$
for $p,q = 0,1$, and $\gamma_{00} = \gamma_{11} = i$, $\gamma_{01} = \gamma_{10} = 1$. In particular, $\rho(j(z)) = j(z + \omega)$.

Now take $\omega = \frac{1}{2}$. Then,
$$\textstyle g_{pq}(z+\frac{1}{2}) ~=~ \theta_{pq}(2z+1)\theta_{pq}(\frac{1}{2}) \gamma_{pq} ~=~ (-1)^p\theta_{pq}(2z)\theta_{pq}(\frac{1}{2}) \gamma_{pq} ~=~ (-1)^p g_{pq}(z).$$
Therefore, we get that
$$\rho([g_{11}(z) : g_{00}(z) : g_{01}(z) : g_{10}(z)]) = [-g_{11}(z) : g_{00}(z) : g_{01}(z) : -g_{10}(z)],$$
which is realized as the coset of $\epsilon_1^2 \epsilon_2^2$ in $N_4/\langle \epsilon^2 \rangle$ by Proposition~\ref{pH4}.

For $\omega = \frac{\Lambda}{2}$, we have that
$$\textstyle g_{pq}(z+\frac{\Lambda}{2}) ~=~ \theta_{pq}(2z+\Lambda)\theta_{pq}(\frac{\Lambda}{2}) \gamma_{pq} ~=~\text{exp}(-4\pi iz) \text{exp}(-\pi i \Lambda) \text{exp}(-\pi i q)g_{pq}(z).$$
So,
$$\rho([g_{11}(z) : g_{00}(z) : g_{01}(z) : g_{10}(z)]) = [-g_{11}(z) : g_{00}(z) : -g_{01}(z) : g_{10}(z)]$$
which is realized as the coset of $\epsilon_2^2$ in $N_4/\langle \epsilon^2 \rangle$ by Proposition~\ref{pH4}.

For $\omega = \frac{1+\Lambda}{2}$, we have that
$$\textstyle g_{pq}(z+\frac{1+\Lambda}{2}) ~=~ \theta_{pq}(2z+1+\Lambda)\theta_{pq}(\frac{1+\Lambda}{2}) \gamma_{pq} ~=~(-1)^p\text{exp}(-4\pi iz) \text{exp}(-\pi i \Lambda) \text{exp}(-\pi i q)g_{pq}(z).$$
So, $\rho$ is realized as the coset of $\epsilon_1^2$ in $N_4/\langle \epsilon^2 \rangle$ by Proposition~\ref{pH4}.

Finally, collect all of the structure constants to define the automorphism $\rho = \sigma^{n/2}$ in $\text{Aut}(S_1)$ in an algebraically closed field 
$\kk'$ over $\mathbb{Q}$. Then a standard base change argument shows that $\rho$ is realized as the graded automorphisms $\pm \epsilon_1^2 \epsilon_2^2$, $\pm \epsilon_2^2$, $\pm \epsilon_1^2$ of the $\kk'$-algebra $S$ and of the $\kk$-algebra $S$, respectively, for $\omega = \frac{1}{2}, ~\frac{\Lambda}{2}, ~\frac{1+\Lambda}{2}$.

 \smallbreak
(4) This is a direct calculation using  \eqref{g-elem}, \eqref{albega-cond1}, and Proposition \ref{pH4}. 

\smallbreak
(5) This follows from part (4) and Lemma~\ref{lfactor-g}.

\smallbreak
(6) By part (5) and Proposition~\ref{pcenterS}(1), we can consider the action of $H_4$ on $x_0^n$, $x_1^n$, $x_2^n$, $x_3^n$. Now this part follows from Proposition \ref{pH4}.  
\end{proof}

Note that the subgroup $\langle \epsilon_1^2, \epsilon_2^2 \rangle$ plays a crucial role in constructing Chirvasitu-Smith's  {\it exotic elliptic algebras of dimension 4} \cite{ChirSmith}; see also the work of Davies  \cite{Davies} where these algebras appeared independently. On the other hand, similarly to Corollary 2.21(4), the space $\kk g_1 \oplus \kk g_2$ was realized as a representation of $H_4$ in work of Kevin DeLaet \cite{DeLaet-Heis} in the case when $S$ is generic.

\bnota{rhos} [$\rho_1$, $\rho_2$, $\rho_3$] Recall from the proof of Corollary~\ref{cfixZ}(3) that if $n=|\sigma|$ is even, then $\rho$ is identified with a nontrivial element of $N_4/\langle \epsilon^2 \rangle$ via an identification with $\{\frac{1}{2}, ~\frac{\Lambda}{2}, ~\frac{1+\Lambda}{2}\}$, the set of nontrivial 2-torsion points $\omega$ of $\mathbb{C}/(\mathbb{Z}+\mathbb{Z}\Lambda) \cong E$. 

We take $\rho_1$ (respectively, $\rho_2$, $\rho_3$)\index{rho3@$\rho_1, \rho_2, \rho_3$}
to be the automorphism $\rho \in$ Aut$(E)$ corresponding to $\omega = \frac{1}{2}$ (respectively, $\omega = \frac{\Lambda}{2}$, $\frac{1+ \Lambda}{2}$), which in turn is identified with the coset of $\epsilon_1^2 \epsilon_2^2$ (respectively, of $\epsilon_2^2$, $\epsilon_1^2$) in $N_4/\langle \epsilon^2 \rangle$.
\enota

\bpr{ZEven} \textnormal{[$a_i$]}
Let $n<\infty$ be even and recall $s=n/2$. Retain the notation above and recall the notation in Proposition \ref{pcenterS}, namely:
$$Z(S)=\kk[z_0,z_1,z_2,z_3,g_1,g_2]/(F_1,F_2),\quad  \text{with $\mathrm{deg}(z_i)=n$, ~$\mathrm{deg}(g_i)=2$, ~ $\mathrm{deg}(F_i)=2n$},$$
and for $i = 1,2$, $$F_i = \Phi_i(z_0,z_1,z_2,z_3) + \ell_i(z_0,z_1,z_2,z_3)h_i(g_1,g_2) + h_i(g_1,g_2)^2.$$
We obtain the precise description of both $\Phi_i$ and $a_i:=-\ell_i/2$\index{a1@$a_1, a_2$},
as given below. Here, $\lambda, \mu \in \kk$ and $\mu\neq 0$.
\begin{itemize}
\item[(1)] If $\rho=\rho_1$, then  
$$\Phi_1 = a_1^2 - z_0z_3,\quad \Phi_2 = a_2^2 - z_1z_2, \quad \text{with}$$
\begin{align*}
a_1&\, =\lambda (z_0+a^sb^s\xi^{3n}z_3)+\mu(z_1+a^sb^{-s}\xi^{n}z_2),\\ 
a_2&\, =\mu (b^{-s}c^{-n}\xi^{2n}z_0+a^sc^{-n}\xi^{3n}z_3)+\lambda (b^sz_1+a^s\xi^{3n}z_2).
\end{align*}
\item[(2)]  If $\rho=\rho_2$, then $4 | n$ and 
$$\Phi_1 = a_1^2 - z_0z_2,\quad \Phi_2 = a_2^2 - z_1z_3, \quad \text{with}$$
\begin{align*}
a_1&\, =\lambda (z_0+a^sc^s\xi^{n}z_2)+\mu(z_1+a^sc^{-s}\xi^{n}z_3),\\ 
a_2&\, =\mu (b^{-n}c^{-s}\xi^nz_0+a^sb^{-n}z_2)+\lambda(c^s\xi^nz_1+a^sz_3).
\end{align*}
\item[(3)]  If $\rho=\rho_3$, then  $4 | n$ and 
$$\Phi_1 = a_1^2 - z_0z_1,\quad \Phi_2 = a_2^2 - z_2z_3, \quad \text{with}$$
\begin{align*}
a_1&\, =\lambda (z_0+b^sc^s\xi^nz_1)+\mu(z_2+b^sc^{-s}\xi^nz_3),\\ 
a_2&\, =\mu (a^{-n}c^{-s}z_0+a^{-n}b^s\xi^nz_1)+\lambda(c^sz_2+b^s\xi^nz_3).
\end{align*}
\end{itemize}
\epr

\begin{proof}
We only show the details for the case (1) and the remaining cases can be verified in a similar fashion. We follow the proof in \cite[Theorem 4.10]{SmithTate} to get the forms of $\Phi_1$ and $\Phi_2$; namely, $\rho_1$ has eigenspaces $\langle x_0, x_3\rangle$ and $\langle x_1, x_2 \rangle$ by Corollary~\ref{cfixZ}(2). Using Corollary~\ref{cfixZ}(6), observe that 
\[
\begin{array}{lll}
\epsilon_1(z_0z_3) = c^n \xi^{-2n} z_1 z_2, \quad
& \epsilon_2(z_0z_3) = c^n z_1 z_2, \quad
& \epsilon_1 \epsilon_2(z_0z_3) = \xi^{2n} z_0 z_3,\\
\epsilon_1(z_1z_2) = c^{-n} \xi^{2n} z_0 z_3, \quad
& \epsilon_2(z_1z_2) = c^{-n} z_0 z_3, \quad
& \epsilon_1 \epsilon_2(z_1z_2) = \xi^{-2n} z_1 z_2.\\
\end{array}
\]
Working in $B$ via Proposition~\ref{pcenterS}(3,4) and Corollary~\ref{cfixZ}(5), we can take $a_1^2 = z_0 z_3$ and $a_2^2 = z_1 z_2$, so
\[
\begin{array}{lllll}
z_0 z_3 &= a_1^2 &= c^n \xi^{-2n} \epsilon_1(a_2)^2 &= c^n \epsilon_2(a_2)^2 &= \xi^{-2n} \epsilon_1 \epsilon_2(a_1)^2,\\
z_1 z_2 &= a_2^2 &= c^{-n} \xi^{2n} \epsilon_1(a_1)^2 &= c^{-n} \epsilon_2(a_1)^2 &= \xi^{2n} \epsilon_1 \epsilon_2(a_2)^2.
\end{array}
\]
If $a_i = \eta_{i0}z_0 + \eta_{i1}z_1 + \eta_{i2}z_2 + \eta_{i3}z_3$ for $\eta_{ij} \in \kk$, for $i =1,2$, then by comparing coefficients of $z_k$ in the equations above we get $a_1$ and $a_2$ as claimed. For instance,
$$\eta_{10} = \pm \eta_{13} a^{-s} b^{-s} \xi^n, \quad
\eta_{11} = \pm \eta_{12} a^{-s} b^{s} \xi^{-n}, \quad 
\eta_{12} = \pm \eta_{11} a^{s} b^{-s} \xi^{n}, \quad
\eta_{13} = \pm \eta_{10} a^{s} b^{s} \xi^{-n},$$
$$\eta_{10} = \pm \eta_{21} b^{-s}, \quad
\eta_{11} = \pm \eta_{20} b^{s} c^n \xi^{-2n}, \quad 
\eta_{12} = \pm \eta_{23} b^{-s} c^n \xi^{-2n}, \quad
\eta_{13} = \pm \eta_{22} b^{s},$$
amongst other conditions on $\eta_{ij}$ that imply that the $\pm$s are unnecessary.

Finally, we need to show that $\mu\neq 0$ in each case. Suppose $\mu=0$, then $\Phi_1=a_1^2-z_0z_3$ becomes a quadratic equation in terms of $z_0,z_3$. Since $\kk$ is algebraically closed, $\mathbb{V}(\Phi_1)$ is an union of two planes, one of which contains $E''$ since $E''$ is irreducible and $E''\subset \mathbb{V}(\Phi_1)$. But $E''$ is not contained in a hyperplane, so we have a contradiction.  
\end{proof}

Now the result below follows immediately from Corollary~\ref{cfixZ}(6) and Proposition~\ref{pZEven}.

\bco{neven-H4} When $n$ is even, the vector spaces $\kk a_1 + \kk a_2$ and $\kk h_1 + \kk h_2$ admit the structure of a 2-dimensional irreducible representation of $H_4$ both via the following actions.
\begin{itemize}
\item[(1)] If $\rho = \rho_1$, then $\epsilon$ acts by $(-1)^s$ and the action of $\epsilon_1$ and $\epsilon_2$ are given by
$$\epsilon_1 \mapsto \begin{pmatrix} 0 & c^{-s}\xi^n\\c^s \xi^{-n} &0 \end{pmatrix} \quad \text{ and } \quad
\epsilon_2 \mapsto \begin{pmatrix} 0 & c^{-s}\\c^s  &0 \end{pmatrix}.$$
\item[(2)] If $\rho = \rho_2$, then $\epsilon$ acts by $(-1)^s$ and the action of $\epsilon_1$ and $\epsilon_2$ are given by
$$\epsilon_1 \mapsto \begin{pmatrix} 0 & b^{-s}\\b^s  &0 \end{pmatrix} \quad \text{ and } \quad
\epsilon_2 \mapsto \begin{pmatrix} 1 & 0\\0  &1 \end{pmatrix}.$$
\item[(3)] If $\rho = \rho_3$, then $\epsilon$ acts by $(-1)^s$ and the action of $\epsilon_1$ and $\epsilon_2$ are given by
$$\epsilon_1 \mapsto \begin{pmatrix} 1 & 0\\0 &1 \end{pmatrix} \quad \text{ and } \quad
\epsilon_2 \mapsto \begin{pmatrix} 0 & a^{-s}\xi^{-n}\\a^s\xi^n  &0 \end{pmatrix}.$$
\end{itemize}

\qed
\eco

\sectionnew{Singular loci of the maximal spectra of the centers of the PI 4-dimensional Sklyanin algebras $S$} \label{sec:sing}

As in Hypothesis~\ref{nhyp}, let $S$ be a 4-dimensional Sklyanin algebra that satisfies a polynomial identity (i.e., is PI) of PI degree $n$. Recall that $n = |\sigma|$, where $\sigma \in$ Aut($E$) and $\left(E,  ~\mathcal{O}_{\mathbb{P}^3}(1)|_E, ~\sigma\right)$ is the projective algebro-geometric data attached to $S$. The purpose of this section is to provide a detailed description of maxSpec of the center $Z:=Z(S)$, and this study depends on the parity of $n$. The main results are given in  Theorems~\ref{tSOdd} and~\ref{tSEven} below for the cases when $n$ is odd and even, respectively.

Let us set some notation that will be used throughout this section and establish a couple of preliminary results.

\bnota{not:Y}[$Y$, $Y^{sing}$, $Y_{\gamma_1, \gamma_2}$, $(Y_{\gamma_1, \gamma_2})^{sing}$, $(Y^{sing})_{\gamma_1, \gamma_2}$, $Y^{symp}_0$] 
Let $Y$\index{Y@$Y$} 
denote the affine variety maxSpec($Z(S)$), let $Y^{sing}$\index{Y1sing@$Y^{sing}$} 
denote the singular locus of $Y$, and let
$$
Y_{\gamma_1, \gamma_2} := Y \cap \mathbb{V}(g_1 - \gamma_1, g_2 - \gamma_2), \quad \text{ for } \gamma_1, \gamma_2 \in \kk. 
\index{Y2gama@$Y_{\gamma_1, \gamma_2}$}
$$
We denote by $(Y_{\gamma_1,\gamma_2})^{sing}$ \index{Y2gamasing@$(Y_{\gamma_1, \gamma_2})^{sing}$}
the singular locus of the subvariety $Y_{\gamma_1,\gamma_2}$ and denote
$$
(Y^{sing})_{\gamma_1,\gamma_2}:=Y^{sing}\cap Y_{\gamma_1,\gamma_2}.
\index{Y3singgama@$(Y^{sing})_{\gamma_1,\gamma_2}$}
$$
Note that $Y^{sing} = \bigcup_{\gamma_1,\gamma_2\in \kk} (Y^{sing})_{\gamma_1,\gamma_2}$. On the other hand, denote 
$$
Y^{symp}_0:=\bigcup_{\gamma_1, \gamma_2 \in \kk} (Y_{\gamma_1, \gamma_2})^{sing}.
\index{Y1Pois@$Y^{symp}_0$}
$$
This will be the variety of {\em{symplectic points}} of the Poisson bracket on $Y$ that we construct in Sections~\ref{sec:constr-Pord} 
and~\ref{sec:bracket} (i.e., where the Poisson structure vanishes, or equivalently the variety of 0-dimensional symplectic cores of the Poisson structure), 
see Corollary~\ref{cYPois} for more details. 
We will also have that $(Y_{\gamma_1, \gamma_2})^{sing} \supseteq (Y^{sing})_{\gamma_1,\gamma_2}$, so $Y^{symp}_0 \supseteq Y^{sing}$, but these containments could be strict (see, Theorem~\ref{tSOdd}(3)).
\enota

Before proceeding with the study of $Y^{sing}$, we recall some results about line modules of $S$ from Levasseur-Smith \cite{LevSmith}. A {\it line module} of $S$ is a graded cyclic module of $S$ of Hilbert series $(1-t)^{-2}$; such modules are in correspondence with secant lines of the elliptic curve $E$ of the point scheme of $S$ \cite[Theorem 4.5]{LevSmith}. 

\bnota{LineA}[$M(p,q)$, $\ell_{p,q}$, $\Omega(z)$] 
For each $z\in E$, let 
$$\{M(p,q)~|~ p,q\in E, ~p+q=z\}
\index{M0pq@$M(p,q)$}
$$
denote the family of line modules of $S$ corresponding to the secant line $\ell_{p,q}$\index{l2pq@$\ell_{p,q}$} 
to $E$ at $z$. All such line modules have a common central annihilator of degree $2$ by \cite[Corollary~6.6]{LevSmith}.
Denote the central degree $2$ annihilator of  $M(p,q)$, with $p+q=z\in E$, by 
$$
\Omega(z)\in \kk g_1+\kk g_2.
\index{O0megaz@$\Omega(z)$}
$$
\enota

By \cite[Corollary 6.9]{LevSmith}, the only equalities among these annihilators are 
\begin{align}\label{EOmega}
\Omega(z)=\Omega(-z-2\tau), \quad \text{ for } z \in E.
\end{align}
Recall that $\tau$ is the point in $E$ such that the automorphism $\sigma\in {\rm Aut}(E)$ is given by $\sigma(p)=p+\tau$ for any $p\in E$.

\smallbreak Now we recall some facts on  the geometry of line modules from both \cite{LevSmith} and \cite{SmithTate}.
 
\bnota{ri}[$\overline{E}$, $\overline{\mathbb{P}}$, $\overline{z}$, $r_i$]
Let $\overline{E}$\index{E1bar@$\overline{E}$} denote $E'' = E/\langle \sigma \rangle$ when $n$ is odd, or denote $E' = E/\langle \sigma^2 \rangle$ when $n$ is even. Likewise, take $\overline{\mathbb{P}}$\index{P0bar@$\overline{\mathbb{P}}$} to be $\mathbb{P}(H^0(E'', \mathcal{L}'')^*)$ when $n$ is odd, or $\mathbb{P}(H^0(E', \mathcal{L}')^*)$ when $n$ is even. Take $\overline{z}$\index{zbar@$\overline{z}$} to be the image of the point $z \in E$ in $\overline{E}$.

Recall that $s = n /(n,2)$. By Proposition \ref{pcenterS}, the subalgebra $\kk[u_0,u_1,u_2,u_3,g_1,g_2]$ of $S$ is subject to degree $2s$ relations 
$$
r_i:=f_i(u_0,u_1,u_2,u_3)+h_i(g_1,g_2),
\index{r1i@$r_i$}
$$
for $i =1,2$, where $f_1,f_2$ are linearly independent quadratic forms in $u_0,u_1,u_2,u_3$ defining the elliptic curve $\overline{E} \subset \overline{\mathbb P}$.  (We get that $u_i = z_i$, $r_i = F_i$, and $f_i = \Phi_i$ in Proposition~\ref{pcenterS}(3), when $n$ is odd.) 
\enota

By \cite[Theorems~3.7 and~4.9]{SmithTate}, the center of the Veronese subalgebra $S^{(n/s)}$ of $S$ is 
\begin{equation} \label{eq:Veronese}
Z(S^{(n/s)})=\kk[u_0,u_1,u_2,u_3]^{(n)}[g_1,g_2], \quad \text{ subject to relations $r_1, r_2$}.
\end{equation}

\bnota{Q}[$Q(\overline{z})$, $\overline{E}_2$, $\overline{\omega}_i$, $\overline{e}_i$]
For $\overline{z}\in \overline{E}$, let 
$$
Q(\overline{z})=\bigcup \{\ell_{\overline{p},\overline{q}} ~|~ \overline{p},\overline{q}\in \overline{E}, ~\overline{p}+\overline{q}=\pm \overline{z}\} 
\index{Qzbar@$Q(\overline{z})$}
$$
denote the quadric in $\overline{\mathbb P}$ containing $\overline{E}$. By \cite[Section~3]{LevSmith} and \cite[paragraph before Theorem 5.9]{SmithTate}, we know $Q(\overline{z})$ is singular, 
if and only if $\overline{z}\in \overline{E}_2$, where $\overline{E}_2$\index{E2bar@$\overline{E}_2$} 
is the $2$-torsion subgroup of $\overline{E}$, and each singular quadric has  rank 3 and  has only one singular point.

Label the four 2-torsion points on $\overline{E}$ as $\overline{E}_2=\{\overline{\omega}_i\}_{0\le i\le 3}$.\index{omegabari@$\overline{\omega}_i$}   
For each corresponding singular quadric $Q(\overline{\omega}_i)$, denote its unique singularity by $\overline{e}_i$.\index{ebari@$\overline{e}_i$}
We choose representatives for all $\{\overline{e}_i\}_{0\le i\le 3}$ as points in the affine space ${\rm maxSpec}(\kk[u_0,u_1,u_2,u_3,g_1,g_2])$ as
\begin{gather}\label{SingEbar}
\mbox{\small{$\overline{e}_0=(1,0,0,0,0,0), \quad \overline{e}_1=(0,1,0,0,0,0), \quad  \overline{e}_2=(0,0,1,0,0,0),\quad \overline{e}_3=(0,0,0,1,0,0)$,}}
\end{gather}
which correspond to the four points $\{e_i\}_{0\le i\le 3}$ in the point scheme of $S$ in Definition-Lemma~\ref{dlgeomS}. 
\enota

This brings us to the first preliminary result, from Smith-Tate \cite{SmithTate}.

\ble{VS2} \textnormal{[$\pi$, $f_{\overline{z}}$, $h_{\overline{z}}$]} \cite[Lemma 5.7]{SmithTate}
There is a morphism 
$$\pi:\overline{E} \to \mathbb P(\kk r_1+\kk r_2)
\index{p3i@$\pi$}$$
such that $\pi(\overline{z})=f_{\overline{z}}+h_{\overline{z}}$, where $f_{\overline{z}}\in k[u_0,u_1,u_3,u_4]_{2s}$\index{fzbar@$f_{\overline{z}}$}
vanishes on the quadric $Q(\overline{z})$, and $h_{\overline{z}}\in \kk[g_1,g_2]_{2s}$\index{hzbar@$h_{\overline{z}}$} 
is a nonzero scalar multiple of 
\begin{align*}
\prod_{z\in E,\, \text{preimage of }\, \overline{z}} \Omega(z) ~= ~\prod_{i=0}^{s-1} \Omega(z+i(n/s)\tau).
\end{align*} 
The morphism $\pi$ is of degree $2$ and $\pi(\overline{p})=\pi(\overline{q})$, if and only if $\overline{p}=\pm \overline{q}$, for $\overline{p}, \overline{q} \in \overline{E}$. \qed
\ele

Finally, we highlight two central elements  of $S$ that will play a key role in the description of $Y^{sing}$ in Theorems~\ref{tSOdd} and~\ref{tSEven} below.

\bnota{central}[$g$, $G$]
Consider the elements of $\kk[g_1,g_2]$: 
$$g:= \prod_{\substack{\omega \in E_2,\\0 \leq k \leq s-2} }\Omega(\omega + k \tau)\quad \text{and}\quad G:=\prod_{\overline{\omega}\in \overline{E}_2} h_{\overline{\omega}}=g\prod_{\omega\in E_2}\Omega(\omega+(s-1)\tau),
\index{G00@$G$}\index{g00a@$g$}$$
which are both central in $S$ of degree $8(s-1)$ and $8s$, respectively (since $\deg(g_i) = 2$).
\enota

\ble{Ramification} 
Retain the notation of Lemma \ref{lVS2}. Let $\overline{p},\overline{q}\in \overline{E}$ such that 
$\overline{p}\neq \pm \overline{q}$.
Then, the following sets are the same in ${\rm Proj}(\kk[g_1,g_2])$: 
\begin{itemize}
\item[(a)] The zero locus $\mathbb{V}(g)$ of $g$; 
\item[(b)] The ramification locus of the map ${\rm Proj}(\kk[g_1,g_2]) \to {\rm Proj}(\kk[h_{\overline{p}},h_{\overline{q}}])$ induced by the natural embedding $\kk[h_{\overline{p}},h_{\overline{q}}]\hookrightarrow \kk[g_1,g_2];$ 
\item[(c)] The zero locus of the determinant of the Jacobian matrix $\frac{\partial(h_{\overline{p}},h_{\overline{q}})}{\partial(g_1,g_2)}$.
\end{itemize}  
\ele
\begin{proof}
By Lemma \ref{lVS2} and \eqref{eq:Veronese}, we know $\pi(\overline{p})$ and $\pi(\overline{q})$ can be realized as two defining relations of $Z(S^{(n/s)})$. Note that the ramification locus in (b) and the zero locus in (c) do not depend on the choice of $h_{\overline{p}}$ and $h_{\overline{q}}$ as long as they are linearly independent, namely as long as $\overline{p}\neq \pm \overline{q}$. 

(a)$\Leftrightarrow$(b) follows from \cite[Lemma 5.8(c)]{SmithTate} and the definition of $g$ in Notation~\ref{ncentral}. 

(b)$\Leftrightarrow$(c) This is a standard algebro-geometric fact; see \cite[Section~6.3]{Shaf}.
\end{proof}



\subsection{For PIdeg($S$) odd}
We assume that the PI degree $n$ of $S$ is odd in this section, so that $n=s$ and $u_i=z_i$ for all $0\le i\le 3$. Fix any $\overline{p}, \overline{q} \in \overline{E} = E/\langle \sigma \rangle$ with $\overline{p} \neq \pm \overline{q}$, and recall by Notation~\ref{nri}, \eqref{eq:Veronese}, and Lemma~\ref{lVS2} that
$$Y:=\text{maxSpec}(Z(S)) = \mathbb{V}(\pi(\overline{p}),\pi(\overline{q})).$$

We first establish that $Y$ is smooth outside of the variety $\mathbb{V}(G)$; see Notation~\ref{ncentral}.

\ble{SingSlice}
The subvariety $Y_{\gamma_1,\gamma_2}$ of $Y$ is smooth of dimension 2 if $(\gamma_1,\gamma_2)\not\in \mathbb{V}(G)$.
\ele
\begin{proof}
The defining ideal of $Y_{\gamma_1,\gamma_2}\subseteq \mathbb A^4_{(z_0,z_1,z_2,z_3)}$ is generated by  elements $f_{\overline{p}}+h_{\overline{p}}(\gamma_1,\gamma_2)$ and $f_{\overline{q}}+h_{\overline{q}}(\gamma_1,\gamma_2)$ in $\kk[z_0,z_1,z_2,z_3]$. Recall that $Y_{\gamma_1,\gamma_2}$ is nonsingular at a point $P\in Y_{\gamma_1,\gamma_2}$ if the rank of the Jacobian matrix $J$ evaluated at $P$ is $4-\dim (Y_{\gamma_1,\gamma_2})$; here,  
\[
J=\begin{bmatrix}\bigskip
\frac{\partial f_{\overline{p}}}{\partial z_0} & \frac{\partial f_{\overline{p}}}{\partial z_1} &\frac{\partial f_{\overline{p}}}{\partial z_2} & \frac{\partial f_{\overline{p}}}{\partial z_3}\\
\frac{\partial f_{\overline{q}}}{\partial z_0} & \frac{\partial f_{\overline{q}}}{\partial z_1} &\frac{\partial f_{\overline{q}}}{\partial z_2} & \frac{\partial f_{\overline{q}}}{\partial z_3}
\end{bmatrix}.
\]
Moreover, $Y_{\gamma_1,\gamma_2}$ is singular at a point $P$ if the rank of  $J$ evaluated at $P$ is less than $4-\dim (Y_{\gamma_1,\gamma_2})$. Hence, it suffices to show that $J$ evaluated at any $P\in Y_{\gamma_1,\gamma_2}$ has rank 2 if $(\gamma_1,\gamma_2)\not\in \mathbb{V}(G)$. 

Suppose that the rank of $J$ evaluated at some point $P\in Y_{\gamma_1,\gamma_2}$ is less than 2. We will then show that the pair $(\gamma_1,\gamma_2)$ lies in $\mathbb{V}(G)$. We can choose $\overline{p},\overline{q}\in \overline{E}$ so that $h_{\overline{p}}$ and $h_{\overline{q}}$ have no common factors (see Lemma \ref{lVS2} and \eqref{EOmega}). Hence $\mathbb{V}(h_{\overline{p}},h_{\overline{q}})=(0,0)$. 

Suppose that $P = (0,0,0,0)$, then  $f_{\overline{p}}(0,0,0,0)+h_{\overline{p}}(\gamma_1,\gamma_2)=h_{\overline{p}}(\gamma_1,\gamma_2)=0$ and $f_{\overline{q}}(0,0,0,0)+h_{\overline{q}}(\gamma_1,\gamma_2)=h_{\overline{q}}(\gamma_1,\gamma_2)=0$. So, $(\gamma_1,\gamma_2)\in \mathbb{V}(h_{\overline{p}},h_{\overline{q}})$. Hence, $(\gamma_1,\gamma_2) = (0,0) \in \mathbb{V}(G)$, and we are done.

Next, take $P \neq (0,0,0,0)$ and we will show that $P$ is the singularity of the following quadric containing $\overline{E}$ defined by 
$$f:=h_{\overline{q}}(\gamma_1,\gamma_2)f_{\overline{p}}-h_{\overline{p}}(\gamma_1,\gamma_2)f_{\overline{q}}=0.$$
Consider $\widetilde{P}=(P,\gamma_1,\gamma_2)\in Y$. We only need to check that $\frac{\partial f}{\partial z_i}(P)=0$ for all $0\le i\le 3$. We have that:
\[
\begin{array}{rl}
\smallskip

\textstyle \frac{\partial f}{\partial z_i}(P)&\,= \textstyle\left(h_{\overline{q}}(\gamma_1,\gamma_2)\frac{\partial f_{\overline{p}}}{\partial z_i}-h_{\overline{p}}(\gamma_1,\gamma_2)\frac{\partial f_{\overline{q}}}{\partial z_i}\right)(P)\\

\smallskip

 &=\left(h_{\overline{q}}\frac{\partial f_{\overline{p}}}{\partial z_i}-h_{\overline{p}}\frac{\partial f_{\overline{q}}}{\partial z_i}\right)(\widetilde{P})\\
 
 \smallskip
 
&\,= \textstyle \left(-f_{\overline{q}}\frac{\partial f_{\overline{p}}}{\partial z_i}+f_{\overline{p}}\frac{\partial f_{\overline{q}}}{\partial z_i}\right)(P) \\

\smallskip

&=2\left(-\sum_{0\le j\le 3}z_j\frac{\partial f_{\overline{q}}}{\partial z_j}\frac{\partial f_{\overline{p}}}{\partial z_i}+\sum_{0\le j\le 3}z_j\frac{\partial f_{\overline{p}}}{\partial z_j}\frac{\partial f_{\overline{q}}}{\partial z_i}\right)(P)\\

\smallskip

&\,= \textstyle 2\sum_{0\le j\le 3}z_j\left(\frac{\partial f_{\overline{p}}}{\partial z_j}\frac{\partial f_{\overline{q}}}{\partial z_i}-\frac{\partial f_{\overline{q}}}{\partial z_j}\frac{\partial f_{\overline{p}}}{\partial z_i}\right)(P)\\

\smallskip

 &=0.
\end{array}
\]
We use the identity $f_{\overline{p}}=2\sum_{0\le j\le 3}z_j\frac{\partial f_{\overline{p}}}{\partial z_j}$ since $f_{\overline{p}}$ is homogenous of degree two, and the similar identity for $f_{\overline{q}}$. Also, the last equality above comes from the Jacobian matrix $J$ has rank $\le 1$ at $P$. Therefore we can set $P=\lambda \overline{e}_k$ for some $\lambda\in \kk^\times$, where $\overline{e}_k$ is the singularity of the corresponding quadric $Q(\overline{\omega}_k)$ containing $\overline{E}$ according to Notation \ref{nQ}. Since $\widetilde{P}\in Y$, we then get by Lemma \ref{lVS2} that
\[
\begin{array}{rll}

\smallskip

0&=\pi (\overline{\omega}_k)(\widetilde{P}) ~=~ (f_{\overline{\omega}_k}+h_{\overline{\omega}_k})(\widetilde{P}) &=f_{\overline{\omega}_k}(P)+h_{\overline{\omega}_k}(\gamma_1,\gamma_2)\\

&=\lambda^2f_{\overline{\omega}_k}(\overline{e}_k)+h_{\overline{\omega}_k}(\gamma_1,\gamma_2) &=h_{\overline{\omega}_k}(\gamma_1,\gamma_2).
\end{array}
\]
By Notation \ref{ncentral}, we conclude that $(\gamma_1,\gamma_2)\in \mathbb{V}(h_{\overline{\omega}_k})\subset \mathbb{V}(G)$, as desired. 
\end{proof}

Next, we show that the union of singular loci of the subvarieties $Y_{\gamma_1,\gamma_2}$ for certain points $(\gamma_1, \gamma_2) \in \mathbb{V}(G)$ is equal to a cuspidal curve, denoted by $C(\omega+k\tau)$ below.

\ble{Wcusp}\textnormal{[$C(\omega+k\tau)$]}
For any $\omega\in E_2$ and $0\le k\le n-1$, there exists a nonzero point $\overline{p}_{\omega, k}\in \mathbb{V}(\Omega(\omega+k\tau))\cap \mathbb{V}(z_0,z_1,z_2,z_3)$ such that 
$$\bigcup_{(\gamma_1,\gamma_2)\in \mathbb V(\Omega(\omega+k\tau))}(Y_{\gamma_1,\gamma_2})^{sing}~=~\{t^n\overline{e}_\omega+t^2\overline{p}_{\omega, k}\,|\, t\in \kk\} ~=: ~C(\omega+k\tau),
\index{C0omegak@$C(\omega+k\tau)$}
$$
where $\overline{\omega}\in \overline{E}_2$ is the image of $\omega$ under the isogeny $E\twoheadrightarrow \overline{E}=E/\langle \sigma\rangle$, and $\overline{e}_\omega$ is the singularity of the quadric $Q(\overline{\omega})$. Moreover, we can take $\overline{p}_{\omega, k}=\overline{p}_{\omega, n-2-k}$ with the second index modulo $n$.
\ele

\begin{proof}
We can choose the generators of the defining ideal of $Y=\mathbb{V}(\pi(\overline{p}), \pi(\overline{q}))$, where $\overline{p}=\overline{\omega}\in \overline{E}_2$ and $\overline{q}$ generic, so that $\pi(\overline{p})=f_{\overline{\omega}}+h_{\overline{\omega}}$ with $h_{\overline{\omega}}=\prod_{k=0}^{n-1}\Omega(\omega+k\tau)$, and $\pi(\overline{q})=f_{\overline{q}}+h_{\overline{q}}$. Let $P \in (Y_{\gamma_1,\gamma_2})^{sing}$ and take $\widetilde{P}=(P,\gamma_1,\gamma_2)\in Y$. Since $(\gamma_1,\gamma_2)$ is in $\mathbb{V}(\Omega(\omega+k\tau))\subset \mathbb{V}(h_{\overline{\omega}})$, we have that
$$0=\pi(\overline{\omega})(\widetilde{P})=f_{\overline{\omega}}(P)+h_{\overline{\omega}}(\gamma_1,\gamma_2)=f_{\overline{\omega}}(P).$$ 
By a similar argument in Lemma \ref{lSingSlice}, one can further show that $P$ is the singularity of the quadric $Q(\overline{\omega})=\mathbb{V}(f_{\overline{\omega}})$ containing $\overline{E}$. Hence $P=\lambda \overline{e}_\omega$ for some $\lambda\in \kk$. Pick some nonzero point $\overline{p}_{\omega, k}\in \mathbb{V}(\Omega(\omega+k\tau))\cap \mathbb{V}(z_0,z_1,z_2,z_3)$. Since $\overline{q}\in \overline{E}$ is generic, by a possible rescaling of $\overline{p}_{\omega, k}$, we can assume that $f_{\overline{q}}(\overline{e}_\omega)=-h_{\overline{q}}(\overline{p}_{\omega, k})\neq 0$. Then $\widetilde{P}$ is contained in the following subset of $Y$:
\[
\begin{array}{rl}
\smallskip

 \{\alpha \overline{e}_\omega+\beta \overline{p}_{\omega, k}\, |\, \alpha,\beta\in \kk\}\cap Y &=\{\alpha \overline{e}_\omega+\beta \overline{p}_{\omega, k}\, |\, \alpha,\beta\in \kk\}\cap \mathbb{V}(f_{\overline{q}}+h_{\overline{q}})\\
 
 \smallskip
 
&= \{\alpha \overline{e}_\omega+\beta \overline{p}_{\omega, k}\, |\, \alpha^2f_{\overline{q}}(\overline{e}_\omega)+\beta^nh_{\overline{q}}(\overline{p}_{\omega, k})=0\}\\

\smallskip

&=\{ t^n\overline{e}_\omega+t^2\overline{p}_{\omega, k}\,|\, t\in \kk\}.
\end{array}
\]
Therefore, we have obtained that $\bigcup_{(\gamma_1,\gamma_2) \in \mathbb{V}(\Omega(\omega + k\tau))} (Y_{\gamma_1,\gamma_2})^{sing} \subset \{ t^n\overline{e}_\omega+t^2\overline{p}_{\omega, k}\,|\, t\in \kk\}$.

Conversely, take $P=t^n\overline{e}_\omega$ and $\widetilde{P}=t^n\overline{e}_\omega+t^2\overline{p}_{\omega,k}$ for any $t\in \kk$. It is straight-forward to check that $P\in (Y_{\gamma_1,\gamma_2})^{sing}$ with $(\gamma_1,\gamma_2)=t^2\overline{p}_{\omega, k}$ since the the first row of the Jacobian matrix vanishes at $P$:
\[
\begin{bmatrix}\bigskip
\frac{\partial f_{\overline{\omega}}}{\partial z_0} & \frac{\partial f_{\overline{\omega}}}{\partial z_1} &\frac{\partial f_{\overline{\omega}}}{\partial z_2} & \frac{\partial f_{\overline{\omega}}}{\partial z_3}\\
\frac{\partial f_{\overline{q}}}{\partial z_0} & \frac{\partial f_{\overline{q}}}{\partial z_1} &\frac{\partial f_{\overline{q}}}{\partial z_2} & \frac{\partial f_{\overline{q}}}{\partial z_3}
\end{bmatrix}.
\]

Finally, we can set $\overline{p}_{\omega, k}=\overline{p}_{\omega, n-2-k}$ since $\Omega(\omega+k\tau)=\Omega(\omega+(n-2-k)\tau)$ by \eqref{EOmega}.
\end{proof}

Now we show that $Y^{sing} \subset \mathbb{V}(g)$ and we can use the lemma above to get a detailed description of $Y^{sing}$.

\bth{SOdd}
For $n$ odd, we have
$$Y^{sing}~=~Y^{sing}\cap \mathbb{V}(g)~=\bigcup_{\substack{\omega\in E_2\\0\le k\le n-2} } C(\omega+k\tau),$$
where $C(\omega+k\tau)$ is the cuspidal curve $\left\{ t^n\overline{e}_\omega+t^2\overline{p}_{\omega k}\, |\, t\in \kk \right\}$ defined in Lemma \ref{lWcusp}. As a consequence, $Y^{sing}$ is a union of $2(n-1)$ cuspidal curves in $\mathbb{V}(g)$ meeting at the origin as depicted in Figure~1. 

We also have the following results on the subvarieties $(Y_{\gamma_1,\gamma_2})^{sing}$, ~ $(Y^{sing})_{\gamma_1,\gamma_2}$ of $Y$:
\smallskip

\begin{enumerate}
\item If $(\gamma_1,\gamma_2)\not \in \mathbb{V}(G)$, then $(Y_{\gamma_1,\gamma_2})^{sing}=(Y^{sing})_{\gamma_1,\gamma_2}=\varnothing$.
\smallskip

\item If $(\gamma_1,\gamma_2)\in \mathbb{V}(g)$, then $(Y_{\gamma_1,\gamma_2})^{sing}=(Y^{sing})_{\gamma_1,\gamma_2}$, which consists of 1 point if $(\gamma_1,\gamma_2)=(0,0)$, and consists of 2 points otherwise.
\smallskip

\item If $(0,0) \neq (\gamma_1,\gamma_2)\in \mathbb{V}(G/g)$, then $(Y_{\gamma_1,\gamma_2})^{sing}$ has 2 points while $(Y^{sing})_{\gamma_1,\gamma_2}$ is empty. 
\smallskip

\item The group $H_4$ in Proposition~\ref{pH4} fixes $Y^{sing}$ and $Y^{symp}_0$ and
\begin{align*}
\epsilon_1\,&: C(\omega+k\tau)\mapsto C(\omega+\xi_1+k\tau),&\epsilon_2&\,: C(\omega+k\tau)\mapsto C(\omega+\xi_2+k\tau),\\
\epsilon_1\epsilon_2\,&: C(\omega+k\tau)\mapsto C(\omega+\xi_3+k\tau),&\epsilon&\,: C(\omega+k\tau)\mapsto C(\omega+k\tau),
\end{align*}
where $\{0,\xi_1,\xi_2,\xi_3\}$ is the 2-torsion subgroup of $E$ as described in \cite[Section~2.6]{ChirSmith}.
\end{enumerate}
\smallskip

As a consequence, $(Y_{\gamma_1,\gamma_2})^{sing}$ consists of only finitely many points, and  
$$Y^{symp}_0 = \bigcup_{\omega \in E_2} C(\omega+(n-1)\tau) \cup Y^{sing}.$$
\eth

 \bigskip
 
 \usetikzlibrary{positioning}
\usetikzlibrary{decorations.pathreplacing}

\begin{flushright}
\begin{tikzpicture}[scale=0.9][every node/.style={minimum size=1cm},on grid]
\begin{scope}[every node/.append style={yslant=-0.5},yslant=-0.5]
  \shade[right color=gray!80, left color=black!10] (1,0) rectangle +(4,4);
  \node at (2.7,0.3) {\footnotesize{$\mathbb{V}(\Omega(\omega+(n-1)\tau))$}};
    \node at (2.5,3.3) {\tiny{\textcolor{rgb:red,1;blue,1}{$C(\omega+(n-1)\tau)$}}};
  \draw[<->][shift={(5,2)}, scale=0.01, draw opacity=.9, line width=\pgflinewidth+.7pt, smooth,variable=\t, color={rgb:red,1;blue,1}] plot ({-15.5*\t*\t},{\t*\t*\t});
\end{scope}
\begin{scope}[every node/.append style={yslant=-0.2, xslant=0},yslant=-0.2, xslant=0]
  \shade[right color=gray!10,left color=gray!80] (5,-1.5) rectangle +(6,4);
  \node at (10.0,-1.2) {\footnotesize{\textcolor{rgb:green,1;blue,1}{$\mathbb{V}(\Omega(\omega))$}}};
\draw[<->][shift={(5.7,0.6)}, scale=0.011, draw opacity=.9, line width=\pgflinewidth+.7pt, smooth,variable=\t, color={rgb:red,1;blue,1}] plot ({18*\t*\t},{\t*\t*\t});
\end{scope}
\begin{scope}[every node/.append style={yslant=0.0, xslant=0},yslant=0.0, xslant=0]
  \shade[right color=gray!40,left color=gray!90] (5,-2.5) rectangle +(5.7,4);
\draw [gray,decorate,decoration={brace,amplitude=5pt},
   xshift=-4pt,yshift=-9pt] (1,4.3)  -- (11,4.3) 
   node [black,midway,above=4pt,xshift=-2pt] {\small $\mathbb{V}(G) \subset \mathbb{A}^6_{(z_0, \dots, z_3, g_1,g_2)}$};
  \draw [gray,decorate,decoration={brace,amplitude=5pt},
   xshift=-4pt,yshift=-9pt] (10,4.0)  -- (12,0) 
   node [black,midway,right=4pt,xshift=-2pt] {\small $\begin{array}{c} \hspace{-.2in} \textcolor{rgb:green,1;blue,1}{\mathbb{V}(g)}\\ \hspace{-.1in}\cap\\ \mathbb{V}(G)\end{array}$}; 
\end{scope}
\begin{scope}[every node/.append style={yslant=0.1, xslant=0},yslant=0.1, xslant=0]
  \shade[right color=gray!30,left color=gray!90] (5,-3.0) rectangle +(5.5,4);
\end{scope}
\begin{scope}[every node/.append style={yslant=0.2, xslant=0},yslant=0.2, xslant=0]
  \shade[right color=gray!10,left color=gray!80] (5,-3.5) rectangle +(5.2,4);
  \node at (8.5,-3.2) {\footnotesize{\textcolor{rgb:green,1;blue,1}{$\mathbb{V}(\Omega(\omega+(n-3)\tau))$}}};
      \draw[<->][shift={(5.7,-1.5)}, scale=0.0102, draw opacity=.9, line width=\pgflinewidth+.7pt, smooth,variable=\t, color={rgb:red,1;blue,1}] plot ({17*\t*\t},{\t*\t*\t});
\end{scope}
\begin{scope}[every node/.append style={yslant=0.5},yslant=0.5]
  \shade[right color=gray!10,left color=gray!60] (5,-5) rectangle +(4.2,4);
  \node at (7.5,-4.7) {\footnotesize{\textcolor{rgb:green,1;blue,1}{$\mathbb{V}(\Omega(\omega+(n-2)\tau))$}}};
  \node at (7.7,-1.7) {\tiny{\textcolor{rgb:red,1;blue,1}{$C(\omega+(n-2)\tau)$}}};
  \draw[->, dotted]  (3.4,-0.3) -- (2.9,-0.3);
    \draw[<->][shift={(5,-3)}, scale=0.01, draw opacity=.9, line width=\pgflinewidth+.7pt, smooth,variable=\t, color={rgb:red,1;blue,1}] plot ({16*\t*\t},{\t*\t*\t});
\end{scope}
  \draw[->, dotted]  (6.8,1.45) -- (7.1,1.1);
\end{tikzpicture}
\end{flushright}

\vspace{-1.3in}

\begin{flushleft}
\begin{tikzpicture}[scale=0.3][every node/.style={minimum size=1cm},on grid]
\draw[dashed] (5,5) ellipse (6.5cm and 7cm);
\begin{scope}[every node/.append style={yslant=0},yslant=1.0]
  \shade[right color=gray!80, left color=black!10] (2,-1) rectangle +(3,3);
\node at (5,-3) {\footnotesize{$\mathbb{V}(G)$}};  
\end{scope}
\begin{scope}[every node/.append style={yslant=-0.4, xslant=0},yslant=-0.4, xslant=0]
  \shade[right color=gray!50,left color=gray!100] (5,6) rectangle +(4,3);
\end{scope}
\begin{scope}[every node/.append style={yslant=-0.6, xslant=0},yslant=-0.6, xslant=0]
  \shade[right color=gray!40,left color=gray!90] (5,7) rectangle +(3.6,3);
\end{scope}
\begin{scope}[every node/.append style={yslant=-0.8, xslant=0},yslant=-0.8, xslant=0]
  \shade[right color=gray!30,left color=gray!80] (5,8) rectangle +(3.4,3);
\end{scope}
\begin{scope}[every node/.append style={yslant=-1.0, xslant=0},yslant=-1.0, xslant=0]
  \shade[right color=gray!20,left color=gray!70] (5,9) rectangle +(3.2,3);
\end{scope}
\begin{scope}[every node/.append style={yslant=0},yslant=-1.2]
  \shade[right color=gray!10,left color=gray!60] (5,10) rectangle +(3,3);
\end{scope}
   \node at (.8,1.0) {\footnotesize{$g_1$}};
 \draw[->, dotted] (5,5.2) -- (1.2,1.4);
   \node at (10.5,5.2) {\footnotesize{$g_2$}};
 \draw[->, dotted] (5,5.2) -- (10,5.2);
   \node at (5,9.8) {\footnotesize{$\mathbb{A}^4_{(z_0,\dots,z_3)}$}};
 \draw[->, dotted] (5,5.2) -- (5,9.2);
\end{tikzpicture}
\end{flushleft}

\vspace{.3in}

\begin{center}
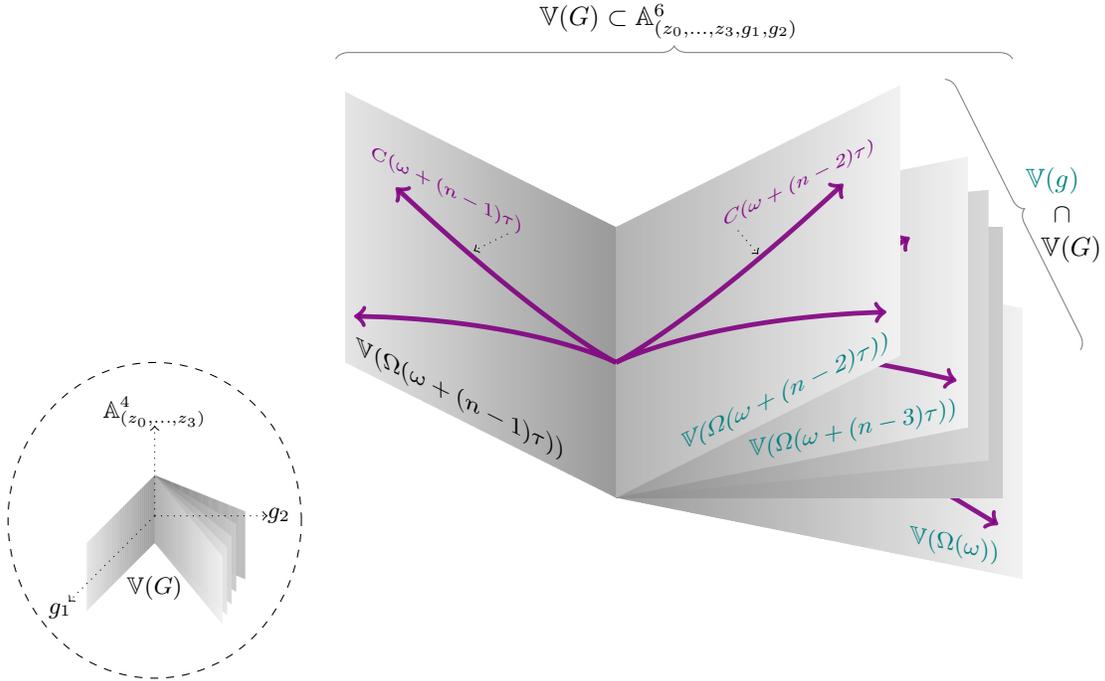
{\textsc{Figure 1.} For $n$ odd: $Y^{sing} =   \bigcup_{\gamma_1, \gamma_2 \in \kk} (Y^{sing})_{\gamma_1,\gamma_2} = \bigcup_{\omega\in E_2, ~0\le k\le n-2 }C(\omega+k\tau)$,\\
\smallskip
\hspace{1.1in} and  $Y^{symp}_0 =   \bigcup_{\gamma_1, \gamma_2 \in \kk} (Y_{\gamma_1,\gamma_2})^{sing} = \bigcup_{\omega\in E_2, ~0\le k\le n-1 }C(\omega+k\tau)$, \\
\smallskip
 union of $2n-2$ (resp., $2n$) cuspidal curves meeting at $\{\underline{0}\}$, each with multiplicity 2}
\end{center}

\vspace{.4in}

\noindent{\it Proof of Theorem~\ref{tSOdd}}.
The Jacobian matrix of $Y=\mathbb{V}(\pi(\overline{p}),\pi(\overline{q}))$ is 
\[
\begin{bmatrix}\bigskip
\frac{\partial f_{\overline{p}}}{\partial z_0} & \frac{\partial f_{\overline{p}}}{\partial z_1} &\frac{\partial f_{\overline{p}}}{\partial z_2} & \frac{\partial f_{\overline{p}}}{\partial z_3}& \frac{\partial h_{\overline{p}}}{\partial g_1} &\frac{\partial h_{\overline{p}}}{\partial g_2}\\
\frac{\partial f_{\overline{q}}}{\partial z_0} & \frac{\partial f_{\overline{q}}}{\partial z_1} &\frac{\partial f_{\overline{q}}}{\partial z_2} & \frac{\partial f_{\overline{q}}}{\partial z_3} & \frac{\partial h_{\overline{q}}}{\partial g_1} &\frac{\partial h_{\overline{q}}}{\partial g_2}
\end{bmatrix}.
\]
It is clear that if it has rank $\leq 1$ at some point $\widetilde{P}=(P,\gamma_1,\gamma_2)\in Y$, then the Jacobian matrix of $Y_{\gamma_1,\gamma_2}$, which is one of its $2\times 4$ minors containing all the partial derivatives of $z_i$, also has rank $\leq 1$ at $P$. So 
\begin{equation} \label{eq:contains}
(Y^{sing})_{\gamma_1,\gamma_2}\subseteq (Y_{\gamma_1,\gamma_2})^{sing}.
\end{equation} 

Now if $\widetilde{P}\in Y^{sing}$ then each $2\times 2$ minors of the matrix above  vanishes at $\widetilde{P}$. In particular, we have that
\[
{\rm det} \left(\begin{matrix}
\bigskip
\frac{\partial h_{\overline{p}}}{\partial g_1} &\frac{\partial h_{\overline{p}}}{\partial g_2}\\
\frac{\partial h_{\overline{q}}}{\partial g_1} &\frac{\partial h_{\overline{q}}}{\partial g_2}
\end{matrix}\right)(\widetilde{P})={\rm det}\left(\frac{\partial(h_{\overline{p}},h_{\overline{q}})}{\partial(g_1,g_2)}\right)(\gamma_1,\gamma_2)=0.
\]
Hence $Y^{sing}\subset \mathbb{V}(g)$ by Lemma \ref{lRamification} and 
$$Y^{sing}~=\bigcup_{(\gamma_1,\gamma_2)\in\mathbb{V}(g)}(Y^{sing})_{\gamma_1,\gamma_2} ~\subseteq \bigcup_{(\gamma_1,\gamma_2)\in \mathbb{V}(g)}(Y_{\gamma_1,\gamma_2})^{sing}~= \bigcup_{\substack{\omega\in E_2\\0\le k\le n-2} }C(\omega+k\tau)$$
by Lemma \ref{lWcusp}. 

Conversely, we verify that $C(\omega+k\tau) \subset Y^{sing}$ for all $\omega\in E_2$ and $0\le k\le n-2$ as follows. Let $\pi(\overline{\omega})$ and $\pi(\overline{q})$ be the two defining relations of $Y$ with $q \in E$ generic. It suffices to show that the first row of the Jacobian matrix of $Y$ above vanishes, i.e., we want to show that
\begin{align*}\small
 &\,\textstyle \left[\frac{\partial f_{\overline{\omega}}}{\partial z_0},\frac{\partial f_{\overline{\omega}}}{\partial z_1}, \frac{\partial f_{\overline{\omega}}}{\partial z_2}, \frac{\partial f_{\overline{\omega}}}{\partial z_3},\frac{\partial h_{\overline{\omega}}}{\partial g_1}, \frac{\partial h_{\overline{\omega}}}{\partial g_2}\right](t^n\overline{e}_\omega+t^2\overline{p}_{\omega, k})\\
&\, = \textstyle t^{2n}\left[\frac{\partial f_{\overline{\omega}}}{\partial z_0}(\overline{e}_\omega),\frac{\partial f_{\overline{\omega}}}{\partial z_1}(\overline{e}_\omega), \frac{\partial f_{\overline{\omega}}}{\partial z_2}(\overline{e}_\omega), \frac{\partial f_{\overline{\omega}}}{\partial z_3}(\overline{e}_\omega),0,0\right]+t^{2n}\left[0,0,0,0,\frac{\partial h_{\overline{\omega}}}{\partial g_1}(\overline{p}_{\omega k}), \frac{\partial h_{\overline{\omega}}}{\partial g_2}(\overline{p}_{\omega k})\right]
\end{align*} 
vanishes. Since $\overline{e}_\omega$ is the singularity of $Q(\overline{\omega})$, the first summand vanishes. Note that when $0\le k\le n-2$ we have $\overline{p}_{\omega, k}\in \Omega(\omega+k\tau)=\Omega(\omega+(n-2-k)\tau)$ is a multiple root of $h_{\overline{\omega}}=\prod_{k=0}^{n-1}\Omega(\omega+k\tau)$. So, the second summand vanishes as well. Thus, the first part of the result holds. 

\medskip

(1) This follows from Lemma \ref{lSingSlice}.

\medskip

(2) Suppose $(\gamma_1,\gamma_2)\in \mathbb{V}(g)$. By Lemma \ref{lWcusp} and the beginning of the result, we have that
$$(Y^{sing})_{\gamma_1,\gamma_2}=(Y_{\gamma_1,\gamma_2})^{sing}=C(\omega+k\tau) \cap \mathbb{V}(g_1-\gamma_1,g_2-\gamma_2)$$
for some $(\gamma_1,\gamma_2)\in \mathbb V(\Omega(\omega+k\tau))$ with $0\leq k \leq n-2$. If $(\gamma_1,\gamma_2)=(0,0)$, then $(Y^{sing})_{\gamma_1,\gamma_2}=(Y_{\gamma_1,\gamma_2})^{sing}=(0,0,0,0)$. If $(\gamma_1,\gamma_2)\neq (0,0)$, then since $n$ is odd there are two choices of $t$ satisfying the defining equation of $C(\omega+k\tau)$ yielding two different points in $(Y^{sing})_{\gamma_1,\gamma_2}=(Y_{\gamma_1,\gamma_2})^{sing}$. 

\medskip

(3) The argument for $(Y_{\gamma_1,\gamma_2})^{sing}$ is similar as in (2) noting that 
$$\bigcup_{(\gamma_1,\gamma_2)\in \mathbb V(\Omega(\omega+(n-1)\tau))}(Y_{\gamma_1,\gamma_2})^{sing}~=~\{t^n\overline{e}_\omega+t^2\overline{p}_{\omega, n-1}\,|\, t\in \kk\}$$
for any $\omega\in E_2$ by Lemma \ref{lWcusp}. 

For $(Y^{sing})_{\gamma_1,\gamma_2}$, recall from  Notation \ref{ncentral} that $G/g=\prod_{\omega\in E_2} \Omega(\omega+(n-1)\tau)$. Since $g=\prod_{\omega\in E_2, 0\le k\le n-2} \Omega(\omega+k\tau)$, we can see that $G/g$ and $g$ as polynomials in $\kk[g_1,g_2]$ have no common factors by using the only non-trivial identity \eqref{EOmega} between these central annihilators $\Omega(\omega+k\tau)$. By the argument in the beginning of the proof of this theorem, we know $Y^{sing}\subseteq \mathbb{V}(g)$. This implies that 
$$(Y^{sing})_{\gamma_1,\gamma_2}~=~Y^{sing}\cap \mathbb{V}(g_1-\gamma_1,g_2-\gamma_2)~\subset~ \mathbb{V}(g)\cap \mathbb{V}(g_1-\gamma_1,g_2-\gamma_2).$$
By the assumption on $(\gamma_1,\gamma_2)$, we get $\mathbb{V}(g)\cap \mathbb{V}(g_1-\gamma_1,g_2-\gamma_2)=\varnothing$ as gcd$(G/g,g)=1$ in $\kk[g_1,g_2]$. Therefore, we have $(Y^{sing})_{\gamma_1,\gamma_2}=\varnothing$ if $(0,0)\neq (\gamma_1,\gamma_2)\in \mathbb{V}(G/g)$.

\medskip

(4) By Corollary~\ref{cfixZ}(4), the group $H_4\subset {\rm Aut}_{gr}(S)$ fixes $Y^{sing}$ and $Y^{symp}_0$. By \cite[Corollary 2.10]{ChirSmith}, we know $\epsilon_i(\Omega(\omega+k\tau))=\Omega(\omega+\xi_i+k\tau)$ for $i=1,2$. Note that $C(\omega+k\tau)\subset \mathbb V(\Omega(\omega+k\tau))$. So we have $\epsilon_i(C(\omega+k\tau))$ and $C(\omega+\xi_i+k\tau)$ are included in $\mathbb V(\Omega(\omega+\xi_i+k\tau))$. By Lemma \ref{lWcusp}, $\mathbb V(\Omega(\omega+\xi_i+k\tau))\cap Y^{symp}_0$ contains only one cuspidal curve, namely $\{t^n\overline{e}_{\omega+\xi_i}+t^2\overline{p}_{\omega+\xi_i, k}\,|\, t\in \kk\}$. So all three cuspidal curves are the same. Finally, $\epsilon$ is just a rescaling of the variables of $S$. So it will fix all the cuspidal curves $C(\omega+k\tau)$. This proves part (4).

\medskip

Finally, we have
\[
\begin{array}{rl}
Y^{symp}_0&=~Y^{symp}_0\cap \mathbb V(G)\\
&=~(Y^{symp}_0\cap \mathbb V(G/g))\cup (Y^{symp}_0\cap \mathbb V(g))\\
&=~ \bigcup_{\omega \in E_2} C(\omega+(n-1)\tau) \cup Y^{sing}.
\end{array}
\]
So, the result follows.

\qed 

\subsection{For PIdeg($S$) even} 
We now assume in this part that the PI degree $n$ of $S$ is even; here, $s=n/2$ and $u_i=z_i^2$ with ${\rm deg}\,u_i=s$ for all $0\le i\le 3$. Recall that the center of the Veronese subalgebra $S^{(2)}$ is 
$$Z(S^{(2)})=\kk[u_0,u_1,u_2,u_3]^{(n)}[g_1,g_2]$$ 
subject to two defining relations $\pi(\overline{p})=f_{\overline{p}}+h_{\overline{p}}$ and $\pi(\overline{q})=f_{\overline{q}}+h_{\overline{q}}$, where the points $\overline{p},\overline{q}\in \overline{E}=E/\langle \sigma^2\rangle$ satisfy $\overline{p}\neq \pm \overline{q}$. Note that $f_{\overline{p}}$ and $f_{\overline{q}}$ are two linearly independent quadrics in terms of $u_0,u_1,u_2,u_3$ defining the elliptic curve $\overline{E}\subset \mathbb \overline{\mathbb P}$. The center $Z$ of $S$ isomorphic to $\kk[z_0,z_1,z_2,z_3,g_1,g_2]$ subject to two defining relations $F_1,F_2$ of degree $2n$. Moreover, we can write the two defining relations of $Z$ as
$$F_1=(a_1-h_1)^2-z_iz_j,\quad F_2=(a_2-h_2)^2-z_kz_l,\quad \{i,j,k,l\}=\{0,1,2,3\},$$
where $h_1,h_2\in \kk[g_1,g_2]_n$ have no common factors and $a_1,a_2$ are linear forms in terms of $z_0,z_1,z_2,z_3$ given in Proposition~\ref{pZEven}.

\ble{SSliceE}
For each pair $(\gamma_1,\gamma_2) \in \kk^2$, we have that $\left(Y_{\gamma_1,\gamma_2}\right)^{sing}\subset \mathbb{V}(z_0z_1z_2z_3)$.
\ele

\begin{proof}
By way of contradiction, suppose there exists a point $$\widetilde{P} = (P,\gamma_1, \gamma_2) \in (Y_{\gamma_1, \gamma_2})^{sing}~ \backslash ~\mathbb{V}(z_0z_1z_2z_3)$$ for some $\gamma_1, \gamma_2 \in \kk$. Let $\mathfrak{m}_{\widetilde{P}}$ be the maximal ideal of $Z(S)$ corresponding to $\widetilde{P}$. Then by the Lying Over Theorem (for the integral extension $Z(S)\subset Z(S^{(2)})$), there exists a maximal ideal $\mathfrak{n}$ of $Z(S^{(2)})$ so that $\mathfrak{n} \cap Z(S) = \mathfrak{m}_{\widetilde{P}}$. Moreover, $\mathbb{V}(\mathfrak{n})$ is a point $\widetilde{Q} = (Q, \gamma_1, \gamma_2)$ of maxSpec($Z(S^{(2)})$).

Using the fact that $z_i = u_i^2$ for all $0\le i\le 3$, we get $Q\not\in \mathbb{V}(u_0u_1u_2u_3)$ since $P$ is not in $\mathbb{V}(z_0z_1z_2z_3)$. In particular, we have $Q \neq \underline{0}$. Moreover, since $P\in (Y_{\gamma_1,\gamma_2})^{sing}$, we  have that $Q$ is a singular point of ${\rm maxSpec}(Z(S^{(2)}))\cap \mathbb V(g_1-\gamma_1,g_2-\gamma_2)$ via the Implicit Function Theorem. Note that the defining relations of $Z(S^{(2)})$ are of the same form as the defining relations of $Z(S)$ in the case when $n$ is odd, we obtain by the argument in the proof of Lemma~\ref{lSingSlice} that the point $Q$ is the singular point of some quadric containing $\overline{E}$, namely $\overline{e}_k$ for some $k=0, \dots, 3$. This contradicts $Q\not\in \mathbb{V}(u_0u_1u_2u_3)$; see~\eqref{SingEbar}.
\end{proof}

\bth{SEven}\textnormal{[$Y_1^{sing}$, $Y_2^{sing}$]} 
Take $n$ even and recall the notation of Proposition~\ref{pZEven}. Then, we obtain that 
the singular locus $Y^{sing}$ of $Y = \mathbb{V}(F_1, F_2) \subset \mathbb{A}^6_{(z_0, z_1,z_2,z_3,g_1,g_2)}$ 
is the union of subvarieties $Y_1^{sing}$ and $Y_2^{sing}$ of $Y$ defined by
\begin{itemize}
\smallskip

\item[(i)] $Y_1^{sing}=\mathbb{V}(a_1-h_1, z_0, z_3)\cap Y$ and $Y_2^{sing}=\mathbb{V}(a_2-h_2, z_1, z_2)\cap Y$, when $\rho = \rho_1$,
\smallskip

\item[(ii)] $Y_1^{sing}=\mathbb{V}(a_1-h_1, z_0, z_2)\cap Y$ and $Y_2^{sing}=\mathbb{V}(a_2-h_2, z_1, z_3)\cap Y$, when $\rho = \rho_2$,
\smallskip

\item[(iii)] $Y_1^{sing}=\mathbb{V}(a_1-h_1, z_0, z_1)\cap Y$ and $Y_2^{sing}=\mathbb{V}(a_2-h_2, z_2, z_3)\cap Y$, when $\rho = \rho_3$.
\end{itemize}
Moreover, we have $Y_1^{sing}~\cap~ Y_2^{sing}=\{\underline{0}\}$. 

We also have the following results on the subvarieties $(Y_{\gamma_1,\gamma_2})^{sing}$,~ $(Y^{sing})_{\gamma_1,\gamma_2}$ of $Y$:
\begin{enumerate}
\item The varieties $Y_1^{sing}$ and $Y^{sing}_2$ from Theorem~\ref{tSEven} are
\[
\text{permuted by }
\begin{cases}
\epsilon_1 \text{ and } \epsilon_2,\\
\epsilon_1 \text{ and } \epsilon_1 \epsilon_2,\\
\epsilon_2 \text{ and } \epsilon_1 \epsilon_2, 
\end{cases}
\text{ and fixed by }
\begin{cases}
\epsilon_1 \epsilon_2, & \text{ when } \rho = \rho_1,\\
\epsilon_2, & \text{ when } \rho = \rho_2,\\
 \epsilon_1, & \text{ when } \rho = \rho_3.
\end{cases}
\]
where $\epsilon_i$  are the group actions in Proposition~\ref{pH4}.
\item $(Y^{sing})_{\gamma_1,\gamma_2}$ has 4 points generically (counting multiplicity), and 1 point,
if and only if $(\gamma_1,\gamma_2)=(0,0)$.
\item  $(Y_{\gamma_1,\gamma_2})^{sing}=(Y^{sing})_{\gamma_1,\gamma_2} = (Y_1)^{sing}_{\gamma_1,\gamma_2}\,\cup\, (Y_2)^{sing}_{\gamma_1,\gamma_2}$, where $$(Y_i)^{sing}_{\gamma_1,\gamma_2} := Y_i^{sing} \cap \mathbb{V}(g_1 - \gamma_1, g_2-\gamma_2)$$ for $i =1,2$.
\end{enumerate}
\eth

\vspace{.3in}
 
 \usetikzlibrary{positioning}
\usetikzlibrary{decorations.pathreplacing}

\begin{center}
\begin{tikzpicture}
   \node at (-1.2,3.8) {\normalsize{$g_1$}};
 \draw[->, dotted] (5,5.2) -- (-1,4);
   \node at (5,10.8) {\normalsize{$g_2$}};
 \draw[->, dotted] (5,5.2) -- (5,10.5);
 
    \coordinate (orig) at (5,5.2);
    
     \node at (orig) {$\textcolor{rgb:red,1;blue,1}{\bullet}$};
   \node at (5,4.6) {\normalsize{\Ovalbox{At $(\gamma_1,\gamma_2) = (0,0)$}}};

    \node at (1,9) {$\textcolor{blue}{\bullet}$};
     \node at (1.25,9.25) {$\textcolor{blue}{\bullet}$};
      \node at (0,8) {$\textcolor{red}{\bullet}$};
       \node at (0.25,8.25) {$\textcolor{red}{\bullet}$};
       
 \node at (-2,10.2) {\normalsize{\Ovalbox{At $(\gamma_1,\gamma_2) \neq (0,0)$}}};
 
 \draw[thin] (-1,8.7) rectangle (.7,7.3);
  \node at (-.1,6.9) {\normalsize{$\mathbb{A}^2_{(z_1,z_2,\gamma_1,\gamma_2)}$}};)
  \node at (-.4,8.3) {\normalsize{\red{$Y_1^{sing}$}}};)
 
  \draw[thin] (0.2,9.8) rectangle (1.9,8.5);
  \node at (1.1,10.2) {\normalsize{$\mathbb{A}^2_{(z_0,z_3,\gamma_1,\gamma_2)}$}};)
  \node at (.8,9.4) {\normalsize{\blue{$Y_2^{sing}$}}};)

    \node at (6,9) {$\textcolor{blue}{\bullet}$};
     \node at (5.75,9.25) {$\textcolor{blue}{\bullet}$};
      \node at (7,8) {$\textcolor{red}{\bullet}$};
       \node at (6.75,8.25) {$\textcolor{red}{\bullet}$};

\node at (9.3,10.2) {\normalsize{\Ovalbox{At $(\gamma_1',\gamma_2') \neq (\gamma_1,\gamma_2)$}}};

 \draw[thin] (6.3,8.7) rectangle (8.1,7.3);
  \node at (7.5,6.9) {\normalsize{$\mathbb{A}^2_{(z_1,z_2,\gamma_1',\gamma_2')}$}};)
  \node at (7.5,7.7) {\normalsize{\red{$Y_1^{sing}$}}};)
 
  \draw[thin] (5.3,9.8) rectangle (7,8.5);
  \node at (6.2,10.2) {\normalsize{$\mathbb{A}^2_{(z_0,z_3,\gamma_1',\gamma_2')}$}};)
  \node at (6.5,9.4) {\normalsize{\blue{$Y_2^{sing}$}}};)
  
    \draw[<-,blue,line width=1.2pt] (1,9) parabola (orig);
  \draw[<-, dashed,blue,line width=1.2pt](1.25,9.25)  parabola (orig);
    \draw[<-,blue,line width=1.2pt] (6,9) parabola (orig);
  \draw[<-, dashed,blue,line width=1.2pt]  (5.75,9.25) parabola (orig);
    \draw[<-,red,line width=1.2pt] (0,8) parabola (orig) ;
  \draw[<-, dashed,red,line width=1.2pt] (0.25,8.25) parabola  (orig);
      \draw[<-,red,line width=1.2pt] (7,8) parabola (orig) ;
  \draw[<-, dashed,red,line width=1.2pt] (6.75,8.25) parabola (orig) ;
\end{tikzpicture}
\end{center}

 \vspace{-2.33in}

\begin{center}
\textcolor{white}{.}
\hspace{-.25in} 
\begin{tikzpicture}[blend group=screen]

 \coordinate (orig) at (5,5.2);

    \begin{scope}
        \shade[left color=blue, right color={rgb:blue,1;white,3}, opacity=0.6]
           (orig) parabola bend (1.0,9.0)  (5,9.0);
    \end{scope}
             \begin{scope}
        \shade[left color=blue, right color={rgb:blue,1;white,3}, opacity=0.3]
            (orig) parabola bend  (1.25,9.25) (5,9.25);
    \end{scope}
     \begin{scope}
    \shade[left color=red, right color={rgb:red,1;white,3}, opacity=0.7]
            (orig) parabola bend  (0,8) (5,8)  ;
    \end{scope}
     \begin{scope}
        \shade[left color=red, right color={rgb:red,1;white,3}, opacity=0.4]
          (orig)  parabola bend (0.25,8.25)  (5,8.25);
    \end{scope}
    
       \begin{scope}
        \shade[right color=blue, left color={rgb:blue,1;white,3}, opacity=0.6]
           (orig) parabola bend (6,9.0)  (5,9.0);
    \end{scope}
             \begin{scope}
        \shade[right color=blue, left color={rgb:blue,1;white,3}, opacity=0.3]
            (orig) parabola bend  (5.75,9.25) (5,9.25);
    \end{scope}
     \begin{scope}
    \shade[right color=red, left color={rgb:red,1;white,3}, opacity=0.7]
            (orig) parabola bend  (7,8) (5,8)  ;
    \end{scope}
     \begin{scope}
        \shade[right color=red, left color={rgb:red,1;white,3}, opacity=0.4]
          (orig)  parabola bend (6.75,8.25)  (5,8.25);
    \end{scope}

\end{tikzpicture}
\end{center}

\vspace{.7in}

\begin{center}
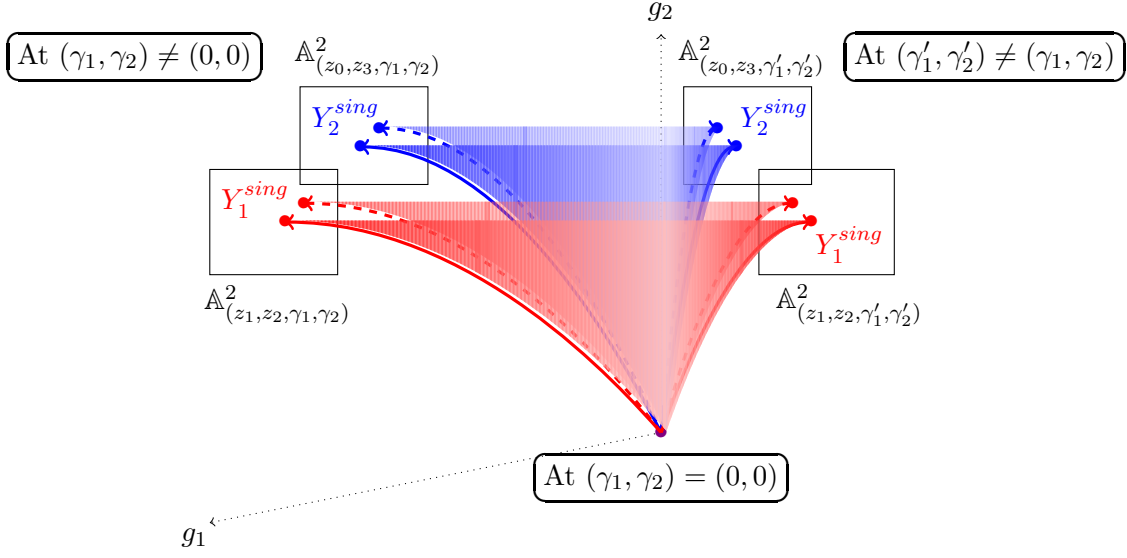
\textsc{Figure 2.} For $n$ even with $\rho = \rho_1$: $Y^{sing} =   Y^{symp}_0$, union of surfaces $Y_1^{sing}$ and $Y_2^{sing}$
\end{center}

\vspace{.5in}

\noindent{\it Proof of Theorem~\ref{tSEven}}.
We will only treat the case for $\rho=\rho_1$, and other cases are similar. We use the presentation of $Z$ in  Proposition~\ref{pZEven}, where we can write the two defining relations of $Y$ as
$$F_1=(a_1-h_1)^2-z_0z_3,\quad F_2=(a_2-h_2)^2-z_1z_2.$$
Moreover, we can assume that $h_1,h_2\in \kk[g_1,g_2]_n$ have no common factors and 
\begin{align*}
a_1&\, =\lambda (z_0+a^sb^s\xi^{3n}z_3)+\mu(z_1+a^sb^{-s}\xi^{n}z_2)\\ 
a_2&\, =\mu (b^{-s}c^{-n}\xi^{2n}z_0+a^sc^{-n}\xi^{3n}z_3)+\lambda (b^sz_1+a^s\xi^{3n}z_2),
\end{align*}
where $\lambda,\mu\in \kk$ with $\mu\neq 0$ by Proposition~\ref{pZEven}. 

Next, it is clear that $(Y^{sing})_{\gamma_1,\gamma_2}\subseteq (Y_{\gamma_1,\gamma_2})^{sing}\subset \mathbb{V}(z_0z_1z_2z_3)$ using Lemma \ref{lSSliceE} for the latter. We now show that $Y^{sing}$ is contained in the union of $Y_1^{sing}$ and $Y_2^{sing}$. Take a point $\widetilde{P} = (P, \gamma_1, \gamma_2) \in (Y^{sing})_{\gamma_1,\gamma_2}$, where $P= (p_0, \dots, p_3)$. Then without loss of generality, we can take $p_0=0$ and hence $a_1-h_1=0$ at $\widetilde{P}$ in $Y^{sing}$. 
Now, let $$A:=a_2-h_2$$ and consider the following Jacobian matrix (it has rank $\le 1$ at $\widetilde{P}$)

\[
\small
\frac{\partial(F_1,F_2)}{\partial(z_0, \dots, z_3,g_1,g_2)}\Big|_{a_1=h_1} = 
\begin{bmatrix}
-z_3 & 0 & 0 & -z_0 & 0 & 0\\
2A\frac{\partial a_2}{\partial z_0} & 2A\frac{\partial a_2}{\partial z_1}-z_2 
&2A\frac{\partial a_2}{\partial z_2}-z_1 & 2A\frac{\partial a_2}{\partial z_3} & -2A\frac{\partial h_2}{\partial g_1} & -2A\frac{\partial h_2}{\partial g_2}\end{bmatrix}
\]
\medskip

Recall that $p_0 = 0$. Now if $p_3=0$, then $\widetilde{P}$ is contained in $Y_1^{sing}$. On the other hand, suppose $p_3\neq 0$. Then the determinant of the first two columns of the matrix above, evaluated at $\widetilde{P}$ and set equal to 0, implies that $$p_2 = 2A\frac{\partial a_2}{\partial z_1}(\widetilde{P}).$$ Likewise, using the first and third columns, and also the first and fourth columns, respectively, we get
$$p_1 = 2A\frac{\partial a_2}{\partial z_2}(\widetilde{P}) \quad \quad \text{and} \quad \quad 
2A\frac{\partial a_2}{\partial z_3}(\widetilde{P}) = 0.$$
Since $\frac{\partial a_2}{\partial z_3} \neq 0$, we have that $A =0$ and thus $p_1 = p_2 =0$. Therefore, $\widetilde{P} \in Y_2^{sing}$.

Conversely, it is straight-forward to check that $Y_1^{sing}$ and $Y_2^{sing}$ are contained in $Y^{sing}$. So, $Y^{sing}=Y_1^{sing}\cup Y_2^{sing}$. 

Finally, let $\widetilde{P}=(p_0,p_1,p_2,p_3,\gamma_1,\gamma_2)\in Y_1^{sing}\cap Y_2^{sing}$. Then we have $p_k = 0$ for all $k$ (by Lemma~\ref{lSSliceE}) and $(\gamma_1,\gamma_2)\in \mathbb{V}(h_1,h_2)$. Since $h_1,h_2$ have no common factors, we get $\gamma_1=\gamma_2=0$ and $\widetilde{P}=0$, as desired.

\medskip

(1) follows from Corollary \ref{cneven-H4}. 

\medskip

(2) First, assume that $(\gamma_1,\gamma_2)=(0,0)$. Then $(Y^{sing})_{0,0}\subseteq (Y_{0,0})^{sing}$, where $Y_{0,0} = \kk[z_0,z_1,z_2,z_3]/(\Phi_1,\Phi_2)$ is the affine cone of the smooth elliptic curve $E''=E/\langle \sigma \rangle$ (see Lemma~\ref{lcenterB}). So we get $(Y^{sing})_{0,0}=(Y_{0,0})^{sing}=\{\underline{0}\}$. 

Now suppose $(\gamma_1,\gamma_2)\neq (0,0)$ and $\rho=\rho_1$. By the beginning of the statement, we have $Y_1^{sing}=\mathbb{V}(a_1-h_1,z_0,z_3,(a_2-h_2)^2-z_1z_2)$. So $(Y_1^{sing})_{\gamma_1,\gamma_2}$ is the intersection points of a line $a_1=h_1$ with a conic $(a_2-h_2)^2=z_1z_2$ in the affine space $\mathbb A^2_{(z_1,z_2,g_1=\gamma_1,g_2=\gamma_2)}$, which has two points (counting multiplicity). The same argument applies to $(Y_2^{sing})_{\gamma_1,\gamma_2}$ as well. Since $Y_1^{sing}~\cap~ Y_2^{sing}=\{\underline{0}\}$ by the beginning of the statement, we get that $(Y^{sing})_{\gamma_1,\gamma_2}=(Y_1^{sing})_{\gamma_1,\gamma_2}
\cup (Y_2^{sing})_{\gamma_1,\gamma_2}$ has 4 points generically.

The argument for $\rho = \rho_2, \rho_3$ follows similarly.

\medskip

(3) Recall that $(Y_{\gamma_1,\gamma_2})^{sing}\subset \mathbb{V}(z_0z_1z_2z_3)$  by Lemma \ref{lSSliceE}. Then, without loss of generality take $\rho = \rho_1$, and take $z_0 =0$ so that $a_1=h_1(\underline{\gamma})$, and consider the Jacobian matrix of $Y_{\gamma_1,\gamma_2}$
\[
\small
\frac{\partial(F_1(\underline{z}, \underline{\gamma}),F_2(\underline{z}, \underline{\gamma}))}{\partial(z_0, z_1,z_2, z_3)}\Big|_{z_0=0,a_1=h_1}= 
\begin{bmatrix}
-z_3 & 0 & 0 & 0 \\
2A\frac{\partial a_2}{\partial z_0} & 2A\frac{\partial a_2}{\partial z_1}-z_2 
&2A\frac{\partial a_2}{\partial z_2}-z_1 & 2A\frac{\partial a_2}{\partial z_3} \end{bmatrix}
\]
where $A=a_2-h_2(\overline{\gamma})$. With this matrix, we can show $(Y_{\gamma_1,\gamma_2})^{sing}\subseteq (Y_1)^{sing}_{\gamma_1,\gamma_2}\,\cup\, (Y_2)^{sing}_{\gamma_1,\gamma_2}$. Moreover, we can conclude that
$$(Y_1)^{sing}_{\gamma_1,\gamma_2}\,\cup\, (Y_2)^{sing}_{\gamma_1,\gamma_2} \quad = \quad (Y^{sing})_{\gamma_1,\gamma_2}\quad \subseteq\quad  (Y_{\gamma_1,\gamma_2})^{sing}\quad \subseteq\quad (Y_1)^{sing}_{\gamma_1,\gamma_2}\,\cup\, (Y_2)^{sing}_{\gamma_1,\gamma_2}.$$
This proves the result.
\qed


\sectionnew{Background material on Poisson orders and specialization} 
\label{sec:bkgdPorder}

We discuss briefly in this section background material on Poisson orders, including the process of specialization mentioned in the introduction, as well as material on symplectic cores. More details can be found in \cite[Sections~2.1 and~2.2]{WWY} and the references therein.


\subsection{Poisson orders and specialization} \label{sec:P-order}

Here we collect some definitions and facts about Poisson orders and describe an extension of the specialization technique 
for obtaining such structures. The following definition is due to Brown and Gordon \cite{BrownGordon}.

\bde{Pois-order} [$\Der(A/C)$, $\partial$, $\partial_z$] Let $A$ be a $\kk$-algebra which is module-finite over a central subalgebra $C$. Denote by $\Der(A/C)$\index{Der@$\Der(A/C)$}
the algebra of $\kk$-derivations of $A$ that preserve $C$.

The algebra $A$ is called a {\it Poisson $C$-order} if there exists a $\kk$-linear map
$$\del : C \to \Der(A/C)\index{p2artial@$\del, \del_z$}$$ such that the induced bracket $\{.,.\}$ on $C$, given~by 
\begin{equation}
\label{Poisson}
\{z, z' \}:= \partial_{z}(z'), \quad z, z' \in C,
\end{equation}
makes $C$ a Poisson algebra. 
The triple $(A,C, \partial \colon C \to \Der(A/C))$ will be also called a {\it Poisson order} in places where the role of $\partial$ 
needs to be emphasized.
\smallskip
\ede

As discussed in \cite[Section~2.2]{BrownGordon}, specializations of families of algebras give rise to Poisson orders. 
In our previous work \cite{WWY} we generalized this construction to obtain Poisson orders from higher degree terms in the derivation $\partial$; this is reviewed as follows.

\bde{spec'n} 
Let $R$ be an algebra over $\kk$ and $\hbar$ be a central element of $R$ which is regular, i.e., not a zero-divisor of $R$. We refer to the $\kk$-algebra $R_0 := R/\hbar R$ as the {\it specialization of $R$ at $\hbar \in Z(R)$}.
\ede

\bnota{theta}[$\theta$, $\iota$, $N$] Retain the notation of Definition~\ref{dspec'n}. Let $[\text{-},\text{-}]$ denote the commutator of elements of $R$.
Let $\theta : R \twoheadrightarrow R_0$\index{t2heta@$\theta$} be the canonical projection; so, $\ker \theta~=~\hbar R$.
Fix a linear map $\iota : Z(R_0) \hra R$\index{i@$\iota$} such that $\theta \circ \iota = \id_{Z(R_0)}$. Let $N \in \Zset_+$ be such that 
\begin{equation}
\label{N}
[\iota(z), y] \in \hbar^N R \quad \mbox{for all} \quad z \in Z(R_0), ~y \in R.
\end{equation}
\enota

Note that \eqref{N} holds for $N=1$: take $\wt{y} \in \theta^{-1}(y)$ for $y \in R_0$ and we get  $\theta([\iota(z), \wt{y}]) = [z, y] =0$.

\bde{spec-deriv} Retain the notation above. For $y \in R_0$ and $z \in Z(R_0)$, the {\it special derivation of level $N$}\index{N@$N$} is defined (in fact, well-defined) as 
\begin{equation}
\label{def-special}
\partial_z (y) := \theta \left( \frac{ [\iota(z), \wt{y}]}{\hbar^N}\right), \quad \mbox{where} \; \; \wt{y} \in \theta^{-1}(y). 
\end{equation}
\ede

The next result states that $\partial_z$ is indeed a derivation, and thus specializations yield Poisson orders.

\bpr{special1} \cite[Proposition~2.7 and Corollary~2.8]{WWY} Let $R$ be a $\kk$-algebra and $\hbar \in Z(R)$ be a regular element. Assume that $\iota : R_0 := R/(\hbar R) \hra R$ 
is a linear section of the specialization map $\theta : R \twoheadrightarrow R_0$ such that \eqref{N} holds for some $N \in \Zset_+$. 
Assume that $R_0$ is module-finite over $Z(R_0)$. 
\begin{enumerate}
\item If, for all $z \in Z(R_0)$, $\partial_z$ is a special derivation of level $N$, then  $$(R_0, Z(R_0), \partial \colon Z(R_0) \to \Der(R_0/ Z(R_0)))$$ is a Poisson order and the map $\partial$ is a homomorphism of Lie algebras.
\smallskip

\item If $C \subset Z(R_0)$ is a Poisson subalgebra of $Z(R_0)$ with respect to the Poisson structure \eqref{Poisson}
and $R_0$ is module-finite over $C$, then $R_0$ is a Poisson $C$-order via the restriction of $\partial$ to $C$. 
\smallskip

\item If, in addition to (2), the restricted section $\iota : C \hra R$ is an algebra homomorphism, then
\[
\partial_{z z'}(y) = z \partial_{z'}(y) + z' \partial_{z}(y) \quad
\mbox{for} \quad z, z' \in C, y \in R_0.
\]
\end{enumerate}
\epr

We coined such a construction with the following terminology.

\bde{N-Pois-ord} \cite[Definition~2.9]{WWY} The Poisson order produced in Proposition~\ref{pspecial1} is  a {\it Poisson order of level $N$} when the level of the special derivation needs to be emphasized.
\ede


\subsection{Symplectic cores and the Brown-Gordon theorem}
\label{sec:BGor-thm}
Poisson orders can be used to establish isomorphisms for different central quotients 
of a PI algebra via the result of Brown and Gordon \cite{BrownGordon} provided below. The result relies on the notion of {\em{symplectic core}}, introduced in \cite{BrownGordon}. We recall some terminology from \cite[Section 3.2]{BrownGordon}.

\bde{symplcore}[$\mathcal{P}(I)$]  
Let $(C, \{.,.\})$ be an affine Poisson algebra over a field $\kk$ of characteristic $0$. For every ideal $I$ of $C$, there exists a unique maximal Poisson ideal contained in $I$, to be denoted by $\PP(I)$\index{P1mathcal@$\PP(I)$}. If $I$ is prime, then  $\PP(I)$ is Poisson prime, \cite[Lemma~6.2]{Goodearl}. 

\begin{enumerate}
\item We refer to $\mathcal{P}(I)$ above as the {\it Poisson core} of $I$.
\smallskip
\item We say that two maximal ideals $\mm, \nn \in \maxSpec C$ of an affine Poisson algebra $(C, \{.,.\})$ are {\it equivalent} if $\PP(\mm) = \PP(\nn)$. 
\smallskip
\item The equivalence class of $\mm \in \maxSpec C$ is referred to as the  {\em{symplectic core}} of $\mm$. 
The corresponding partition of $\maxSpec C$ is called {\em{symplectic core partition}}. 
\end{enumerate}
\ede

\smallbreak One main benefit of using the symplectic core partition is the powerful result below.

\bth{BrGor} \cite[Theorem~4.2]{BrownGordon} Assume that $\kk = \mathbb{C}$ and that $A$ is a Poisson $C$-order which is an affine $\Cset$-algebra. If $\mm, \nn \in \maxSpec C$ are in the same symplectic core, 
then there is an isomorphism between the corresponding finite-dimensional $\mathbb{C}$-algebras
\[
A/(\mm A) \cong A /(\nn A).
\]

 \vspace{-.25in}
 
  \qed
\eth


\sectionnew{A specialization setting for 4-dimensional Sklyanin algebras} 
\label{sec:formal}
The goal of this section is to produce a setting so that the PI 4-dimensional Sklyanin algebras arise as Poisson orders via specialization; see Section~\ref{sec:P-order}. The section also sets up  
some of the notation regarding Poisson orders that we will use throughout this work.

Recall that $S:= S(\alpha,\beta,\gamma)$ is a 4-dimensional Sklyanin algebra  and {\it we do not necessarily need that  $S$  is module-finite over its center $Z$}. In any case, recall that $B ~(\cong S/(Sg_1+Sg_2))$ is the corresponding twisted homogeneous coordinate ring.

\smallbreak
The reader may wish to view Figure~3 at this point for a preview of the setting that we will construct for $S$. Our objective is to produce a degree 0 deformation ${S_\hbar}$ of $S$ using a formal parameter $\hbar$. The specialization map for $S$ will 
be realized via a canonical projection $\theta_S : {S_\hbar} \to S$ given by $\hbar \mapsto 0$. Moreover, ${S_\hbar}$ will have the structure of a $\kk[[\hbar]]$-algebra.

\bnota{hbar} [$\hbar$, $\wt{\alpha}$, $\wt{\beta}$, $\wt{\gamma}$]
\index{h0bar@$\hbar$}
\index{alphatilde@$\wt{\alpha}$}
\index{betatilde@$\wt{\beta}$}
\index{g3ammatilde@$\wt{\gamma}$}
To begin, we fix a formal parameter $\hbar$ and let
$$\wt{\alpha}:= \alpha + \alpha_1 \hb+\alpha_2\hb^2+\cdots, \;
\wt{\beta}:= \beta + \beta_1 \hb+\beta_2\hb^2+\cdots, \;
\wt{\gamma}:= \gamma + \gamma_1 \hb+\gamma_2\hb^2+\cdots,$$
in $\kk[[\hbar]]$ satisfying $\wt{\alpha}+\wt{\beta}+\wt{\gamma}+\wt{\alpha}\wt{\beta}\wt{\gamma}=0$ (a version of $\eqref{albega-cond1}$).
\enota

By our choice \eqref{albega-cond2}, we know $(\alpha, \beta, \gamma) \not \in \{(-1,1,\gamma),(\alpha,-1,1),(1,\beta,-1)\}$. Hence, it is clear that 
\begin{align*}
(\widetilde{\alpha}, \widetilde{\beta}, \widetilde{\gamma}) \not \in \{(-1,1,\widetilde{\gamma}),(\widetilde{\alpha},-1,1),(1,\widetilde{\beta},-1)\}.
\end{align*}

\bde{ext-formalS}[$\wh{S}_\hbar$, ${S}_\hbar$] Denote by $\wh{S}_\hbar$\index{S3hbar@$S_\hbar, \wh{S}_\hbar$} the 4-dimensional Sklyanin algebra over $\ol{\kk((\hb))}$
with parameters $(\wt{\alpha}, \wt{\beta}, \wt{\gamma})$. Define the {\em{formal Sklyanin algebra}} to be the 
$\kk[[\hb]]$-subalgebra ${S}_\hb$\index{S3hbar@$S_\hbar, \wh{S}_\hbar$} of $\wh{S}_\hb$ generated by $x_0,x_1,x_2,x_3$, that is,
\[
{S}_\hb := \kk[[\hb]]\lcor  x_0, x_1, x_2,x_3 \rcor \subset \wh{S}_\hb.
\] 
\ede

It is important to point out that ${S}_\hb$ is a graded $\kk[[\hb]]$-algebra with the grading inherited from $\wh{S}_\hb$, such that $\text{deg}(\hb)=0$ and $\text{deg}(x_i) =1$ for $0\le i\le 3$. Notice that by our choice, $1\pm\widetilde{\alpha}, 1\pm\widetilde{\beta},  1\pm \widetilde{\gamma}\not\in (\hb)$. So, these elements are all invertible in the formal power series $\kk[[\hb]]$.  Recall by Lemma~\ref{lfactor-g}, we obtain the following result.

\ble{centtilde} \textnormal{[$\wt{g_1}$, $\wt{g_2}$]}
\index{g120tilde@$\wt{g_1}, \wt{g_2}$}
The elements 
$$
\wt{g_1}= -x_0^2+x_1^2+x_2^2+x_3^2  \quad \text{ and } \quad 
  \wt{g_2}=x_1^2 + (1+\wt{\alpha})(1-\wt{\beta})^{-1}x_2^2 +  (1-\wt{\alpha})(1+\wt{\gamma})^{-1}x_3^2 $$
form a central regular sequence in  $\wh{S}_\hb$.  \qed
\ele

\ble{Specialization} 
The following statements hold for the formal Sklyanin algebra ${S}_\hbar$. 
\begin{enumerate}
\item $\wh{S}_\hb\cong \ol{\kk((\hb))}\otimes_{\kk[[\hb]]} {S}_\hb$.
\item At each degree $d$ of ${S}_\hb$, we get that $({S}_\hb)_d$ is a free $\kk[\hb]]$-module of rank ${ d+3 \choose 3}$.
\item The elements $\wt{g_1},\wt{g_2}$ belong to the center of ${S}_\hb$. 
\item There is a natural surjection from ${S}_\hb\twoheadrightarrow S$ via $\hb\mapsto 0$ with kernel equal to $\hb {S}_\hb$.
\end{enumerate}
\ele
\begin{proof}
(1) This is clear from the definitions of $\wh{S}_\hbar$ and $S_\hbar$. 

\smallskip

(2) Since $\wh{S}_\hb$ is a domain, each graded piece $({S}_\hb)_d\subset (\wh{S}_\hb)_d$ is a finitely generated torsion-free module over $\kk[[\hb]]$. 
Because $\kk[[\hb]]$ is a PID, this implies that $({S}_\hb)_d$ is a free $\kk[[\hbar]]$-module. By (1), $(S_\hbar)_d$ has rank equal to $\dim (\wh{S}_\hb)_d$, which in turn is equal to ${ d+3 \choose 3}$ for $\wh{S}_\hb$ has Hilbert series $1/(1-t)^4$. 

\smallskip

(3) It is easy to check that $\wt{g_1},\wt{g_2}\in Z(\wh{S}_\hb)\cap {S}_\hbar\subset Z({S}_\hb)$.

\smallskip

(4) It suffices to show that ${S}_\hb/\hb{S}_\hb\cong S$. Clearly there is a surjection $S \twoheadrightarrow {S}_\hb/\hb{S}_\hb$. Moreover, it is an isomorphism since on each degree $\dim S_d=\dim ({S}_\hb)_d/\hb({S}_\hb)_d={ d+3 \choose 3}$ by (2). The kernel part follows directly. 
\end{proof}

\bnota{thetaS}[$\theta_S$] Denote by $\theta_S$ the corresponding {\it specialization map for the formal Sklyanin algebra ${S}_\hbar$}, namely
$$
\theta_S : {S}_\hb \to S \quad \mbox{given by} \quad \hb \mt 0.
\index{t3hetaS@$\theta_S$} 
$$
\enota

So, the first column in Figure 3 below is established and we now turn our attention to the second column of that figure.

\bde{formalthcr}[$E_\hb$, $\mathcal{L}_\hb$, $\sigma_\hb$, $\wh{B}_\hb$, ${B}_\hb$]
Denote by $E_\hb$\index{E3hbar@$E_\hb$} 
the elliptic curve over $\ol{\kk((\hb))}$, by $\mathcal{L}_\hb=\mathcal O_{\mathbb P^3}(1)|_{E_\hb}$\index{L1hbar@$\mathcal{L}_\hb$} 
the invertible sheaf over $E_\hb$, and by $\sigma_\hb$\index{s6igmahbar@$\sigma_\hb$} 
the automorphism of $E_\hb$
corresponding to $\wh{S}_\hb$ as in \dlref{geomS} with $(\alpha,\beta,\gamma)$ replaced by $(\wt{\alpha}, \wt{\beta}, \wt{\gamma}$) from Notation~\ref{nhbar}. Let
$$
\wh{B}_\hb:=B(E_\hb, \mathcal{L}_\hb, \sigma_\hb)      
\index{B0hbar@$B_\hb, \wh{B}_\hb$}
$$
be the corresponding twisted homogeneous coordinate ring. Its $\kk[[\hb]]$-subalgebra 
$$
{B}_\hb := \kk[[\hb]] \lcor x_0,x_1,x_2,x_3 \rcor \subset \wh{B}_\hb,
\index{B0hbar@$B_\hb, \wh{B}_\hb$}
$$
generated by 
$x_0,x_1,x_2,x_3$, will be called {\em{formal twisted homogeneous coordinate ring}}.
\ede

\ble{thetaB} \textnormal{[$\psi_\hbar$]} The canonical projection $\wh{S}_\hb\twoheadrightarrow \wh{B}_\hb$ 
induces a surjection $\psi_\hb :{S}_\hb\twoheadrightarrow {B}_\hb$, \index{p6sihbar@$\psi_\hbar$}
whose kernel is generated by $\wt{g_1},\wt{g_2}$. As a consequence, the composition ${S}_\hb \stackrel{\theta_S}{\twoheadrightarrow} S\twoheadrightarrow B$ factors through the map $\psi_\hb$. 
\ele

\begin{proof}
Note that $\psi_\hb$ is the composition ${S}_\hb\hookrightarrow \wh{S}_\hb\twoheadrightarrow \wh{B}_\hb$, whose image is ${B}_\hb$ by definition. So it remains to show 
$$\mathrm{ker}(\psi_\hb)~=~(\wt{g_1}\wh{S}_\hb+\wt{g_2}\wh{S}_\hb)\cap S_\hb~=~\wt{g_1}S_\hb+\wt{g_2}S_\hb.$$
For simplicity, we write $M:=\wt{g_1}{S}_\hb$ and $N:=\wt{g_1}\wh{S}_\hb\cap {S}_\hb$. It is clear that $M\subseteq N$, which are both graded ideals in ${S}_\hb$. By Lemma \ref{lSpecialization}(1), we have $M \otimes_{\kk[[\hb]]}\ol{\kk((\hb))}=N$. Hence, as finitely generated free modules over $\kk[[\hb]]$, we get that $\mathrm{rank}(M_d)=\mathrm{rank}(N_d)$  in each degree $d$. So $(N/M)_d$ is a torsion module over $\kk[[\hb]]$. Then for any $r\in N_d$, there exists an integer $m$ such that $\hb^mr=y\wt{g_1}$ for some $y\in ({S}_\hb)_{d-2}$. Since $\hb\nmid \wt{g_1}$, we get that $\hb^d\mid y$ and $r\in M_d$. This implies $M_d=N_d$ for every $d$, and $M=N$. By a similar argument, we can then conclude that $(\wt{g_1}\wh{S}_\hbar+\wt{g_2}\wh{S}_\hbar) \cap S_\hbar = \wt{g_1}S_\hbar + \wt{g_2}S_\hbar$ where the left-hand side is equal to ker($\psi_\hbar$). Finally, the factorization through $\psi_\hb$ is straight-forward. 
\end{proof}

\bde{thetaB} \textnormal{[$\theta_B$]}
Let $\theta_B: B_\hbar \twoheadrightarrow B$\index{t3hetaB@$\theta_B$} 
be the map induced by Lemma~\ref{lthetaB}, which we call the {\it specialization map for the formal twisted homogeneous coordinate ring $B_\hbar$}.
\ede

Now we complete the verification of Figure~3 as follows.  

\label{sec:4Col}
\bde{Ahbar} [$w$, $L_\hbar$, $v_0'$, $v_1'$, $v_2'$, $v_3'$, $R_\hbar$]
Fix a nonzero section $w\in H^0(E,\mathcal L)$,\index{w1@$w$} 
which is also realized as an element of $H^0(E_\hb,\mathcal L_\hb)$ via the vector space isomorphism $H^0(E_\hb,\mathcal L_\hb)\cong \ol{\kk((\hb))}\otimes_\kk H^0(E,\mathcal L)$. 
Denote by $L_\hbar$\index{L0hbar@$L_\hbar$} 
the $\kk[[\hbar]]$-subalgebra of the function field $\ol{\kk((\hbar))}(E_\hbar)$ generated by 
$v_0':= v_0/w$, $v_1':= v_1/w$, $v_2':=v_2 /w$ and $v_3'=v_3/w$ 
\index{vsprime@$v'_0, \ldots, v'_3$}
satisfying the dehomogenized relations: 
\[
\begin{array}{ll}
\smallskip

\phi_1(v_0', v_1', v_2', v_3') &=v_0'^2+v_1'^2+v_2'^2+v_3'^2=0\\
\phi_2(v_0',v_1',v_2',v_3') &= (1-\wt{\gamma})(1+\wt{\alpha})^{-1}v_1'^2 + (1+\wt{\gamma})(1-\wt{\beta})^{-1}v_2'^2 + v_3'^2=0.
\end{array}
\]
We call $R_\hbar := (L_\hbar)_{(\hbar)}$\index{R0hbar@$R_\hbar$}
the {\it integral form} of the field $\ol{\kk((\hbar))}(E_\hbar)$.
\ede

Note that the quotient field of $L_\hbar$ is $Q(L_\hbar) \cong \ol{\kk((\hbar))}(E_\hbar)$. Using \eqref{sigma} with replacing $(\alpha,\beta,\gamma)$ by $(\wt{\alpha}, \wt{\beta}, \wt{\gamma})$, one sees that the automorphism $\sigma_\hbar \in \Aut \left( \ol{\kk((\hbar))}(E_\hbar) \right)$ restricts to an automorphism of $R_\hbar$, 
given by 
\[
\sigma_\hb(v_i')=\sigma_\hb(v_i/w)=\frac{\sigma_\hb(v_i)/w^2}{\sigma_\hb(w)/w^2},\ \text{for all}\ 0\le i\le 3.\]
Let $D$ and $D_\hb $ be the divisors of zeros of $w$ for $E$ and $E_\hb$, respectively. In view of Lemma~\ref{lembedB}, we have the following commutative diagram.
\[
\xymatrix{
(\wh{B}_\hb)_1=H^0(E_\hb,\mathcal L_\hb)\ar[r]^-{\cong} & H^0(E_\hb,\mathcal O(D_\hb))t\ar@{^{(}->}[r]&\ol{\kk((\hb))}(E_\hb)[t^{\pm 1};\sigma_\hb]\\
({B}_\hb)_1=\kk[[\hb]]\otimes_\kk H^0(E,\mathcal L)\ar[r]^-{\cong}\ar@{^{(}->}[u]\ar@{->>}[d]_-{\hb\;\mapsto0}&\kk[[\hb]]\otimes_\kk H^0(E,\mathcal O(D_\hb))t\ar@{^{(}->}[r]\ar@{^{(}->}[u]\ar@{->>}[d]_-{\hb\;\mapsto 0}&R_\hb[t^{\pm 1}; \sigma_\hb]\ar@{^{(}->}[u]\ar@{->>}[d]_-{\hb\;\mapsto 0}\\
B_1=H^0(E,\mathcal L)\ar[r]^-{\cong} & H^0(E,\mathcal O(D))t\ar@{^{(}->}[r] &\kk(E)[t^{\pm 1};\sigma]
}
\]
Moreover, we have the canonical embeddings
\begin{equation}
\label{BtoR}
{B}_\hbar \hra R_\hbar [t^{\pm 1}; \sigma_\hbar] \hra \ol{\kk((\hb))}(E_\hbar)[ t^{\pm 1}; \sigma_\hbar].
\end{equation}
The ring $R_\hbar [t^{\pm 1}; \sigma_\hbar]$ is a graded localization of ${B}_\hbar$ by an Ore set 
which does not intersect the kernel $\ker \theta_B = \hbar {B}_\hbar$. Therefore, the following map is well-defined.

\bde{thetaR}[$\theta_R$] Let $\theta_R : R_\hbar[t^{\pm 1}; \sigma_\hbar] \twoheadrightarrow \kk(E)[t^{\pm 1}, \sigma]$\index{t3hetaR@$\theta_R$}  
be defined by $$\theta_R(v_i') =v_i',\ \text{for all}\ 0\le i\le 3, \quad \theta_R(t)=t, \quad \theta_R(\hbar)=0,$$
which is the extension of $\theta_B$ via localization. We also denote by $\theta_R$ its restriction to the specialization map $R_\hbar \twoheadrightarrow \kk(E)$. These maps are referred to as the {\it specialization maps for the integral form of the formal twisted homogeneous coordinate $B$.}
\ede

The commutativity of the cells in Figure~3 between the second and third column 
follows directly from the definitions of the maps within this part of the diagram and from the previous commutative diagram. 

\vspace{.2in}

\[\small
\xymatrix@C=2.1em@R=3em{
 \widehat{S}_\hb=S(\wt{\alpha},\wt{\beta},\wt{\gamma})\ar@{->>}[rr]^-{\text{mod}(\wt{g_1},\wt{g_2})}&& \widehat{B}_\hb=B(E_\hb,\mathcal L_\hb,\sigma_\hb) \ar@{^{(}->}[rr]^-{\text{gr. quot. ring}} &&  \overline{\kk((\hb))}(E_\hb)[\,t^{\pm 1};\sigma_\hb]\\
 {S}_\hb=k[[\hb]]\langle x_0,x_1,x_2,x_3\rangle  \ar@{->>}[rr]^-{\text{mod}(\wt{g_1},\wt{g_2})}\ar@{^{(}->}[u]\ar@{->>}[d]^-{\theta_S} && {B}_\hb=k[[\hb]]\langle x_0,x_1,x_2,x_3\rangle \ar@{^{(}->}[rr]\ar@{^{(}->}[u]\ar@{->>}[d]^-{\theta_B} && R_\hb[\,t^{\pm 1};\sigma_\hb]\ar@{^{(}->}[u]\ar@{->>}[d]_-{\theta_R}\\
S=S(\alpha,\beta,\gamma) \ar@{->>}[rr]^-{\text{mod}(g_1,g_2)} && B=B(E,\mathcal L,\sigma)\ar@{^{(}->}[rr]^-{\text{gr. quot. ring}}  && \kk(E)[\,t^{\pm 1};\sigma]   \\
Z\left(S\right) \ar@{->>}[rr]^-{\text{mod}(g_1,g_2)}\ar@{^{(}->}[u]&&Z(B)\ar@{^{(}->}[u] \ar@{^{(}->}[rr]^-{\text{gr. quot. ring}} && \kk(E)^{\sigma}[\,t^{\pm n};\sigma]\ar@{^{(}->}[u] 
}
\]
\vspace{.1in}
\begin{center}
\textsc{Figure 3.} Specialization setting for Sklyanin algebras:\\
integral forms, Poisson orders, and centers are respectively in the last three rows.
\end{center}

\vspace{.2in}

\sectionnew{Poisson orders on PI 4-dimensional Sklyanin algebras $S$}
\label{sec:constr-Pord}
This section establishes a construction of Poisson orders on PI 4-dimensional Sklyanin algebras
with the property that the induced Poisson structure on its center is nontrivial. 

Recall that $S$ is a PI 4-dimensional Sklyanin algebra of PI degree $n = |\sigma|$ corresponding 
to the parameters $\al, \be, \ga \in \kk$ satisfying \eqref{albega-cond1} and \eqref{albega-cond2}.
Let the scalars $\al_i, \be_i, \ga_i$ defining $\wt{\al}, \wt{\be}, \wt{\ga}$ in Notation~\ref{nhbar}, for $i \geq 1$, be such that 
\begin{equation}
\label{cond-hdef}
\wt{\alpha}+\wt{\beta}+\wt{\gamma}+\wt{\alpha}\wt{\beta}\wt{\gamma}=0
\quad \quad 
\mbox{and} 
\quad \quad 
| \sig_\hb| = \infty.
\end{equation}
Such scalars exist because 
$
\{ (\al, \be, \ga) \in {\mathbb{A}}^3 \mid |\sig_{\al \be \ga}|=n \} 
$
is a closed subset of the surface 
$
{\mathbb{V}}(\al + \be + \ga + \al \be \ga) \subset {\mathbb{A}}^3.
$


\subsection{Construction of Poisson orders}
\label{Constr-Theorem} 

Throughout the section we will identify the first graded pieces $S_1$ of $S$ and $B_1$ of $B$ with each other through vector space isomorphism (Lemma~\ref{lfactor-g}). In particular, 
 $\{x_0, x_1, x_2, x_3\}$ is a good basis of $B_1$ as in \leref{actionrho}. 
 
 \bnota{xtilde} [$\wt{x}_i$, $\wt{z}_i$]
Denote by $\wt{x}_i$\index{x3stilde@$\wt{x}_0, \ldots, \wt{x}_3$} the preimage of $x_i$ under the specialization map
$\theta_S : {S}_\hbar \twoheadrightarrow S$ which is given by the same linear combinations of the standard generators of ${S}_\hbar$ 
as is $x_i$ given in terms of the standard generators of $S$. Moreover, set
$$
\wt{z}_i := \wt{x}_i^n + \textstyle \sum_{1 \leq j < n/2} c_{ij} \wt{g}^j \wt{x}_i^{n-2j} \in {S}_\hbar, \quad 0 \leq i \leq 3
\index{z3stilde@$\wt{z}_0, \ldots, \wt{z}_3$}
$$
for the polynomials $c_{ij} \in \kk[g_1,g_2]_{2j}$ from \prref{centerS}(1). 
\enota

\bde{good-ect} 
Consider the terminology below.
\begin{enumerate} 
\item A degree 0 section $\iota : Z \hra {S}_\hbar$ of the map $\theta_S : {S}_\hbar \twoheadrightarrow S$ will be called {\it good} if 
\begin{enumerate}
\item[(i)] $\iota(z_i) - 
\wt{z}_i \in \wt{g}_1 \kk [ \wt{x}_i, \wt{g}_1, \wt{g}_2, \hbar ]
+  \wt{g}_2 \kk [\wt{x}_i, \wt{g}_1, \wt{g}_2, \hbar ]$,
and
\item[(ii)] $\iota(g_1) = \wt{g}_1$ and $\iota(g_2) = \wt{g}_2$.  
\end{enumerate}

\smallskip

\item A specialization map $\theta_S:~ S_\hb~\to~S$ will be called a {\it good specialization of $S$ of level $N$} if there exists a good section 
$\iota : Z \hra S_\hbar$  such that 
\begin{equation}
\label{NN}
[\iota(z), y] \in \hbar^N {S}_\hbar \quad \mbox{for all} ~ z \in Z,~ y \in {S}_\hbar. 
\end{equation}
\end{enumerate}
\ede

Given any section $\iota : Z \hra {S}_\hbar$ of $\theta_S$,
\[
[\iota(z), y] \in \hbar {S}_\hbar \quad \mbox{for all}~ z \in Z, ~y \in {S}_\hbar. 
\]
Therefore, $N \geq 1$. Next, we prove that for a given PI 4-dimensional Sklyanin algebra~$S$, 
the levels $N$ of good specializations for $S$ are bounded from above.

\ble{exist} Let $S$ be a PI Sklyanin algebra corresponding 
to the parameters $\al, \be, \ga \in \kk$ satisfying \eqref{albega-cond1} and \eqref{albega-cond2}, and choose any formal parameters $\wt{\alpha}$, $\wt{\beta}$, $\wt{\gamma}$ satisfying \eqref{cond-hdef}. 
Then, the set of levels $N$ for good sections of $\theta_S$ has an upper bound. 
\ele

\begin{proof} Consider the algebras $L_\hbar$ and $R_\hbar :=(L_\hbar)_{(\hbar)}$ from \deref{Ahbar} for $w= \wt{x}_0$, identified as the preimage of $v_0'$ under $\theta_B$.
Since the condition \eqref{cond-hdef} is satisfied, the automorphism $\sigma_\hbar$ of $R_\hbar$ has infinite order.
This gives
\begin{equation}
\label{Nbound1}
R_\hbar^{\sigma_\hbar^n} \subsetneq R_\hbar.
\end{equation}
On the other hand, $\hb \in R_\hb$ is a regular element, from which one gets 
\begin{equation}
\label{Nbound2}
\textstyle \bigcap_{l \in \Zset_+} \hbar^l R_\hbar =0.
\end{equation}
Combining \eqref{Nbound1} and \eqref{Nbound2}, we obtain that there exists a positive integer $M$ such that 
\begin{equation}
\label{MM}
r^{\sigma_\hbar^n} -r \notin \hbar^M R_\hbar \quad \mbox{for some} \quad r \in R_\hbar.
\end{equation}
Denote by $M_{\min}$ the least such positive integer, and by $r_{\min} \in R_\hbar$ an element satisfying~\eqref{MM}
for this integer $M_{\min}$. 
Although, it will not play a role in the proof, we note that $M_{\min} \geq 2$ because 
$r^{\sigma_\hbar^n} - r \in \hbar R_\hbar$ for all $r \in R_\hbar$.

In the remainder of the proof, we will show that
\[
N<M_{\min}
\]
which gives the stated upper bound for levels of good sections. Let $\iota : Z \hra {S}_\hbar$ be one such good section of $\theta_S$, 
with respect to which the specialization map $\theta_S$ has level $N$. Recall the map $\psi_\hb :{S}_\hb\twoheadrightarrow {B}_\hb$
from \leref{thetaB}.
Then, the condition \eqref{NN} implies that
\[
[\psi_\hb \iota(z), y] \in \hbar^N {B}_\hbar, \quad \quad \forall \; z \in Z, ~y \in {B}_\hbar. 
\]
From \eqref{BtoR}, we have ${B}_\hb \hookrightarrow R_\hb[t^{\pm 1}; \sigma_\hb]$ with respect to which $R_\hbar[t^{\pm 1}; \sigma_\hbar]$ 
is a localization of ${B}_\hbar$. Hence,  
\begin{equation}
\label{zw-comm}
[\psi_\hb \iota(z), y] \in \hbar^N R_\hbar[t^{\pm 1}; \sigma_\hbar], \quad \quad \forall \; 
z \in Z, ~y \in R_\hbar[t^{\pm 1}; \sigma_\hbar]. 
\end{equation}

\prref{centerS}, combined with the facts that $\ker \psi_\hb = \wt{g}_1 {S}_\hbar+ \wt{g}_2 {S}_\hbar$ (from Lemma~\ref{lthetaB}) and that $\iota$ is a good section, leads to 
$\psi_\hb \iota(z_1) = t^n \in R_\hbar[t^{\pm 1}; \sigma_\hbar]$. By applying~\eqref{zw-comm} to $z = z_1$ and $y=r_{\min} \in R_\hbar$, 
we obtain that
\[
[t^n,r_{\min}]=(r_{\min}^{\sigma_\hbar^n}-r_{\min})t^n\in \hbar^N R_\hbar[t^{\pm 1}; \sigma_h] \quad \mbox{for all} \quad r \in R_\hbar.
\]
Since $t$ is a unit of $R_\hbar[t^{\pm 1}; \sigma_h]$, 
\[
r_{\min}^{\sigma_\hbar^n} -r_{\min} \in \hbar^N R_\hbar.
\]
The last equation implies that $N < M_{\min}$ because $r_{\min}^{\sigma_\hbar^n} -r_{\min} \notin \hbar^{M_{\min}} R_\hbar$. 
\end{proof}

The main theorem of this section establishes the structure of Poisson order on each PI 4-dimensional Sklyanin algebra $S$
for which the Poisson structure on the center of $S$ is nontrivial.

\bth{constr-Pord} Let $S$ be a PI 4-dimensional Sklyanin algebra  satisfying \eqref{albega-cond1}--\eqref{albega-cond2} as above. 
 Each Poisson 
order $(S, Z, \partial: Z \to \Der(S/Z))$ of level~$N$, coming from a good specialization of maximum level $N$,  has the property that the induced Poisson structure on $Z$ is non-zero.
\eth

One interesting direction for further investigation is as follows.

\bqu{max-level} Given  PI 4-dimensional Sklyanin algebra $S$ as used in \thref{constr-Pord}, what is the maximal level $N$ of a good specialization for $S$?
\equ

One approach to this question is to determine the minimal positive integer $M_{\min}$ from the proof 
of \leref{exist} which is explicitly defined in \eqref{MM} and then to use the upper bound 
$N \leq M_{\min}$.

\subsection{Derivations of PI $S$} 
\label{sec:deri}
We classify a certain type of derivations of $S$ which will play a role in the proof of \thref{constr-Pord}. 
This description will be obtained in three stages where similar types of derivations of the algebras in Notation~\ref{nSprime} are classified.

In fact, the techniques below can be generalized easily to obtain similar results for PI algebras $T$, 
where some factor of $T$ by a regular central sequence is the twisted homogeneous coordinate ring of an elliptic curve; see Theorem~\ref{tderT}.

\bnota{Sprime} [$S'$, $\psi$, $\psi_1$, $\psi_2$] For $S$ a PI 4-dimensional Sklyanin algebra of PI degree $n < \infty$ with center $Z := Z(S)$. Recall
$B = S/(g_1 S + g_2 S)$ and denote $S':=S/(g_1 S)$\index{S0pimes@$S'$}. Moreover, denote the canonical projections
$$
\psi_1 : S \to S', \quad \quad
\psi_2 : S' \to B, \quad \quad \mbox{and} \quad
\psi := \psi_2 \psi_1 : S \to B.
\index{p3si@$\psi, \psi_1, \psi_2$}
$$
So, 
\[
\ker \psi_1 = g_1 S,  \quad \quad \ker \psi_2 = g_2 S', \quad \quad \mbox{and} \quad \ker \psi = g_1 S + g_2 S.
\]
\enota

\ble{centers} In the  notation above, we have 
\[
\psi_1(Z) = Z(S') \quad \quad \mbox{and} \quad \quad \psi_2(Z(S')) = Z(B). 
\]
\ele
\begin{proof} Since $\psi(Z) = Z(B)$,
\[
\psi_2(Z(S')) = Z(B). 
\]
The map $\psi_1$ is surjective, thus
\begin{equation}
\label{eq-2cent}
\psi_1(Z) \subseteq Z(S').
\end{equation}
Both sides of the inclusion are $\Nset$-graded algebras with
\[
\psi_1(Z) = \textstyle \bigoplus_{k \in \Nset} \psi_1(Z_k), \quad \quad
Z(S') =  \bigoplus_{k \in \Nset} Z(S')_k.
\]
To prove that \eqref{eq-2cent} is an equality, one needs to prove that
\begin{equation}
\label{cent-k}
\psi_1(Z_k) = Z(S')_k,
\end{equation}
which we show by induction on $k$. 

The equality is obvious for $k=0$. Let $k$ be a positive integer. Assume that \eqref{cent-k} holds for indices less than $k$. 
It follows from the equality  $\psi(Z) = Z(B)$ that for each $z' \in Z(S')_k$, there exists $z \in Z_k$ such that
\[
z' - \psi_1(z) \in \ker \psi_2 = g_2  S'.
\]
Write $z' - \psi_1(z) = g_2 u'$ for some $u' \in S'_{k-2}$. The assumptions on $z$ and $z'$ imply that $g_2 u' \in Z(S')_k$. Therefore, 
$u' \in Z(S')_{k-2}$ because $g_2$ is a regular central element of $S'$ of degree 2.  
It follows from the induction assumption that
$u' = \psi_1(u)$ for some $u \in Z_{k-2}$. Finally, using that $g_2 \in Z_2$, we obtain
\[
z+ g_2 u \in Z(S)_k \quad \quad \mbox{and} \quad \quad
z' = \psi_1(z+ g_2 u) \in \psi_1(Z_k),
\] 
which completes the induction and the proof of the lemma.
\end{proof}

\bnota{adr}[$\ad_r$]
As usual, for an algebra $R$ and $r \in R$, $\ad_r$ \index{adr@$\ad_r$}will denote the inner derivation of $R$ given by $\ad_r(r'):= [r,r']=rr' -r'r$. 
\enota

\bpr{deriv1} \textnormal{[$x$, $\delta$]} Assume that $S$ is a PI 4-dimensional Sklyanin algebra of PI degree $n=|\sigma|<\infty$. 
Let $x \in S_1 \cong B_1$\index{x1@$x$} 
be a good element. 
\begin{enumerate} 
\item 
If $\delta \in \Der B$\index{delta@$\delta$} is such that
\begin{center}
\textnormal{(i)} $\delta|_{Z(B)} =0$, \quad \quad \quad 
\textnormal{(ii)}  $\delta(x) =0$, \quad \quad \quad 
\textnormal{(iii)} $\deg \delta =l$ with $l \leq n$, 
\end{center}
then 
\[
\delta = 
\begin{cases} 
0, & \mbox{if} \; \; l \leq 0 \; \; \mbox{or} \; \; l=n
\\
\lambda \ad_{x^l} , & \mbox{if} \; \; 0< l <n.
\end{cases}
\]
for some $\lambda \in \kk$. 

\medskip

\item If $\delta \in \Der S'$ is such that
\begin{center}
\textnormal{(i)} $\delta|_{Z(S')} =0$, \quad \quad \quad 
\textnormal{(ii)}  $\delta(x) =0$,   \quad \quad \quad 
\textnormal{(iii)} $\deg \delta =l$ with $l \leq n$, 
\end{center}
then $\delta = \ad_{p'}$ for some 
$p' \in \kk[x,g_2]_l$ such that $\deg_x p' < n.$

\medskip

\item If $\delta \in \Der S$ is such that
\begin{center}
\textnormal{(i)} $\delta|_{Z} =0$,  \quad \quad \quad 
\textnormal{(ii)}  $\delta(x) =0$,   \quad \quad \quad 
\textnormal{(iii)} $\deg \delta =l$ with $l \leq n$, 
\end{center}
then $\delta = \ad_p$ for some $p \in \kk[x,g_1,g_2]_l$ such that $\deg_x p < n$.
\end{enumerate}
\epr
\begin{proof} (1) For $l < n$, the statement of this part of the proposition (even without the assumption that 
$\delta|_{Z(B)}=0$) was proved by Smith and Tate in \cite[Theorem 3.3]{SmithTate}. Suppose that $l=n$. Recall that the 
graded quotient ring of $B$ is identified with $\kk(E)[t^{\pm 1}; \sigma]$ by sending $x \mt t$. The derivation 
$\delta \in \Der(B)$ can be extended to a derivation of  $\kk(E)[t; \sigma]$, to be denoted by the same symbol. 
The conditions (i)--(ii) give that
\[
\delta(t) =0 \quad \mbox{and} \quad \delta(\kk(E)^\sigma) =0.
\]
Choose $y \in \kk(E)$ and let $q(t) \in \kk(E)^\sigma$ be its minimal polynomial over $\kk(E)^\sigma$. 
Since the extension $\kk(E)/ \kk(E)^\sigma$ is separable, $q'(y) \neq 0$. 
The assumption that $\deg \delta =n$ implies  that $\delta(y) \in \kk(E) x^n$ commutes with $y$. Therefore
$0 = \delta (q(y)) = \delta(y) q'(y)$ and, hence, $\delta(y)=0$. We have that $\delta(\kk(E))=0$ and $\delta(t) =0$, so $\delta =0$. 

\smallskip

(2) Since $\delta(g_2)=0$, $\psi \delta : S' \to B$ descends to a derivation of $B$, which will be denoted by the same 
symbol.

We prove part (2) by induction on $l$. If $l \leq 0$, then $\psi_2 \delta$ is a homogeneous derivation of $B$ such that 
\begin{center}
$\psi_2 \delta|_{Z(B)} =0$, \quad \quad \quad 
$\psi_2 \delta(x) =0$, \quad \quad \quad  and \quad
$\deg \psi_2 \delta  \leq 0$
\end{center}
because $\psi_2(Z(S')) = Z(B)$ (see Lemma~\ref{lcenterB}). 
Part (1) implies that $\psi_2 \delta =0$. Thus,
$\delta(S'_1) \subseteq \ker \psi_2 = g_2 S'  \subseteq S'_{\geq 2}$.
On the other hand, it follows from $\deg \delta \leq 0$ that
\[
\delta(S'_1) \subseteq S'_{\leq 1}. 
\]
This implies that $\delta(S'_1)=0$ and, consequently that $\delta=0$ because $S'$ is generated in degree 1. 

Let $1 \leq l \leq n$. Assume by way of induction that the statement of part (2) holds for derivations of degree less than $l$. 
Let $\delta$ be a derivation of $S'$ of degree $l \in [1,n]$ such that $\delta|_{Z(S')}=0$ and $\delta(x) =0$.  
As explained earlier, $\psi_2 \delta : S' \to B$ descends to a derivation of $B$ which satisfies
\begin{center}
$\psi_2 \delta|_{Z(B)} =0$, \quad \quad \quad 
$\psi_2 \delta(x) =0$, \quad \quad \quad  and \quad
$\deg \psi_2 \delta  = l \in [1,n]$
\end{center}
because $\psi_2(Z(S')) = Z(B)$. It follows from part (1) that
\[
\psi_2 \delta = \lambda \ad_{x^l} 
\]
for some $\lambda \in \kk$ such that $\lambda =0$ if $l=n$. Therefore,
\[
\text{im} ( \delta - \lambda \ad_{x^l}) \subseteq S' g_2.  
\]
Using that $g_2$ is a regular element of $S'$, we define
\[
\delta^{-} := \textstyle\frac{1}{g_2} (\delta - \lambda \ad_{x^l}) \in \Der(S').
\]
The assumptions on $\delta$ imply at once that
\begin{center}
$\delta^{-}|_{Z(S')} =0$, \quad \quad \quad 
$\delta^{-}(x) =0$, \quad \quad \quad  and \quad
$\deg \delta^{-} =l-2$.
\end{center}
By the inductive assumption there exists $h' \in \kk[x,g_2]_{l-2}$ such that
$\delta^{-} = \ad_{h'}$, and $\deg_x h' <n.$
Therefore, 
\[
\delta = \ad_{p'} \quad \quad \mbox{for} \quad \quad p': = g_2 h' + \lambda x^l \in \kk[x,g_2]_l.
\]
Since $\lambda=0$ if $l=n$, we have $\deg_x p' < n$. This completes the proof of part (2).

\smallskip

Part (3) is derived from part (2) in the same fashion as (2) is derived from (1) making use of the fact that $g_1$ is a 
regular element of $S$.
\end{proof}

Now the techniques above  yield the general statements about derivations of certain PI algebras as follows.

\bnota{notT} [$T$, $\Omega_k$, $T^{[k]}$, $\psi_k$]
Let $T$\index{T0@$T$} be a connected $\mathbb{N}$-graded algebra, and $\Omega_1, \dots, \Omega_m$\index{O0megak@$\Omega_k$} 
a regular central sequence of elements of $T$ of positive degree. Let 
$$T^{[k]} := T/(\Omega_1 T + \cdots + \Omega_k T)
\index{T0k@$T^{[k]}$}$$
for all $1 \leq k \leq m$ so that $T^{[0]} := T$. Note that $T^{[k]} \cong T^{[k-1]}/(\Omega_k T^{[k-1]})$ and let $\psi_k: T^{[k-1]} \twoheadrightarrow T^{[k]}$\index{p3sik@$\psi_k$} 
be the canonical projection.
\enota

Observe that if $T^{[k]}$ is module-finite over its center, then so is $T^{[i]}$ for all $i <k$ by \cite[Lemma~3.6]{SmithTate}.

\bth{derT}
Retain the notation above. Then the following statements hold.
\begin{enumerate}
\item $\psi_k(Z(T^{[k-1]})) = Z(T^{[k]})$  for all $1 \leq k \leq m$. 

\smallskip

\item Suppose that $T^{[m]}$ is isomorphic to $B_T:=B(E_T, \mathcal{L}_T, \sigma_T)$, a twisted homogeneous coordinate ring of an elliptic curve $E_T \subset \mathbb{P}^d$ with $\mathcal{L}_T = \mathcal{O}_{\mathbb{P}^d}(1)|_{E_T}$ and $\sigma_T$ in $\Aut(E_T)$ given by translation. Assume, further, that $B_T$  is module-finite over its center so that $|\sigma_T| < \infty$. 
 Let $x_T$ be a good element of $B_T$ of degree 1. 
 
 If  $\delta$ is a derivation of $T^{[k]}$ such that
\begin{center}
\textnormal{(i)} $\delta|_{Z(T^{[k]})} =0$,  \quad \quad \quad 
\textnormal{(ii)}  $\delta(x_T) =0$,  \quad \quad \quad 
\textnormal{(iii)} $\deg \delta =l$ with $l \leq |\sigma_T|$, 
\end{center}
then $\delta = \ad_p$ for some $p \in \kk[x_T,\Omega_{k+1}, \dots, \Omega_m]_l$ such that $\deg_{x_T} p < |\sigma_T|$. \qed
\end{enumerate}
\eth

\subsection{Proof of \thref{constr-Pord}, nontriviality of Poisson order structure on $S$}
\label{sec:proof-tPord}
Let $\iota : Z \to {S}_\hb$ be a good section of maximal level $N$. Let $\partial: Z \to \Der(S/Z)$ be the corresponding Poisson order. 
We need to prove that the induced Poisson structure on $Z$ is nontrivial. By way of contradiction, assume that this is not the case. Then, for 
$0 \leq i \leq 3$, we have 
\begin{center}
\textnormal{(i)} $\partial_{z_i} |_{Z} =0$, \quad \quad \quad
\textnormal{(ii)} $\partial_{z_i}(x_i) = 0$, \quad \quad  \quad 
\textnormal{(iii)} $\deg \partial_{z_i} =n$.
\end{center}
Property (i) is a restatement of the vanishing of the Poisson structure. Property (ii) follows from the definition of good section, and more precisely the condition in Definition~\ref{dgood-ect}(1.i) and Notation~\ref{nxtilde}. Namely, we have $[\iota(z_i), \wt{x}_i]=0$ and thus, 
\[
\partial_{z_i} (x_i) = \theta ([\iota(z_i), \wt{x}_i]/\hbar^N) =0.  
\]
The third condition is a consequence of the fact that $\iota(z_i) \in ({S}_\hb)_n$. It follows from \prref{deriv1}(3) that there exist polynomials $p_i$ for $0 \leq i \leq 3$
\begin{equation}
\label{pi}
p_i(x_i, g_1, g_2) \in \kk[x_i,g_1,g_2]_n \quad \quad \mbox{and} \quad \quad \deg_{x_i} p_i(x_i, g_1, g_2) < n,
\end{equation}
such that $\partial_{z_i} = \ad_{p_i(x_i, g_1, g_2)}$.
We modify the section $\iota : Z \hra {S}_\hbar$ to form a new good section $\iota\spcheck : Z \hra {S}_\hbar$ by
first setting
\[
\iota\spcheck(g_1) := \wt{g}_1, \quad
\iota\spcheck(g_2) := \wt{g}_2
\quad \text{ and } \quad 
\iota\spcheck(z_i) := \iota(z_i) - \hbar^N p_i (\wt{x}_i, \wt{g}_1, \wt{g}_2)
\]
for $0 \leq i \leq 3$. Then we choose a $\kk$-basis of~$Z$ of the form 
\[
\mathcal{B}:=\{g_1^{m_1} g_2^{m_2} z_0^{l_0} z_1^{l_1} z_2^{l_2} z_3^{l_3} 
\mid (m_1, m_2, l_0, l_1, l_2 , l_3) \in L \}
\]
for some $L \subset \Nset^6$ and complete the definition of $\iota\spcheck$ by setting
\begin{equation}
\label{iota3}
\iota\spcheck \left( g_1^{m_1} g_2^{m_2} z_0^{l_0} z_1^{l_1} z_2^{l_2} z_3^{l_3} \right) = 
\iota\spcheck(g_1)^{m_1} \iota\spcheck(g_2)^{m_2} \iota\spcheck(z_0)^{l_0} 
\iota\spcheck(z_1)^{l_1} \iota\spcheck(z_2)^{l_2} \iota\spcheck(z_3)^{l_3} 
\end{equation}
for all $(m_1, m_2, l_0, l_1, l_2, l_3) \in L$.

The two properties of the polynomials $p_i$ in \eqref{pi} imply that 
$\iota\spcheck : Z \hra {S}_\hbar$ is a good section of $\theta_S$.
Furthermore, we have that 
\begin{equation}
\label{contra}
[\iota\spcheck(z_i), y]  \in \hbar^{N+1} S_\hbar, \quad \quad \forall \; 0 \leq i \leq 3, \; y \in {S}_\hbar.
\end{equation}
This follows from the equalities 
\begin{align*}
&[\iota\spcheck(z_i), y] = \hb^N\left(\frac{[\iota(z_i),y]}{\hbar^N} + [p_i(\wt{x_i}, \wt{g}_1 ,\wt{g}_2), y]\right),
\\
&\theta_S \left(\frac{[\iota(z_i),y]}{\hbar^N} + [p_i(\wt{x_i}, \wt{g}_1 ,\wt{g}_2), y]\right)
= \partial_{z_i}(\theta(y)) - \partial_{z_i}(\theta(y)) = 0
\end{align*}
and the fact that $\ker \theta_S = \hb {S}_\hb$, proved in \leref{Specialization}(4). 

Combining \eqref{iota3} and \eqref{contra}, we obtain that 
\[
[\iota\spcheck(z), y]  \in \hbar^{N+1} {S}_\hbar, \quad \quad \forall z \in \mathcal{B}, \; 0 \leq i \leq 3, \; y \in {S}_\hbar.
\]
This is a contradiction, since $N$ equals the maximum level of a good section of the projection $\theta_S : {S}_\hb \to S$. 
The contradiction implies that the induced Poisson structure on $Z$ from the Poisson order $\partial$ is nontrivial.
\qed

\sectionnew{The Jacobian structure of Poisson orders on PI 4-dimensional Sklyanin algebras} \label{sec:bracket}

The goal of this part is to describe the Jacobian structure of the (nontrivial) Poisson bracket on the center of the PI 4-dimensional Sklyanin algebras that arise as a Poisson order via good specialization of maximal level.
We begin this section with a straight-forward  result about  PI algebras that are Poisson orders via specialization.

\bpr{lemT} Recall Notation~\ref{nnotT}. For a formal parameter $\hbar$, suppose that there exists an $\kk[[\hbar]]$-algebra $T_\hbar$ which is an $\kk[[\hbar]]$-torsionfree degree 0 deformation of $T$ so that $T_\hbar / (\hbar T_\hbar) \cong T$. Take $\theta_T: T_\hbar \twoheadrightarrow T$, the canonical projection, and let $\wt{\Omega}_k \in Z(T_\hbar)$ be a lift of $\Omega_k \in Z(T)$ via $\theta_T$ for all $1 \leq k \leq m$.

Then, for every specialization of $T$ of level $N$ equipped with a section $\iota: Z(T) \hookrightarrow T_\hbar$ with $\iota(\Omega_k) = \wt{\Omega}_k$ for all $1\leq k \leq m$, we have that $\partial_{\Omega_k} = 0$ for the corresponding Poisson order on $T$. In particular, $\Omega_1, \dots, \Omega_m$ are in the Poisson center of $Z(T)$.
\epr

\begin{proof}
For $y \in T$ and $\wt{y} \in \theta_T^{-1}(y)$, we get that $\partial_{{\Omega}_k}(y) = \theta_T\left([\wt{\Omega}_k,y]/\hbar^N\right)$, which is equal to 0. Therefore, $\partial_{\Omega_k} = 0$, and the second statement holds by the definition of a Poisson order. 
\end{proof}

Applying this to the PI 4-dimensional Sklyanin algebras $S$, we have the statement below.

\bco{cor-Z(S)}
The central elements $g_1, g_2$ lie in the Poisson center of $Z=Z(S)$ when $S$ arises as a Poisson $Z$-order of level $N$ via specialization.
\eco 

\begin{proof}
Apply Proposition~\ref{plemT} by taking $T=S$, $T_\hbar = {S_\hbar}$, $\{\Omega_1, \dots, \Omega_m\} = \{g_1, g_2\}$ and by using Lemma~\ref{lSpecialization}(3).
\end{proof}

Next, we have a  general result pertaining to the Poisson structure of the PI algebras $T$ from above. In the following proposition, we use notation $J(\underline{F})_{k,l}$ to denote the determinant of the Jacobian matrix  by taking all $z_i$ derivatives except with respect to variables $z_k$ and $z_l$ up to a sign $(-1)^{k+l}$.

\bpr{ps-Z(T)} \textnormal{[$J(\underline{F})_{k,l}$]} Recall Notation~\ref{nnotT} and we abuse some notation utilized above as follows. 
Suppose that $T$ is module-finite over its center with $Z(T)$ generated by algebraically independent homogeneous elements $z_1, \dots, z_d$, along with homogenous elements $\Omega_1, \dots, \Omega_m$,
subject to $d-2$ homogeneous relations $F_1, \dots, F_{d-2}$.

Suppose that $Y:= \mathbb{V}(F_1, \dots, F_{d-2}) \subset \mathbb{A}^{d+m}$ is an irreducible affine variety and define the Jacobian matrix to be
$$
J(\underline{F})_{k,l}:= (-1)^{k+l}\textnormal{det}\left(\frac{\partial(F_1,F_{2}, \dots,F_{d-3},F_{d-2})}{\partial(z_1,\dots,\widehat{z_k},\dots,\widehat{z_l},\dots,z_d)}\right),\ \text{for $k<l$}
\index{J0Fkl@$J(\underline{F})_{k,l}$}
$$
and $J(\underline{F})_{l,k}=-J(\underline{F})_{k,l}$. Moreover, suppose that $Z(T) = \kk[Y]$ admits a (homogeneous) Poisson structure of degree 0 so that 
\smallskip
\begin{enumerate}
\item[(i)] $\Omega_1, \dots, \Omega_m$ are in the Poisson center,
\item[(ii)] $J(\underline{F})_{k,l}\neq 0$ for some $1\le k,l\le d$,
\item[(iii)] $\bigcap_{1\le k,l\le d} \kk[Y][(J(\underline{F})_{k,l})^{-1}] = \kk[Y]$, where the intersection is taken over $J(\underline{F})_{k,l}\neq$~\hspace{-.1in}~$0$.
\end{enumerate}

\smallskip

\noindent Then, the Poisson bracket on $Z(T)$ is determined by Jacobian matrices above as follows:
\begin{equation}
\label{Jacobian-Pbrack}
\{z_k, z_l\} = \eta J(\underline{F})_{k,l}, \quad \text{ for all $1 \leq k,l \leq d$, and for some}\ \eta\in \kk [Y].
\end{equation}
If, further, $\eta$ has degree 0, then  $\eta\in \kk [Y]_0=\kk$.
\epr

\begin{proof}
By condition (i) we have that
$$\{z_k,F_i\}~=~\sum_{l=1}^d\{z_k,z_l\}\partial F_i/\partial z_l+\sum_{l=1}^m\{z_k,\Omega_l\}\partial F_i/\partial \Omega_l~=~\sum_{l=1}^d\{z_k,z_l\}\partial F_i/\partial z_l~=~0.$$
Now consider the vector space of all $(d-1)$-tuples over the base field $\kk(Y)$, which is the fraction field of the domain $\kk[Y]$. Note that the vector 
$$V=(\{z_k,z_1\},\dots,\widehat{\{z_k,z_k\}},\dots,\{z_k,z_d\})$$ 
is perpendicular to vectors $(\partial F_i/\partial z_1,\dots,\widehat{\partial F_i/\partial z_{k}},\dots,\partial F_i/\partial z_d)$ for all $1\le i\le d-2$. 
By condition (ii) and Cramer's rule, we know $V$ has to be proportional to the vector 
$$
\left(J(\underline{F})_{k,1},\dots, J(\underline{F})_{k,k-1},J(\underline{F})_{k,k+1},\dots,J(\underline{F})_{k,d}\right).
$$
So we can write 
$$\{z_k, z_l\} = \eta_{kl}J(\underline{F})_{k,l},\quad \text{for some $\eta_{kl}\in \kk(Y)$.}$$
 By construction, we get $\eta_{kl}=\eta_{kl'}$ and $\eta_{kl}=\eta_{lk}$. Hence, $\eta_{kl}=\eta\in \kk(Y)$ for all $1\le k, l \le d$. 
Finally, the condition (iii) implies that $\eta \in \kk[Y].$ 
\end{proof}

Now we arrive at the result below.

\bco{ps-Z(S)}
Recall the notation of  Section~\ref{sec:back} and set $\{0,1,2,3\}=\{i,j,k,l\}$ with $i<j$ and $k<l$. Then, the Poisson structure of the center $Z$ of $S$  obtained via good specialization is given by
\begin{align*}
\{z_k, z_l\} &\,= (-1)^{k+l}\eta J(F_1,F_2)_{i,j} \\
&\,= 
\begin{cases}
\bigskip
(-1)^{k+l}\eta \det \begin{pmatrix}
\frac{\partial \Phi_p}{\partial z_q}
\end{pmatrix}_{p=1,2;~q=i,j}, & \text{ for $|\sigma|=:n$ is odd}\\
(-1)^{k+l}\eta \det \begin{pmatrix}
\frac{\partial \Phi_p}{\partial z_q} + h_p \frac{\partial \ell_p}{\partial z_q}
\end{pmatrix}_{p=1,2;~q=i,j}, & \text{ for $|\sigma|=:n$ is even}
\end{cases}
\end{align*}
for some $\eta\in \kk^\times$, with $g_1, g_2$ in the Poisson center of the Poisson algebra $Z$. 
\eco

\noindent {\it Proof}.
We apply Proposition~\ref{pps-Z(T)} for $T=S$ module-finite over its center $Z$. In this case, $\{\Omega_1, \dots, \Omega_m\} = \{g_1, g_2\}$ and $Y = \mathbb{V}(F_1, F_2) \subset \mathbb{A}^6$ from Proposition~~\ref{pcenterS}. The affine variety $Y$ is irreducible because $S$ is a domain (see Proposition~\ref{philb}). First, the degree of $\eta$ is 0 since the degree of $\det(J(F_1,F_2)_{i,j})$ is $2n$, which is equal to the degree of $\{z_k, z_l\}$. Moreover, Theorem~\ref{tconstr-Pord} yields the nontriviality of the Poisson structure on $Z(S)$ via good specialization. So it suffices to verify conditions (i)-(iii) in Proposition~\ref{pps-Z(T)}. First, Corollary~\ref{ccor-Z(S)} yields condition~(i). Next, for condition (ii), note that 
$$Y^{symp}_0 \cap \mathbb{V}(g_1 - \gamma_1, g_2 - \gamma_2)~=~(Y_{\gamma_1,\gamma_2})^{sing}.$$
It follows from Theorem~\ref{tSOdd} ($n$ odd) and Theorem~\ref{tSEven} ($n$ even) that $(Y_{\gamma_1,\gamma_2})^{sing}$ only consists of finitely many points. 
The second equality in \eqref{Yeqeq} implies that there exist indices $0 \leq k \neq l \leq 3$ such that $J(F_1,F_2)_{k,l}$ does not vanish 
identically on $Y$.

Finally, condition~(iii)  of Proposition~\ref{pps-Z(T)} depends on the singular locus on each slice $Y_{\gamma_1,\gamma_2}$ as follows.  We employ work of Stafford \cite{Stafford} to understand the structures of $Y$ and its
coordinate ring $Z$; namely, $S$ is a maximal order, so $Z$ is integrally closed and $Y$ is a normal affine variety. So, to get condition~(iii), it remains to show that 
$Y^{symp}_0=\bigcup_{\gamma_1,\gamma_2 \in \Bbbk} (Y_{\gamma_1,\gamma_2})^{sing}$
has codimension $\geq 2$ in $Y$ due to \cite[discussion after proof of Corollary~11.4]{Eisenbud}. Indeed, we have that $Y^{symp}_0$ is a union of $2n$ cuspidal curves for $n$ odd by Lemma~\ref{lWcusp} and Theorem~\ref{tSOdd} (see Figure~1); thus $Y^{symp}_0$ has codimension $\geq 2$ in $Y$. For $n$ even, condition~(iii) of Proposition~\ref{pps-Z(T)} follows from Theorem~\ref{tSEven} since $Y^{symp}_0=Y^{sing}$ is a union of two rational surfaces $Y_1^{sing}$ and $Y_2^{sing}$, which have codimension $\geq 2$ in $Y$. 
\qed
\medbreak 

Recall that one says that a Poisson structure on an affine variety vanishes at a point $y \in Y$ if the maximal ideal $\mathfrak{m}_y$ of $y$ is a Poisson ideal of $\kk[Y]$.

\bco{YPois} The variety $Y^{symp}_0$ from Notation~\ref{nnot:Y} is precisely the subvariety of $Y=\text{maxSpec}(Z(S))$ consisting of points 
on which the Poisson bracket of Corollary~\ref{cps-Z(S)} vanishes.
\eco
\noindent {\it Proof}. By definition, 
\begin{equation}
Y^{symp}_0 = \bigcup_{\gamma_1, \gamma_2 \in \kk} (Y_{\gamma_1, \gamma_2})^{sing} = 
\left(\textstyle \bigcap_{1 \leq k, l \leq 4} \mathbb{V}(  J(F_1,F_2)_{k,l} ) \right) \cap Y.
\label{Yeqeq}
\end{equation}
Corollary~\ref{cps-Z(S)} implies that the right hand side is precisely the subset of points of $Y$ on which 
the Poisson bracket of Corollary~\ref{cps-Z(S)} vanishes.
\qed
\medbreak

The following result directly follows from the proof of Corollary \ref{cps-Z(S)}, which will be useful in Section \ref{sec:repthy} for determining the Azumaya locus of $S$.

\bco{PoissonZero}
Let $Y={\rm maxSpec}(Z(S))$ be the Poisson variety with the Poisson bracket given in Corollary \ref{cps-Z(S)}. Then, for $(\gamma_1, \gamma_2) \in \Bbbk^2$, we have that $Y^{symp}_0 \cap Y_{\gamma_1,\gamma_2} = (Y_{\gamma_1,\gamma_2})^{sing}$, which consists of finitely many points, and that $Y^{symp}_0$ has codimension $\ge 2$ in $Y$. \qed
\eco

\sectionnew{On the representation theory of PI 4-dimensional Sklyanin algebras $S$} 
 \label{sec:repthy}

    The goal of this section is to use the algebro- and Poisson-geometric results presented in the previous sections to study the irreducible representations of PI 4-dimensional Sklyanin algebras $S$ of PI degree $n< \infty$. Recall that such representations are finite-dimensional and their maximum dimension is equal to $n$ \cite[Proposition~3.1]{BrownGoodearl}. Moreover, the isomorphism classes of irreducible representations (or, simple modules) of $S$ are governed by their central annihilators. We refer the reader to Brown-Goodearl's text \cite[Chapter~III]{BrownGoodearl:book} for details. In summary, there exists a finite-to-1 map on isomorphism classes of simple $S$-modules $[M]$ to the vanishing of $\text{ann}_S(M) \cap Z(S)$ in  $\text{maxSpec}(Z(S))=:Y$.   The {\it Azumaya locus} of $S$ consists of maximal ideals $\mathfrak{m} \in Y$ that annihilate irreducible representations of maximum dimension ($=n$). To study the representation theory of $S$ geometrically, it is advantageous to have that the Azumaya locus of $S$ and the smooth part of $Y$ coincide-- this  is established in Section~\ref{sec:Azumaya}. In Section~\ref{sec:fat} we will then use the singular locus of $Y$ to investigate the irreducible representations of $S$ of intermediate dimension. In Section~\ref{sec:discr} we apply these representation theoretic results to describe the zero sets of the discriminant ideals of the PI 4-dimensional Sklyanin algebras.


\subsection{Irreducible representations of PI $S$ of maximum dimension}
\label{sec:Azumaya}
We will verify the following result, with the proof presented at the end of this section.

\bth{Azumaya} The Azumaya locus of each PI 4-dimensional Sklyanin algebra $S$ is the smooth part of $Y= \maxSpec (Z(S))$.
\eth

\bnota{kappa} [$[\kappa_1:\kappa_2]$, $S_{[\kappa_1:\kappa_2]}$, $g'$, $\psi'$, $\psi''$]
For $[\kappa_1 : \kappa_2 ] \in \Pset^1$, 
\index{kappa@$[\kappa_1: \kappa_2]$}
consider the factors
$$
S_{[\kappa_1:\kappa_2]} := S/ (\kappa_1 g_1 - \kappa_2 g_2) S.
\index{S4kappa@$S_{[\kappa_1:\kappa_2]}$}
$$
Set
$$
g':= 
\begin{cases}
g_1, & \kappa_2 \neq 0
\\
g_2, & \mbox{otherwise}.
\end{cases}
\index{g0prime@$g'$}
$$
The images of $g_1, g_2, g'$ in $S_{[\kappa_1:\kappa_2]}$ will be denoted by the same symbols. 
Denote the projections 
$$
\psi' : S \twoheadrightarrow S_{[\kappa_1:\kappa_2]} \quad \quad \text{and} \quad \quad
\psi'' : S_{[\kappa_1:\kappa_2]} \twoheadrightarrow S/(g_1 S + g_2 S) \cong B.
\index{p4siprime@$\psi', \psi''$}
$$
\enota

Note that the kernel of $\psi'$ is $g' S_{[\kappa_1:\kappa_2]}$. Moreover, the composition $\psi'' \psi'$ 
equals the projection $\psi : S \twoheadrightarrow B$ from Notation~\ref{nSprime} because
$g_1 S + g_2 S = g'S + (\kappa_1 g_1 - \kappa_2 g_2) S.$

Our first preliminary result is given below.

\bpr{3d-quo} For every $[\kappa_1:\kappa_2 ] \in \Pset^1$, the factor $S_{[\kappa_1:\kappa_2]}$ of $S$ is a 
PI domain of PI degree $n = | \sigma |$ and its center is given by
$Z(S_{[\kappa_1:\kappa_2]}) = \psi'( Z(S) )$.
\epr

\begin{proof}
First, $S_{[\kappa_1:\kappa_2]}$ is a domain by \cite[Proposition~6.2]{LevSmith}. Since the algebra $S_{[\kappa_1:\kappa_2]}$ 
is a homomorphic image of $S$, 
we obtain that the PI degree of $S_{[\kappa_1:\kappa_2]}$ is less or equal 
to $n$. From the isomorphism $S_{[\kappa_1:\kappa_2]} /( g' S_{[\kappa_1:\kappa_2]}) \cong B$,
we obtain that the PI degree of $S_{[\kappa_1:\kappa_2]}$ is greater or equal than that 
of $B$ which equals $n$ by \prref{PI}. Thus, $S_{[\kappa_1:\kappa_2]}$ is a PI domain of PI degree~$n$. 
The result on the center of $S_{\kappa_1:\kappa_2]}$ follows from a proof similar to that of Lemma~\ref{lcenters}.
\end{proof}

We turn our attention to maxSpec of the center of $S_{[\kappa_1:\kappa_2]}$ next.

\bnota{Ykappa} [$Y_{[\kappa_1:\kappa_2]}$] 
\index{Y4kappa@$Y_{[\kappa_1:\kappa_2]}$}
Denote by $Y_{[\kappa_1:\kappa_2]}$ the subvariety $Y \cap \mathbb{V}( \kappa_1 g_1 - \kappa_2 g_2)$ of $Y$, which  is the disjoint union of subvarieties $\{ Y_{\gamma_1, \gamma_2} \mid \gamma_1, \gamma_2 \in \kk, ~\kappa_1 \gamma_1 - \kappa_2 \gamma_2 =0 \}$~of~$Y$.
\enota

It also follows from \prref{3d-quo} that
\begin{equation}
\label{Ykappa}
Y_{[\kappa_1:\kappa_2]} \cong \maxSpec (Z(S_{[\kappa_1:\kappa_2]})).
\end{equation}

\bco{irred} The varieties $Y_{[\kappa_1:\kappa_2]}$ and $Y_{\gamma_1, \gamma_2}$
are irreducible for all $[\kappa_1:\kappa_2] \in \Pset^1$ and $\gamma_1, \gamma_2 \in \kk$.
\eco

\begin{proof} The first fact follows from the isomorphism \eqref{Ykappa} and the fact that  $S_{[\kappa_1:\kappa_2]}$ is a domain 
(see \prref{3d-quo}). The set
\[
Y_{[\kappa_1:\kappa_2]}^* : = Y_{[\kappa_1:\kappa_2]} \backslash \mathbb{V}(g') = Y_{[\kappa_1:\kappa_2]} \backslash \mathbb{V}(g_1, g_2)
\]
is dense in $Y_{[\kappa_1:\kappa_2]}$. Thus, it is irreducible too. The $\Nset$-grading of $S$ gives rise to a $\kk^\times$-action 
on $Y$ that preserves $Y_{[\kappa_1:\kappa_2]}$. It is easy to see that 
\[
Y_{[\kappa_1:\kappa_2]}^* \cong \kk^\times \times Y_{\gamma_1, \gamma_2} 
\]
for $\gamma_1, \gamma_2 \in \kk$ such that $\kappa_1 \gamma_1 - \kappa_2 \gamma_2 =0$. This implies that the 
varieties $Y_{\gamma_1, \gamma_2}$ are irreducible for all $\gamma_1, \gamma_2 \in \kk$.
\end{proof}

Next, we study the symplectic cores of $Y$ in terms of the varieties $Y_{\gamma_1, \gamma_2}$. Recall the notation from Section~\ref{sec:bkgdPorder} and consider the following notation.

\bnota{mfrak} [$\mathfrak{m}_y$]
For $y \in Y$, denote by $\mm_y$\index{m1y@$\mm_y$} the corresponding maximal ideal of $Z(S)$. 
\enota

\bpr{scores} Consider a Poisson structure on the algebra $Z$ coming from a Poisson $Z$-order on $S$ 
of level $N$ via a good specialization, as in \thref{constr-Pord}.
\begin{enumerate}
\item For all $\gamma_1, \gamma_2 \in \kk$,  $Y_{\gamma_1, \gamma_2}$ is a Poisson subvariety of $Y = \maxSpec (Z)$. 
\smallskip
\item The symplectic cores of $Y$ are
\begin{enumerate}
\item The points in $(Y_{\gamma_1, \gamma_2})^{sing}$ for $\gamma_1, \gamma_2 \in \kk$ \textnormal{(}0-dimensional symplectic cores\textnormal{)};
\item The sets $Y_{\gamma_1, \gamma_2} \backslash (Y_{\gamma_1, \gamma_2})^{sing}$ for  
$\gamma_1, \gamma_2 \in \kk$ \textnormal{(}2-dimensional symplectic cores\textnormal{)}.
\end{enumerate}
\end{enumerate}
\epr
\begin{proof} (1) The ideal $(g_1 - \gamma_1, g_2 - \gamma_2)$ is a Poisson ideal of $Z$ for all $\gamma_1, \gamma_2 \in \kk$
because $g_1$ and $g_2$ are in the Poisson center of $Z(S)$ by \coref{cor-Z(S)}.
This implies that $Y_{\gamma_1, \gamma_2}$ is a Poisson subvariety of $Y$. 
\smallskip

(2) It follows from Corollary~\ref{cPoissonZero} that for $y \in Y$ the ideal $\mm_y$ is Poisson, if and only if $y \in (Y_{\gamma_1, \gamma_2})^{sing}$ for
some $\gamma_1, \gamma_2 \in \kk$. 

Let $y \in Y_{\gamma_1, \gamma_2} \backslash (Y_{\gamma_1, \gamma_2})^{sing}$ 
for  some $\gamma_1, \gamma_2 \in \kk$. It remains to show that the symplectic core containing $y$ is all of $Y_{\gamma_1, \gamma_2} \backslash (Y_{\gamma_1, \gamma_2})^{sing}$, i.e., that
\begin{equation}
\label{2d-s-cores}
\mathbb{V}(\mathcal{P}(\mathfrak{m}_y)) = Y_{\gamma_1, \gamma_2} \quad \text{ for such $y$}.
\end{equation}
Now \cite[Lemma~6.2]{Goodearl} implies that the Poisson core $\PP(\mm_y)$ is a Poisson prime ideal of $Z(S)$.  We have
\[
\{y\} \subsetneq \mathbb{V}(\PP(\mm_y)) \subseteq Y_{\gamma_1, \gamma_2}
\]
The second inclusion follows from part (1) and the first inclusion follows from the fact that 
$\mm_y$ is not a Poisson ideal for $y \notin (Y_{\gamma_1, \gamma_2})^{sing}$. On the one hand, 
$\mathbb{V}(\PP(\mm_y))$ is an irreducible Poisson subvariety of $Y_{\gamma_1, \gamma_2}$, and on the 
other hand $Y_{\gamma_1, \gamma_2} \backslash (Y_{\gamma_1, \gamma_2})^{sing}$
is a smooth irreducible symplectic variety. This implies that 
\[
\mathbb{V}(\PP(\mm_y)) \backslash (Y_{\gamma_1, \gamma_2})^{sing} = 
Y_{\gamma_1, \gamma_2} \backslash (Y_{\gamma_1, \gamma_2})^{sing}
\]
because a smooth irreducible symplectic variety has no nonempty Poisson subvarieties.
Since $(Y_{\gamma_1,\gamma_2})^{sing}$ consists of finitely many points (see Corollary~\ref{cPoissonZero}) and $\mathbb{V}(\mathcal{P}(\mathfrak{m}_y))$ is an irreducible variety properly containing $\{y\}$, we obtain that $\mathbb{V}(\mathcal{P}(\mathfrak{m}_y)) \not \subset (Y_{\gamma_1,\gamma_2})^{sing}$. Therefore, $\mathbb{V}(\mathcal{P}(\mathfrak{m}_y)) = Y_{\gamma_1,\gamma_2}$ for the element $y \in Y_{\gamma_1,\gamma_2} \backslash (Y_{\gamma_1,\gamma_2})^{sing}$ fixed above. So,~\eqref{2d-s-cores} holds as desired.

The dimensions of the symplectic cores follow from Corollary~\ref{cPoissonZero} and Lemma~\ref{lSingSlice}.
\end{proof}

Before giving a complete description of the Azumaya locus of $S$, we prove that it is sufficiently big in the sense 
that its complement in $Y$ is of codimension $\geq$ 2; see Corollary~\ref{ccorAzumaya} below.

\bpr{Azumaya1} The Azumaya locus of $S$ contains the union of 
\[
Y_{\gamma_1, \gamma_2} \backslash (Y_{\gamma_1, \gamma_2})^{sing} 
\]
for $\gamma_1, \gamma_2 \in \kk$. 
\epr
We first prove the proposition for $\kk = \Cset$ using Poisson geometry, and then extend it to arbitrary algebraically closed 
base fields $\kk$ of characteristic 0.

\begin{proof}[Proof of \prref{Azumaya1} for $\kk= \Cset$] Consider first a pair $(\gamma_1, \gamma_2) \in \Cset^2 \backslash (0,0)$. Fix 
\[
[\kappa_1:\kappa_2] \in \Pset^1 \quad \mbox{such that} \quad \kappa_1 \gamma_1 - \kappa_2 \gamma_2 =0.
\]
For $y \in Y_{[\kappa_1:\kappa_2]}$, denote by $\ol{\mm}_y$ the corresponding maximal ideal of $Z(S_{[\kappa_1:\kappa_2]})$. 
It follows from \prref{3d-quo} that
\begin{equation}
\label{2-factors}
S/\mm_y S \cong S_{[\kappa_1:\kappa_2]}/\ol{\mm}_y S_{[\kappa_1:\kappa_2]}
\quad \mbox{for} \quad y \in Y_{[\kappa_1:\kappa_2]}.
\end{equation}
\prref{scores} implies that the set
\[
Y_{[\kappa_1:\kappa_2]}^{**} ~:=~ \bigsqcup \{ Y_{\gamma'_1, \gamma'_2} \backslash (Y_{\gamma'_1, \gamma'_2})^{sing}
\mid (\gamma'_1, \gamma'_2) \in \Cset^2 \backslash (0,0) \} ~=~ 
\Cset^\times \cdot ( Y_{\gamma_1, \gamma_2} \backslash (Y_{\gamma_1, \gamma_2})^{sing} )
\]
is a single $\Cset^\times$-orbit of symplectic cores of $Y_{[\kappa_1:\kappa_2]} \cong \maxSpec (Z(S_{[\kappa_1:\kappa_2]}))$. 
By the Brown--Gordon theorem, Theorem~\ref{tBrGor}, we obtain that
\begin{equation}
\label{iso-fact}
S_{[\kappa_1:\kappa_2]}/\ol{\mm}_{y} S_{[\kappa_1:\kappa_2]} \cong S_{[\kappa_1:\kappa_2]}/\ol{\mm}_{y'} S_{[\kappa_1:\kappa_2]}
\end{equation}
for all $y, y' \in Y_{[\kappa_1:\kappa_2]}^{**}$. The Azumaya locus of $S_{[\kappa_1:\kappa_2]}$ and $Y_{[\kappa_1:\kappa_2]}^{**}$ 
are both dense subsets of the irreducible variety $Y_{[\kappa_1:\kappa_2]}$. Hence, $Y_{[\kappa_1:\kappa_2]}^{**}$ intersects 
nontrivially the Azumaya locus of $S_{[\kappa_1:\kappa_2]}$. Since the PI degree of $S_{[\kappa_1:\kappa_2]}$ equals $n$ by Proposition~\ref{p3d-quo}, 
it follows from~\eqref{iso-fact} that $S_{[\kappa_1:\kappa_2]}/\ol{\mm}_y S_{[\kappa_1:\kappa_2]} \cong M_n(\Cset)$ 
for all $ y \in Y_{\gamma_1, \gamma_2} \backslash (Y_{\gamma_1, \gamma_2})^{sing}$. 
Finally,~\eqref{2-factors} implies that the Azumaya locus of $S$ contains 
$Y_{\gamma_1, \gamma_2} \backslash (Y_{\gamma_1, \gamma_2})^{sing}$ for $(\gamma_1, \gamma_2) \in \Cset^2 \backslash (0,0)$. 

Next, let $(\gamma_1, \gamma_2)= (0,0)$. It follows from \prref{centerS} that
$Y_{0,0} \cong \maxSpec(Z(B))$.
For $y \in Y_{0,0}$, denote by $\ol{\ol{\mm}}_y$ the corresponding maximal ideal of $Z(B)$. Applying \prref{centerS}, gives
that
\begin{equation}
\label{2-factors-b}
S/\mm_y S \cong B/\ol{\ol{\mm}}_y B
\quad \mbox{for} \quad y \in Y_{0,0}.
\end{equation}
The set $Y_{0, 0} \backslash (Y_{0, 0})^{sing}$ is a single symplectic core of $Y_{0,0}$ and thus 
\[
B/ \ol{\ol{\mm}}_y B \cong B/ \ol{\ol{\mm}}_{y'} B \quad \mbox{for} \quad y, y' \in  Y_{0, 0} \backslash (Y_{0, 0})^{sing}.
\]
On the other hand, $Y_{0, 0} \backslash (Y_{0, 0})^{sing}$ is a dense subset of the irreducible variety $Y_{0,0}$; hence, $Y_{0, 0} \backslash (Y_{0, 0})^{sing}$  
intersects the Azumaya locus of $B$ nontrivially. Since the PI degree of $B$ equals $n$, by \prref{PI}, 
$B/\ol{\ol{\mm}}_y B \cong M_n(\Cset)$ for $y \in Y_{0, 0} \backslash (Y_{0, 0})^{sing}$. By \eqref{2-factors-b},
\[
S/\mm_y  S \cong M_n(\Cset) \quad \mbox{for} \quad y \in Y_{0, 0} \backslash (Y_{0, 0})^{sing}
\]
Therefore, the Azumaya locus of $S$ contains $Y_{0, 0} \backslash (Y_{0, 0})^{sing}$.
\end{proof}

\begin{proof}[Proof of \prref{Azumaya1} for  $\kk$] Fix $y \in Y_{\gamma_1, \gamma_2} \backslash (Y_{\gamma_1, \gamma_2})^{sing}$ for some $\gamma_1, \gamma_2 \in \kk$. 
Denote by $\kk'$ the algebraically closed subfield of $\kk$ generated by the coordinates of $y$ and $\alpha, \beta, \gamma \in \kk$. 
In particular, $\gamma_1, \gamma_2 \in \kk'$. Choose a field embedding $\kk' \hookrightarrow \Cset$.

Let $S_\kk=S$, $S_{\kk'}$, and $S_{\Cset}$ be the corresponding Sklyanin algebras over $\kk$, $\kk'$, and $\Cset$, respectively. Let 
$Y_\kk=Y$, $Y_{\kk'}$, and $Y_\Cset$ be the max-spectra of their centers, respectively. Denote by
$\mm_{y, \kk}= \mm_y$, $\mm_{\kk',y}$ and $\mm_{\Cset, y}$, the maximal ideals of 
 $Z(S_\kk)$, $Z(S_{\kk'})$, and $Z(S_{\Cset})$, corresponding to $y$ considered as a point on 
 $Y_\kk=Y$, $Y_{\kk'}$, and $Y_\Cset$, respectively. We have 
 \begin{equation}
\label{ext}
(S_{\kk'}/ \mm_{\kk', y} S_{\kk'}) \otimes_{\kk'} \kk \cong S_\kk / \mm_{\kk, y} S_\kk \quad \mbox{and} \quad 
(S_{\kk'}/ \mm_{\kk', y} S_{\kk'}) \otimes_{\kk'} \Cset \cong S_\Cset / \mm_{\Cset, y} S_\Cset.
\end{equation}

Since $y \in Y_{\gamma_1, \gamma_2} \backslash (Y_{\gamma_1, \gamma_2})^{sing}$, we get that
$y \notin (Y_{\kk', \gamma_1, \gamma_2})^{sing}$. So,
$y \notin (Y_{\Cset, \gamma_1, \gamma_2})^{sing}$. Applying the proposition for $\kk = \Cset$ established above, we obtain that 
\[
S_\Cset / \mm_{\Cset, y} S_\Cset \cong M_n(\Cset)
\]
since all three algebras $S_\kk=S$, $S_{\kk'}$, and $S_{\Cset}$ have PI degree $n$. 
Now it follows from~\eqref{ext} that $S_{\kk'} / \mm_{\kk', y} S_{\kk'} \cong M_n(\kk')$ because $\kk'$ is algebraically closed. 
Applying~\eqref{ext} again gives $S/\mm_y S = S_\kk / \mm_{\kk, y}  S_\kk \cong M_n(\kk)$, which completes the proof of 
the proposition for the field~$\kk$.
\end{proof}

Now it follows from Proposition~\ref{pAzumaya1} and Corollary~\ref{cPoissonZero} that:

\bco{corAzumaya} The complement of the Azumaya locus $\mathcal{A}$ of $S$ in $Y$ is of codimension $\geq 2$, that is, $\mathcal{A}_\mathfrak{p}$ is Azumaya over $Z_\mathfrak{p}$ for all height 1 primes $\mathfrak{p}$ of $Z$. \qed
\eco

Finally, the proof of the main result of this section is  brief, due to the work above.
 
 \medskip
 \noindent {\it Proof of Theorem~\ref{tAzumaya}}. This holds by applying a result of  Brown-Goodearl \cite[Theorem~3.8]{BrownGoodearl}; namely, the hypotheses of that result follow from Proposition~\ref{philb} and Corollary~\ref{ccorAzumaya}. 
Alternatively, we can apply the improvement of the theorem obtained by Brown-MacLeod \cite[Theorem 3.13]{BrownMacLeod}.   
 \qed

\medskip

Now the consequence of Theorem~\ref{tAzumaya} below follows from Theorem~\ref{tSOdd}(2) and Theorem~\ref{tSEven}(2); namely, the slice $Y_{0,0} = Y \cap \mathbb{V}(g_1, g_2)$ has singular locus $\{\underline{0}\}$. 

\bco{repB}
All nontrivial irreducible representations of the twisted homogeneous coordinate ring $B$ of PI degree $n$ have dimension $n$. \qed
\eco

\subsection{Irreducible representations of PI $S$ of intermediate dimension} \label{sec:fat}
Throughout this subsection, we work over complex numbers $\kk=\mathbb C$ (as mentioned in Remark~\ref{rrem:simple}). The goal of this part is to use the geometry of fat points of $S$ to classify irreducible representations of $S$ of intermediate dimension. Recall that {\it fat points} of $S$ are simple objects in the quotient category $S$-qgr. A fat point over $S$ can be represented by a $1$-critical graded $S$-module (called \emph{fat point module}) of Hilbert series $d(1-t)^{-1}$ with multiplicity $d\ge 1$. In particular, when a fat point has multiplicity 1, we usually call it a {\it point} which is represented by a {\it point module} of $S$. We employ work of Smith \cite{Smith1993}  on the geometry of fat points of $S$ to achieve our goal.

\bth{thm:fat}\textnormal{[$M(p)$, $F(\omega+k\tau)$]} \cite{Smith1993} All fat points of $S$ have been classified and they come in one of $3$-parametric families:
\begin{enumerate}
\item for each $p\in E$, there is a point module $M(p)$\index{m@$M(p)$};
\item for any $\omega\in E_2$ and $0\le k\le s-1$, there is a fat point module denoted by $F(\omega+k\tau)$\index{F0omega@$F(\omega+k\tau)$} 
of multiplicity $k+1$; 
\item all other fat point modules over $S$ have multiplicity $s=n/(n,2)$. 
\end{enumerate}
In particular, the fat point $F(\omega+k\tau)$ lies on all the secant lines $\ell_{pq}$ of $E$ such that $p+q=\omega+k\tau$; namely, there exists a short exact sequence of graded $S$-modules 
\begin{align}\smaller\label{FSES}
\hspace{-.15in}
0 \longrightarrow M(p-(k+1)\tau,q-(k +1)\tau)[-k -1] \longrightarrow  M(p,q) \longrightarrow  F(\omega+k \tau) \longrightarrow  0.
\end{align}
\qed
\eth

Note that when $k=0$, the four point modules $F(\omega)$ correspond to the four points $\{e_i\}_{0\le i\le 3}$ in the point scheme of $S$.

Next, we recall some facts about the relationship between fat points and irreducible representations  of $S$. 
\bre{rem:simple} [$V(\omega+k \tau)$] Constructed originally by Sklyanin in \cite{Skly1983} for each $\omega\in E_2$ and $k\in \mathbb N\cup \{0\}$, there exists a representation $V(\omega+k\tau)$ \index{V0omega@$V(\omega+k\tau)$}
over $S$ in a certain $(k+1)$-dimensional subspace of theta functions of order $2k$. (We require $\kk = \mathbb{C}$ here.) These representations were proved later by Smith and Staniszkis to be irreducible whenever $k<s$ \cite[Theorem 3.6]{SmithStan}. According to \cite[Proposition~3.3 and Section 6]{Smith1993}, all fat point modules $F(\omega+k\tau)$ of $S$ of multiplicity $<s$ arise as
\begin{align}\label{FatSimple}
F(\omega+k\tau)= \textstyle \bigoplus_{i\ge 0} V(\omega+k\tau)t^i,
\end{align}
where the graded $S$-action on the right side is given by $s_i\cdot v_jt^j=(s_iv_j)t^{i+j}$ for all $s_i\in S_i$ and $v_j\in V(\omega+k\tau)$. Here, the multiplicity of $F(\omega+k\tau)$ is equal to the dimension of $V(\omega+k\tau)$ which is equal to $k+1$. 
\ere

\bqu{rep-dims}
We can  extend \eqref{FatSimple} for $k=s-1$, where $F(\omega+(s-1)\tau)$ is a fat point module of multiplicity $s$. This yields 4 irreducible representations $V(\omega+(s-1)\tau)$ of dimension $s$ for $\omega\in E_2$. What are the forms of the other  irreducible representations of $S$ of dimension $s$?
\equ

Moreover, we have the following fact. Recall that $S$ is assumed to be module-finite over its center, so all of its irreducible representations are finite-dimensional. 

\ble{lem:quot} \cite[Lemma~4.1]{SmithStan}\cite[Sectiosn~3.3 and~5]{Smith1993}
Every irreducible representation $V$ of $S$ is the quotient of some fat point module $F$ of $S$. Here, $\dim_\kk(V) \geq \text{mult}(F)$. Moreover, if $F$ is isomorphic to a shift $F[e]$ for some $e \geq 1$, then $\dim_\kk(V)\leq e \cdot \text{mult}(F)$. \qed
\ele

\bnota{TwistM}[$V^\lambda$]
Let ${\rm Mod}_d(S)$ be the variety of all $d$-dimensional representations of $S$. The algebraic group ${\rm PGL}_d(\kk)\times \kk^\times$ acts on ${\rm Mod}_d(S)$ via $$((T,\lambda).\varphi_V)(a):=\lambda^iT\varphi_V(a)T^{-1}$$ for any $V\in {\rm Mod}_d(S)$ with corresponding map $\varphi_V: S\to {\rm End}(V)$ and $a\in S_i$. For any $V\in {\rm Mod}_d(S)$ and $\lambda\in \kk^\times$, we define the {\it twisted representation} $$V^\lambda:=(1,\lambda).V.$$ It is important to point out that $V\cong V^\lambda$, 
if and only if there is some $T\in {\rm PGL}_d(\kk)$ such that $(T,\lambda).V=V$.
\enota

\bpr{FatQS}  
Let $F$ be a fat point module over $S$ of multiplicity $d$, and $V$ be a nontrivial simple quotient of $F$. Then we have 
\begin{enumerate}
\item $\dim V=de$, where $e$ is the period of $F$ such that $F\cong F[e]$ and $F\not\cong F[i]$ for $1\le i\le e-1$ in $S$-qgr;
\item the stabilizer of $V$ in ${\rm PGL}_{de}(\kk)\times \kk^\times$ is conjugate to the subgroup generated by $(g_\zeta,\zeta)$ with $g_\zeta={\rm diag}(\underbrace{1,\dots,1}_{d},\underbrace{\zeta,\dots,\zeta}_{d},\dots,\underbrace{\zeta^{e-1},\dots,\zeta^{e-1}}_{d})$ and $\zeta$ is a primitive $e$-th root of unity;
\item all nontrivial simple quotients of $F$ are in the form of $V^{\lambda}$, for some irreducible representation $V$ of $S$ and for some $\lambda\in \kk^\times$. \qed 
\end{enumerate}
\epr
\begin{proof}
Part (1) follows from \cite[Proposition 6 and its proof]{LeBruyn}.
Part (2) comes from \cite[Lemma~4]{BS}.
Part (3) holds by the proof of \cite[Theorem~3.13]{Smith1993}.
\end{proof}

We now turn our attention to the relationship between fat points and Azumaya locus of $S$, and for this we need some results from noncommutative projective algebraic geometry; see \cite{Artin:geometry} \cite[Section~5]{SmithTate} for further details. Let $$\mathcal R:={\rm Proj}(Z(S)).$$ We consider the sheaf $\mathcal S$ of $\mathcal O_\mathcal R$-algebras defined by $\mathcal S(\mathcal R_{(z)})=S[z^{-1}]_0$ for any non-zero homogeneous element $z\in Z(S)$.  Denote by $\mathcal Z$ the {\it central Proj} of $\mathcal S$, which satisfies  $\mathcal Z(\mathcal R_{(z)}) = Z(S[z^{-1}]_0)$ for any non-zero homogeneous element $z\in Z(S)$.  Consider the commutative scheme ${\bf Spec}(\mathcal Z)$ defined in \cite[Chapter~II, Exercise~5.17]{Hartshorne}, and we get by \cite[Theorem 5.2]{SmithTate} that it is isomorphic to ${\rm Proj}(Z(S^{(n/s)}))$. 
Now fat point modules of $S$ have support in ${\bf Spec}(\mathcal Z)$; indeed, a fat point module $F$ of $S$ gives rise to a sheaf $\mathcal F$ of $\mathcal S$-modules  such that $\mathcal F(\mathcal R_{(z)})=F[z^{-1}]_0$.

We get that there is a quotient morphism 
\begin{equation}\label{eq:SpecZ}
{\bf Spec}(\mathcal Z) \twoheadrightarrow \mathcal R,
\end{equation} induced by the embedding $Z(S)\hookrightarrow Z(S^{(n/s)})$ via $z_i=u_i^{(n/s)}$ and $g_j=g_j$. 
We obtain that the support of a fat point module $F$ in $\mathcal R$ is the image of the support of the corresponding sheaf of $\mathcal S$-modules $\mathcal F$ in ${\bf Spec}(\mathcal Z)$ via \eqref{eq:SpecZ}. 

Now, with \eqref{EOmega} and a result in \cite{SmithTate}, we obtain the following results on the fat point modules of $S$ of multiplicity $<s$.

\bpr{NAS}\cite[Theorem 5.11]{SmithTate} When $0 \leq k \leq s-1$, 
the support of the fat point modules $F(\omega+k\tau)$ and $F(\omega'+k'\tau)$ coincide in ${\bf Spec}(\mathcal Z)$ if $\omega+k\tau=-\omega'-k'\tau-2\tau$. \qed
\epr

\bnota{scriptC}[$\mathcal{C}(\omega+k\tau)$] Recall that for any $\omega\in E_2$ and $0\le k\le s-1$, the simple module $V(\omega+k\tau)$ constructed by Sklyanin  gives a non-zero point $(p_0,p_1,p_2,p_3,\gamma_1,\gamma_2)$ in $Y$ via its support, or rather central annihilator, in $Z(S)$. Consider the parametric curve 
$$
\mathcal C(\omega+k\tau) \index{c1@$\mathcal{C}(\omega+k\tau)$} := \left\{(t^np_0,t^np_1,t^np_2,t^p_3,t^2\gamma_1,t^2\gamma_2) ~|~ t\in \kk \right\}.
$$ 
Since scaling of the parameter $t$ in $\mathcal C(\omega+k\tau)$ respects the degree of the variables, the curve $\mathcal C(\omega+k\tau)$ represents a single point in $\mathcal R={\rm Proj}(Z(S))$. 
\enota

\ble{FatSing}
The parametric curve $\mathcal C(\omega+k\tau)\subset Y$ is the support of $F(\omega+k\tau)$ in $\mathcal R$. As a consequence, $\mathcal C(\omega+k\tau)=\mathcal C(\omega'+k'\tau)$ if $\omega+k\tau=-\omega'-k'\tau-2\tau$. Moreover, there is a $(n/s)$-to-$1$ correspondence between the isomorphism classes of nontrivial simple quotients of $F(\omega+k\tau)$ and the points of $\mathcal C(\omega+k\tau)$.
\ele
\begin{proof}
By Proposition \ref{pcenterS}, $Z(S)=\kk[z_0,z_1,z_2,z_3,g_1,g_2]/(F_1,F_2)$ is a graded algebra with ${\rm deg}(z_i)=n$ and ${\rm deg}(g_j)=2$. Hence any scaling of a point $(p_0,p_1,p_2,p_3,\gamma_1,\gamma_2)\in Y$ with respect to the grading still belongs to $Y$. So $\mathcal C(\omega+k\tau)\subset Y$. Now let $\mathcal I$ be the smallest graded ideal in $Z(S)$ containing $\mathcal C(\omega+k\tau)$. By \eqref{FatSimple}, one can check that $\mathcal I={\rm Ann}_{Z(S)}(F(\omega+k\tau))$. So $\mathcal C(\omega+k\tau)$ represents the support of $F(\omega+k\tau)$ in $\mathcal R$ via~\eqref{eq:SpecZ}.

From the definition of $F(\omega+k\tau)$ in \eqref{FatSimple}, we know $V(\omega+k\tau)$ is a simple quotient of $F(\omega+k\tau)$. By Proposition \ref{pFatQS}(3), all the nontrivial simple quotients of $F(\omega+k\tau)$ are given by twisted modules $V(\omega+k\tau)^\lambda$ for some $\lambda\in \kk^\times$. Moreover, it is clear that $F(\omega+k\tau)\cong F(\omega+k\tau)[1]$ in $S$-qgr. So $F(\omega+k\tau)$ has period $1$. (See, also \cite[proof of Corollary~8.7 and Proposition~8.1(c)]{Smith1993}). As a consequence, ${\rm PGL}_{k+1}(\kk)\times \kk^\times$ acts freely on the simple quotients of $F(\omega+k\tau)$ by Proposition \ref{pFatQS}(2). So $V(\omega+k\tau)\not\cong V(\omega+k\tau)^\lambda$ for any $\lambda\in \kk^\times$. Let $(p_0,p_1,p_2,p_3,\gamma_1,\gamma_2)\in Y$ be the point corresponding to $V(\omega+k\tau)$. One can check that the point corresponding to $V(\omega+k\tau)^\lambda$ is given by $(\lambda^np_0,\lambda^np_1,\lambda^np_2,\lambda^np_3,\lambda^2\gamma_1,\lambda^2\gamma_2)$. So $V(\omega+k\tau)^\lambda$ and $V(\omega+k\tau)^{\lambda'}$ share the same central annihilator, if and only if $\lambda=\lambda'$ if $n$ is odd and $\lambda=\pm \lambda'$ if $n$ is even. This establishes the $(n/s)$-to-$1$ correspondence.

Finally, $\mathcal C(\omega+k\tau)=\mathcal C(\omega'+k'\tau)$ if $\omega+k\tau=-\omega'-k'\tau-2\tau$ by Proposition \ref{pNAS} since they represent the same point in $\mathcal R$. 
\end{proof}

This brings us to the main results of this section, first for $n$ odd.

\bth{OddS}
When the PI degree $n$ of $S$ is odd, we have the following statements.
\begin{enumerate}
\item The curve $\mathcal C(\omega+k\tau)$ is equal to the curve $C(\omega+k\tau)$ from Lemma~\ref{lWcusp}, for $\omega \in E_2$ and $0\leq k \leq n-2$. Thus, $Y^{sing} = \bigcup_{\omega \in E_2, ~0 \leq k \leq n-2} \mathcal C(\omega+k\tau)$.
\smallskip

\item If $\underline{0}\neq y\in Y^{sing}$ and $0 \leq k \leq n-2$, then there are exactly 2 irreducible representations of $S$ with central annihilator $\mathfrak{m}_y$. 
They are simple quotients of $F(\omega+k\tau)$ of dimension $k+1$ and of $F(\omega+(n-2-k)\tau)$ of dimension $n-1-k$, 
when $y\in C(\omega+k\tau)\backslash \{ \underline{0} \}$. Moreover, the origin $\underline{0}\in Y$ corresponds uniquely to the trivial module of $S$.
\smallskip

\item If $y\in Y^{smooth}$, then there is exactly 1 irreducible representation of $S$ of dimension $n$ with central annihilator $\mathfrak{m}_y$, which is a simple quotient of either the point module $M(p)$ for some $p\in E$ or some fat point module of multiplicity $n$.
\smallskip

\item $S$ has irreducible representations of each dimension $1,2,3,\dots,n$, where the nontrivial irreducible representations of intermediate dimension $k+1$, for $0\le k\le n-2$, are given by $V(\omega+k\tau)^\lambda$ for $\omega\in E_2$ and $\lambda\in \kk^\times$.
\end{enumerate}
\eth

\begin{proof}
(1) It suffices to show for all $\omega\in E_2$ and $0\le k\le n-2$ that  the curves $\mathcal C(\omega+k\tau)$ and $C(\omega+k\tau)$  share a nonzero common point. Indeed, since the curves are parametrized in the same way, sharing a nonzero point implies that the two sets are equal to each other. Now let $\underline{0} \neq y\in \mathcal C(\omega+k\tau)$ correspond to a nontrivial simple quotient $V(\omega+k\tau)$ of a fat point module $F(\omega+k\tau)$ of intermediate multiplicity. Since $\dim V(\omega+k\tau)=k+1$, which is $<n$ by Remark~\ref{rrem:simple}, the point $y$ lies in the non-Azumaya locus and hence the singular locus of $S$ by Theorem~\ref{tAzumaya}. Then \eqref{FSES} and Lemma \ref{lWcusp} imply that $y\in Y^{sing}\cap \mathbb{V}(\Omega(\omega+k\tau))=C(\omega+k\tau)$.

\smallskip

(2) Let $\underline{0}\neq y\in \mathcal C(\omega+k\tau)=\mathcal C(\omega+(n-2-k)\tau)\subset Y^{sing}$ by Lemma \ref{lFatSing} for some $\omega\in E_2$ and $0\le k\le n-2$. Then there are 2 nontrivial irreducible representations $V_1$ and $V_2$ of $S$ with central annihilator $\mathfrak{m}_y$ such that $V_1$ is a simple quotient of $F(\omega+k\tau)$ and $V_2$ is a simple quotient of $F(\omega+(n-2-k)\tau)$. It is clear that $V_1\not\cong V_2$ since $n$ is odd, $\dim V_1=\dim V(\omega+k\tau)=k+1$ and $\dim V_2=\dim V(\omega+(n-2-k))=n-1-k$ by Remark~\ref{rrem:simple}. Now if $\{V_i\}_{1\le i\le t}$ is a complete set of all non-isomorphic irreducible representations  over $S$ central annihilator $\mathfrak{m}_y$, then applying a result of Braun \cite[Proposition 4]{Braun} yields
$$\text{PIdeg}(S)~=~n~\ge~ \sum_{1\le i\le t} \dim V_i~\ge~ \dim V_1+\dim V_2~=~(k+1)+(n-1-k)~=~n.$$
This implies that $t=2$, as desired. 

The statement for the correspondence between the trivial module and the origin of $Y$ is clear. 

\smallskip
(3) This can be proved similarly to (2), again using Theorem~\ref{tAzumaya}.
\smallskip

(4) This follows from (2), (3) and by Remark~\ref{rrem:simple}.
\end{proof}

For $n$ even case, the picture for irreducible representations of $S$  of intermediate dimension versus singular locus of $Y$ is somewhat elusive  due to the fact that we do not know how the varieties $Y^{sing}_1$ and $Y^{sing}_2$ in Theorem~\ref{tSEven} compare to the curves $\mathcal{C}(\omega + k\tau)$. But we  have the following result.

\bth{even}
Suppose that the PI degree $n=2s$ of $S$ is  even. Then the following statements hold.  
\begin{enumerate}
\item[(1)] If $\underline{0}\neq y\in \mathcal C(\omega+k\tau)$ for some $\omega\in E_2$ and $0\le k\le s-2$, then there are exactly 4 irreducible representations of $S$ with central annihilator $\mathfrak{m}_y$, two of which are simple quotients of $F(\omega+k\tau)$ of dimension $k+1$ and the other two are simple quotients of $F(\omega+s\tau+(s-k-2)\tau)$ of dimension $s-1-k$.
\smallskip

\item[(2)] If $\underline{0}\neq y\in \mathcal C(\omega+(s-1)\tau)$ for some $\omega\in E_2$, then there are exactly 2  irreducible representations of $S$ of dimension $s$ with central annihilator $\mathfrak{m}_y$, which are simple quotients of $F(\omega+(s-1)\tau)$.
\smallskip

\item[(3)] If $\underline{0}\neq y\in Y^{sing}\setminus \bigcup_{\omega\in E_2,0\le k\le s-1} \mathcal C(\omega+k\tau)$, then there are exactly 2 irreducible representations of $S$ of dimension $s$ with central annihilator $\mathfrak{m}_y$, which are simple quotients of a fat point module of multiplicity $s$.
\smallskip

\item[(4)] If $y\in Y^{smooth}$, then there is exactly 1 irreducible representation of $S$ of dimension $n$ with central annihilator $\mathfrak{m}_y$, which is a simple quotient of a generic fat point module of multiplicity $s$. 
\smallskip

\item[(5)] $S$ has irreducible representations of each dimension $1,2,3,\dots,s$ and $n$, where the nontrivial irreducible representations of intermediate dimension $k+1$, for $0\le k\le s-2$, are given by $V(\omega+k\tau)^\lambda$ for $\omega\in E_2$ and $\lambda\in \kk^\times$.
\end{enumerate}
\eth

\begin{proof}
(1) Let $\underline{0}\neq y\in \mathcal C(\omega+k\tau)=\mathcal C(\omega+s\tau+(s-k-2)\tau)$. Then there are two non-isomorphic irreducible representations $V$ and $W$,  with central annihilator $\mathfrak{m}_y$, which are simple quotients of $F(\omega+k\tau)$ and $F(\omega+s\tau+(s-k-2)\tau)$, respectively. Moreover, by Lemma~\ref{lFatSing} there exists another simple quotient $V^{-1}$ (respectively, $W^{-1}$) of $F(\omega+k\tau)$ (respectively, of $F(\omega+s\tau+(s-k-2)\tau)$) which also corresponds to $y$ and is not isomorphic to $V$ (respectively, to $W$). So we have in total four non-isomorphic simple modules $V^{\pm 1}, W^{\pm 1}$ corresponding to the same point $y$. Now applying Braun's result \cite[Proposition 4]{Braun} to get that these are the only four irreducible representations with central annihilator $\mathfrak{m}_y$. 

\smallskip

(2), (3) and (4) can be proved similarly to part (1). 

\smallskip

(5) Let $V$ be a nontrivial irreducible representation over $S$, which is a simple quotient of some fat point $F$ by  Lemma~\ref{llem:quot}. If $F=M(p)$ for some $p\in E$, then $V$ is a module over the twisted homogeneous coordinate ring $B$ and $\dim V=n$ by \cite[Lemma~5.8(c)]{LevSmith} and Proposition~\ref{pFatQS}(3). Suppose $F=F(\omega+k\tau)$. Then $\dim V=\dim V(\omega+k\tau)=k+1$ by Remark~\ref{rrem:simple}. Finally, if $V$ is a simple quotient of a fat point of multiplicity $s$. Then by Proposition \ref{pFatQS}(1), we have $\dim V=es$, where $e$ is the period of $F$. By \cite[Proposition 3.1(a)]{BrownGoodearl}, we know $\dim V\le {\rm PI}{\rm deg}(S)=n$. So we have $e=1,2$. Now our result follows from Proposition \ref{pFatQS}(3). 
\end{proof}


\subsection{Discriminant ideals of the PI 4-dimensional Sklyanin algebras}
\label{sec:discr}
Discriminant ideals of PI algebras play an important role in the study of maximal orders \cite{Reiner}, the automorphism and isomorphism problems for families of algebras \cite{CPWZ}, 
the Zariski cancellation problem \cite{BellZhang}, and the description of dimensions of irreducible representations \cite{BrownYakimov}. 

Let $A$ be an algebra and $C \subseteq Z(A)$ be a central subalgebra. A {\em{trace map}} on $A$ is a nonzero map $\tr \colon A \to C$ which is cyclic 
($\tr( xy) = \tr(yx)$ for $x, y \in A$) and $C$-linear. For a positive integer $\ell$, 
the $\ell$-th {\em{discriminant ideal}} $D_\ell(A/C)$ and the $\ell$-th {\em{modified discriminant ideal}} $MD_\ell(A/C)$
of $A$ over $C$ are the ideals of $C$ with generating sets
\begin{align*}
&\{\det( [\tr(y_i y_j)]_{i,j=1}^\ell) \mid y_1, \ldots, y_\ell \in A\} \; \;  \mbox{and}
\\
&\{ \det( [\tr(y_i y'_j)]_{i,j=1}^\ell) \mid y_1, y'_1, \ldots, y_\ell, y'_\ell \in A \}.
\end{align*}

Every maximal order $A$ in a central simple algebra admits the {\em{reduced trace map }} $\tr : A \to Z(A)$, see \cite[Section~9]{Reiner}.
Stafford \cite{Stafford} proved that the PI 4-dimensional Sklyanin algebras $S$ are maximal orders in central simple algebras. 

The next theorem describes the zero sets of the discriminant ideals of a PI Sklyanin algebra $S$ of PI degree $n$. Denote the quadric function
$$
q(k) := (k^2 + (s-k)^2) n/s,
$$
keeping in mind that $s =n$ if $n$ is odd, and $s = n/2$ if $n$ is even.
\bth{disc} Let $S$ be a PI 4-dimensional Sklyanin algebra of PI degree $n$ with reduced trace map $\tr : S \to Z$. For all positive integers 
$\ell$, the zero sets of the $\ell$-th discriminant and $\ell$-th modified discriminant ideals of $S$ over its center coincide,
\[
\mathbb{V}(D_\ell(S/Z(S)), \textnormal{tr})=\mathbb{V}(MD_\ell(S/Z(S)), \tr);
\]
denote this set by $\mathbb{V}_\ell \subset Y : = {\rm maxSpec}(Z)$. For the top discriminant ideal of $S$, we have
\[
\mathbb{V}_{n^2} = Y^{sing}.
\]
For the lower level discriminant ideals of $S$ and base field $\kk = \Cset$, the following hold:
\smallskip

(1) If $n$ is odd, then
$$
{\mathbb{V}}_{\ell}
= 
\begin{cases}
Y^{sing}, & \ell \in \left[ q \left(n-1\right) + 1,~~ n^2 \right]
\\
\bigcup_{k \in [0,n-2], q(k+1)< \ell} C(\omega + k \tau), & \ell \in \left[ q(\lfloor \frac{n}{2}\rfloor) +1 , ~~ q \left(n-1\right) \right]
\\
\underline{0}, & \ell \in \left[ 2,~~ q (\lfloor \frac{n}{2}\rfloor) \right]
\\
\varnothing, & \ell = 1.
\end{cases}
$$

(2) If $n$ is even, then
$$
{\mathbb{V}}_{\ell}
= 
\begin{cases}
Y^{sing}, & \ell \in \left[ 2s^2+1,~~ n^2 \right]
\\
 \bigcup_{k \in [0,s-2], q(k+1)<  \ell} \mathcal{C}(\omega + k \tau) , & \ell \in \left[q \left( \lfloor \frac{s}{2} \rfloor \right) +1,~~ 2s^2 \right]
\\
\underline{0}, & \ell \in \left[2,~~ q \left( \lfloor  \frac{s}{2} \rfloor \right) \right]
\\
\varnothing, & \ell =1.
\end{cases}
$$
\qed
\eth
\begin{proof} 
Let $A$ be a maximal order in a central simple algebra over a field of characteristic~0 with reduced trace map $\tr : A \to Z(A)$.
For $\mm \in {\rm maxSpec}(Z)$, denote by 
$\Irr_\mm(A)$ the set of isomorphism classes of irreducible representations of $A$ with central annihilator $\mm$, and set
$$
d(\mm) := \sum_{V \in \Irr_\mm(A)} (\dim_\kk V)^2.
$$
By 
\cite[Main Theorem (a),(e)]{BrownYakimov} and Theorem~\ref{tAzumaya}, we have that $\mathbb{V}_{n^2} = Y^{sing}$, and for all $\ell \in \mathbb{Z}$ that
\begin{equation}
\label{dis}
 \mathbb{V}_\ell
=  \Big\{ \mm \in \mathrm{maxSpec} Z(A) \mid
d(\mm)  < \ell \Big\}.
\end{equation}
Theorems \ref{tOddS} and \ref{teven} imply that for the PI 4-dimensional Sklyanin algebra $S$, 
the square-dimension function $d :  {\rm maxSpec}(Z) \to \Zset$ is given by
$$
d(\mm_y) = 
\begin{cases}
n^2, & y \in Y^{smooth}
\\
q(k+1), & y \in C(\om + k \tau) \backslash \{ \underline{0} \}, ~k \in [0, n-2]
\\
1, & y = \underline{0}
\end{cases}
$$
for $n$ odd, and by
$$
d(\mm_y) = 
\begin{cases}
n^2, & y \in Y^{smooth}
\\
q(0)=2 s^2, & y \in Y^{sing} \backslash \left( \bigcup_{k \in [0, s-2]} \mathcal{C}(\om + k \tau) \right)
\\
q(k+1), & y \in \mathcal{C}(\om + k \tau) \backslash \{ \underline{0} \}, ~k \in [0, s-2]
\\
1, & y = \underline{0}
\end{cases}
$$
for $n$ even.
The statement of the theorem follows by combining the formulas for $d(\mm_y)$ with \eqref{dis}. Namely for $n$ odd, the smallest value of $d(\mm_y)$ for $y \in C(\omega + k \tau) \backslash \{\underline{0}\}$, with $0 \leq k \leq n-2$, is  $q(\lfloor \frac{n}{2}\rfloor)$ and the largest value for such $d(\mm_y)$ is $q(n-1)$. Moreover, for $n$ even, the smallest value of $d(\mm_y)$ for  $y \in \mathcal{C}(\omega + k \tau) \backslash \{\underline{0}\}$, with $0 \leq k \leq s-2$, is $q(\lfloor \frac{s}{2}\rfloor)$ and the largest value for such $d(\mm_y)$ is $q(s-1)$; note that $q(0)>q(s-1)$.
\end{proof}


\bibliography{Poisson-Skly4-biblio}

\printindex
\end{document}